\documentclass[11pt]{article}
\usepackage[utf8]{inputenc}

\usepackage{epsfig,epsf,fancybox}
\usepackage{amsmath}
\usepackage{mathrsfs,bbm}
\usepackage{amssymb}
\usepackage{graphicx}
\usepackage{color}
\usepackage{multirow}
\usepackage{paralist}
\usepackage{verbatim}
\usepackage{galois}
\usepackage{algorithm}
\usepackage{algorithmic}
\usepackage{boxedminipage}
\usepackage{booktabs}
\usepackage{accents}
\usepackage{stmaryrd}
\usepackage{subfig}
\usepackage{appendix}
\usepackage{graphicx}
\usepackage{epstopdf}
\usepackage{appendix}
\usepackage{cases}  
\usepackage[top=1in,bottom=1.2in,left=1in,right=1in,xetex]{geometry}
\usepackage{extarrows}

\usepackage{titlesec}
\usepackage{titletoc}

\usepackage[unicode=true,pdfusetitle,
bookmarks=true,bookmarksnumbered=true,bookmarksopen=true,bookmarksopenlevel=3,
breaklinks=false,pdfborder={0 0 1},backref=false,colorlinks=false]
{hyperref}

\numberwithin{equation}{section}

\newtheorem{theorem}{Theorem}[section]
\newtheorem{corollary}[theorem]{Corollary}
\newtheorem{lemma}[theorem]{Lemma}

\newtheorem{remark}[theorem]{Remark}

\newtheorem{problem}{Problem}[section]

\newcommand{\lr}{\mathcal{L}}
\newcommand{\e}{\mathbb{E}}
\newcommand{\br}{\mathbb{R}}
\newcommand{\pr}{\mathcal{P}}
\newcommand{\dd}{\partial}

\newcommand{\brn}{{\mathbb{R}^n}}

\newcommand{\brd}{{\mathbb{R}^d}}
\newcommand{\de}{\Delta}

\newcommand{\hv}{\hat{V}}
\newcommand{\hvv}{\hat{v}}

\newcommand{\tv}{\tilde{V}}
\newcommand{\tx}{\tilde{X}}
\newcommand{\ty}{\tilde{Y}}

\newcommand{\tp}{\tilde{P}}
\newcommand{\tq}{\tilde{Q}}

\newcommand{\cx}{\check{X}}
\newcommand{\cy}{\check{Y}}

\newcommand{\hr}{\mathcal{H}}
\newcommand{\ur}{\mathcal{U}}
\newcommand{\V}{\mathcal{V}}

\newcommand{\bh}{\mathbb{H}}
\newcommand{\bx}{\bar{X}}

\allowdisplaybreaks[2]

\makeatother

\title{Degenerate Mean Field Type Control with Linear and Unbounded Diffusion, and their Associated Equations}

\author{Alain Bensoussan\footnote{E-mail: axb046100@utdallas.edu}\\
	International Center for Decision and Risk Analysis,\\
	Naveen Jindal School of Management, University of Texas at Dallas\\\\
	Ziyu Huang\footnote{E-mail: zyhuang19@fudan.edu.cn}\\
	School of Mathematical Sciences, Fudan University \\\\
	Shanjian Tang\footnote{E-mail: sjtang@fudan.edu.cn}\\
	Department of Finance and Control Sciences,\\ 
	School of Mathematical Sciences, Fudan University \\\\
	Sheung Chi Phillip Yam\footnote{E-mail: scpyam@sta.cuhk.edu.hk}\\
	Department of Statistics, The Chinese University of Hong Kong}
\date{}
\begin{document}

\maketitle

\begin{abstract}
	We study the well-posedness of a system of forward-backward stochastic differential equations (FBSDEs) corresponding to a degenerate mean field type control problem, when the diffusion coefficient depends on the state together with its measure and also the control. Degenerate mean field type control problems are rarely studied in the literature. Our method is based on a lifting approach which embeds the control problem and the associated FBSDEs in Wasserstein spaces into certain Hilbert spaces. We use a continuation method to establish the solvability of the FBSDEs and that of the G\^ateaux derivatives of this FBSDEs. We then explore the regularity of the value function in time and in measure argument, and we also show that it is the unique classical solution of the associated Bellman equation. We also study the higher regularity of the linear functional derivative of the value function, by then, we obtain the classical solution of the mean field type master equation.\\
	
	
	\noindent {\bf 2000 MR Subject Classification } 35R15, 49L25, 49N70, 91A13, 93E20, 60H30, 60H10.
\end{abstract}

\section{Introduction}

Mean field type (or McKean-Vlasov) control problems and mean field games have been widely studied in recent years. In such a problem, the controlled dynamical system is affected not only by the state and the control, but also by the probability distribution of the state. There are numerous works in various settings in this area. For probabilistic approaches to mean field type control problems, we refer to Buckdahn-Li-Peng-Rainer \cite{BR}, Cardaliaguet-Delarue-Lasry-Lions \cite{CDLL}, Carmona-Delarue \cite{CR,book_mfg} and Chassagneux-Crisan-Delarue \cite{CJF}. For the dynamic programming principle and Bellman equation of McKean-Vlasov control problem, we refer to Djete-Possamai-Tan \cite{DMF} and Pham-Wei \cite{PH}. For the lifting method and a Hilbert space approach for mean field type control problem, we refer to Bensoussan-Graber-Yam \cite{AB7}, Bensoussan-Huang-Yam \cite{AB6,AB8}, Bensoussan-Tai-Yam \cite{AB5} and Bensoussan-Yam \cite{AB}. For mean field games and stochastic differential games of McKean-Vlasov type, we refer to Bensoussan-Frehse-Yam \cite{AB_book}, Cardaliaguet-Cirant-Porretta \cite{CP1}, Cosso-Pham \cite{CA}, Ferreira-Gomes \cite{FR}, Gomes-Pimentel \cite{GDA1}, Gomes-Pimentel-Voskanyan \cite{GDA}, Huang-Malham\'e-Caines \cite{HM}, Lauri\`ere-Tangpi \cite{LM}, Porretta \cite{PA} and Ricciardi \cite{RM}. 

In this paper, we study the well-posedness of a degenerate mean field type control problem, and give the unique classical solution of its associated Bellman equation and mean field type master equation. Our approach is to study the regularity of a system of forward-backward stochastic differential equations (FBSDEs), which is arising from the maximum principle \cite{AB8}. This differs from the recently well-received PDE approach via master equation in the contemporary literature, which can be applied mostly to mean field game problems but not for mean field type control problems. Our approach, being more aligned with the stochastic control method, allows us to deal with a linear and unbounded diffusion, and also requires less regularity in the coefficient functions compared to the existing literature. To the best of our acknowledge, the study of degenerate mean field type control problems are rare in the contemporary literature, and the focus is essentially on the case when the diffusion is a constant. When the diffusion coefficient is degenerate, the general belief is that the difficulty in the mathematical analysis is escalated to another level. In our work, we allow the diffusion coefficient to depend on the state together with its measure and also the control. And our running cost function can be non-separable with quadratic growth.

Our method is based on a lifting approach which embeds the control problem and the associated FBSDEs in Wasserstein spaces into those in Hilbert spaces. The idea of lifting the control problem into Hilbert spaces was first introduced by P. L. Lions, and was then used in \cite{AB} to study the control problem in the Wasserstein space. Later, we used a different Hilbert space in \cite{AB6} to recover the first order mean field type control problem as a particular case. And the lifting method in this article for second order mean field type control problem is fundamentally different from the lifting method proposed by P. L. Lions. This kind of Hilbert space is already used in our previous work \cite{AB5} for the case when the diffusion is a matrix, and the cost functional is separable; and in \cite{AB8} for the maximum principle for mean field type control problem with a general diffusion coefficient. One advantage of our lifting approach is that it reconciles the closed/open loop approaches and avoid the difficulties of the traditional feedback control methodology, which does not require high regularity assumptions. 

One major concern of this paper is to study the well-posedness of the associated FBSDEs. We use a continuation method for general FBSDEs defined on Hilbert spaces under a monotonicity condition, and then apply this result for the solvability of the system of FBSDEs corresponding to our mean field type control problem, and also the solvability of its Jacobian flow. This kind of continuation method is first proposed by Hu-Peng \cite{YH2}, and then used by Carmona-Delarue \cite{CR} and Ahuja-Ren-Yang \cite{SA1} for FBSDEs for mean field games. The monotonicity condition here is satisfied when the running cost functional is strongly convex in the control variable. We also refer to  \cite{AB5,AB9} for the relation between the convexity assumption for mean field type control problem and the displacement monotonicity condition and Lasry-Lions monotonicity condition extensively used in the literature \cite{SA1,CP,CDLL,book_mfg,CJF,HZ1}.

Another major concern of this article is the regularity of the value function in time and distribution, and the solvability of the associated Bellman and master equations. For the differentiability in distribution, we first use the Hilbert method to give the G\^ateaux derivatives of the value function, and then establish its linear functional differentiability. To solve the Bellman equation, we  use  It\^o's formula for processes on Hilbert spaces (see \cite{AB5,AB} for more details). To solve the master equation, we  calculate the linear functional derivatives of the solution of the FBSDEs associated with the mean field type control problems. Under mild assumptions, we show that the value function is the unique classical solution of the Bellman equation, and the linear functional derivative of the value function is the unique classical solution of the master equation. We can see that, the well-posedness of the mean field type master equation requires higher regularity of the coefficients than the well-posedness of Bellman equation, and the latter also requires higher regularity than the solvability of the mean field type control problem. Therefore, we use the optimal control method and the FBSDE solution to study the Bellman equation and the master equation, instead of studying the master equation. We also refer to Bensoussan-Tai-Yam \cite{AB5}, Cardaliaguet-Delarue-Lasry-Lions \cite{CDLL}, Gangbo-M\'esz\'aros-Mou-Zhang \cite{GW} and Huang-Tang \cite{HZ} for more discussion on the mean field type master equation.

This article is organized as follows. In Section~\ref{sec:2}, we introduce the Wasserstein spaces, Hilbert spaces and various derivatives of functionals. In Section~\ref{sec:3}, we give the formulation of our mean field type control problem, the FBSDEs and value function, and introduce our own lifting method but different from Lions' one. Section~\ref{sec:FB_Hilbert} studies the solvability and regularity of FBSDEs taking values in Hilbert spaces. In Section~\ref{sec:value_Hilbert}, we study the differentiability of the value function in time and the G\^ateaux differentiability in distribution. In Section~\ref{sec:MFC}, we go back to mean field type control problems and interpret our results for FBSDEs in Hilbert spaces as those for mean field type control problems. Section~\ref{sec:Bellman} gives the linear functional derivatives of the value function in distribution, and also gives the solvability of the Bellman equation. In Section~\ref{sec:Master}, we study the higher regularity of the linear functional derivative of the value function, and hence provide the solvability of the master equation of mean field type. The proof of statements in Sections~\ref{sec:FB_Hilbert}, \ref{sec:value_Hilbert}, \ref{sec:Bellman} and \ref{sec:Master} are shown in Appendices~\ref{app:01}, \ref{app:02}, \ref{app:1} and \ref{app:2}, respectively. 

\section{Problem Setting and Notations}\label{sec:2}

\subsection{Wasserstein space and Hilbert subspaces of $L^2$-random variables}

We denote by $\mathcal{P}_{2}(\brn)$ the space of probability measures on $\brn$, with the second moment, i.e. $\int_{\brn}|x|^{2}dm(x)<+\infty$, equipped with the 2-Wasserstein metric: 
\begin{align*}
	W_2(m,\tilde{m}):=\inf_{\pi\in\Gamma(m,\tilde{m})}\sqrt{\int_{\brn\times\brn}|x-y|^2\pi(dx,dy)},
\end{align*}
where $\Gamma(m,\tilde{m})$ is the set of joint probability measures with respective marginals $m$ and $\tilde{m}$. Recall that a $\{m_{k},\ k\geq 1\}$ converges to $m$ in $\mathcal{P}_{2}(\brn)$ if and only if it converges in the sense of the weak convergence and $W_2(m_k,\delta_0)\to W_2(m,\delta_0)$ as $k\to+\infty$; see Ambrosio-Gigli-Savar\'e \cite{AL} and Villani \cite{VC} for details.

Let $(\Omega,\mathcal{A},\mathbb{P})$ denote an atomless probability space. In this paper, the control is $\brd$-valued, and the state is $\brn$-valued.  For $m\in\pr_2(\brn)$, consider $L^2_m (\brn;\brn)$ the Hilbert space with respect to the following inner product
\begin{align*}
	(X,X'):=\int_{\brn} X_x^* X'_xdm(x),\quad X,X'\in L^2_m (\brn;\brn), 
\end{align*}
where $X^*$ is the transpose of $X$. We also work in $\hr_m:=L^2(\Omega,\mathcal{A},\mathbb{P};L_m^2(\brn;\brn))$ and $\ur_m:=L^2(\Omega,\mathcal{A},\mathbb{P};L_m^2(\brn;\brd))$ for $m\in\pr_2(\brn)$. Elements in $\hr_m$ and $\ur_m$ are generally denoted by $X_x(\omega)$ and $V_x(\omega),\ (x,\omega)\in\brn\times\Omega$; $X_x$ and $V_x$ are random vectors in $L^2(\Omega,\mathcal{A},\mathbb{P};\brn)$ and $L^2(\Omega,\mathcal{A},\mathbb{P};\brd)$, respectively. The inner product over $\hr_m$ is defined as
\begin{align*}
	&(X,X')_{\hr_m}:=\e\left[\int_\brn X_x^* X'_x dm(x)\right],\quad X,X'\in\hr_m.
\end{align*}
The inner product in $\ur_m$ is defined  in a similar way. 

Let $w(t),\ t\geq 0$ be an $n$-dimensional standard Wiener process on $(\Omega,\mathcal{A},\mathbb{P})$. For any $0\le t\le s\le T$, we denote by $\mathcal{W}_t^s:=\sigma(w(\tau)-w(t);\tau\in[t,s])$ the filtration generated by the Wiener process, and abbreviate $\mathcal{W}_t:=\mathcal{W}_t^T$; so $\mathbb{P}$ is referred to the randomness generated purely by Brownian increment. For any $X\in\hr_m$ independent of $\mathcal{W}_t$, we denote by $\mathcal{W}^s_{tX}:=\sigma(X)\bigvee \mathcal{W}_t^s$ the filtration generated by $X$ and the Wiener process, and denote by $\mathcal{W}_{tX}$ the filtration generated by the $\sigma$-algebras $\mathcal{W}_{tX}^s$ for $s\in[t,T]$. Denote by $L^2_{\mathcal{W}_{tX}(t,T;\ur_m)}$ the subspace of $L^2(t,T;\ur_m)$ of all processes adapted to the filtration $\mathcal{W}_{tX}$. 

For $X\in\hr_m$, we denote by $X\#(m\otimes\mathbb{P})$ the push forward of $m\otimes\mathbb{P}$ by the map $(x,\omega)\mapsto X_x(\omega)$. For test function $\varphi$, 
\begin{align*}
	\int_\brn \varphi(y)d(X\#(m\otimes\mathbb{P}))(y)=\e \left[\int_\brn \varphi(X_x)dm(x)\right].
\end{align*}
The measure $X\#(m\otimes\mathbb{P})\in\pr_2(\brn)$ and satisfies
\begin{equation}\label{m2}
	\begin{split}
		&W_2(X\#(m\otimes\mathbb{P}),X'\#(m\otimes\mathbb{P}))\le \|X-X'\|_{\hr_m},\quad X,X'\in \hr_m.
	\end{split}	
\end{equation}
We refer to \cite{AB6,AB8,AB5} for details. For notational convenience, we use the notation $X\otimes m$ instead of $X\#(m\otimes\mathbb{P})$. 

\subsection{Derivatives of functionals}

Our work is based on the derivative in $m\in\pr_2(\brn)$. We use the concept of the linear functional derivative; see Carmona-Delarue \cite{book_mfg}. For a functional $k:\mathcal{P}_{2}(\brn)\to\br$, the linear functional derivative of $k(m)$ at $m$ is a function $\mathcal{P}_{2}(\brn)\times \brn\ni(m,x)\mapsto\dfrac{dk}{d\nu}(m)(x)$ being continuous under the product topology, satisfying $\int_{\brn}\left|\dfrac{dk}{d\nu}(m)(x)\right|^{2}dm(x)\leq c(m)$ for some positive constant $c(m)$ depending locally on $m$, and for any $m'\in\mathcal{P}_{2}(\brn)$,
\begin{align*}
	\lim_{\epsilon\to0}\dfrac{k(m+\epsilon(m'-m))-k(m)}{\epsilon}=\int_\brn\dfrac{dk}{d\nu}(m)(x)(dm'(x)-dm(x)).
\end{align*}
We have not expressed $\dfrac{dk}{dm}(m)(x)$ to make the difference between the notation $\nu$ and the argument $m$. In \cite{AB8}, we have built a connection between the linear functional derivative in $\pr_2(\brn)$ and the G\^ateaux derivative in $\hr_m$. For a differentiable functional $k$ such that the derivative $D_\xi\frac{d f}{d\nu}(\mu)(\xi)$ is continuous in $(\mu,\xi)$ and $\left|D_\xi\frac{d k}{d\nu}(\mu)(\xi)\right|\le c(\mu)(1+|\xi|)$, $(\mu,\xi)\in\pr_2(\brn)\times\brn$, we define $K:\hr_m\to\br$ as $K(X):=k(X\otimes m)$, $X\in \hr_m$. Then, $K$ is G\^ateaux differentiable, and
\begin{align}\label{lem01_1}
	D_X K(X)\Big|_x=D_\xi\frac{d k}{d\nu}(X\otimes m)(X_x),\quad X\in \hr_m.
\end{align}
When $X$ is the identity function $I$, \eqref{lem01_1} is identical to the $L$-derivative $\dd_m k(m)(x)$ of Carmona-Delarue \cite{book_mfg}. We refer to Bensoussan-Graber-Yam \cite{AB7} and Bensoussan-Tai-Yam \cite{AB5} for further discussion.

\section{Problem Formulation}\label{sec:3}

\subsection{Mean field type control problem}

We consider the following coefficient maps
\begin{align*}
	&\text{Drift:}\quad\qquad\qquad  f:\brn\times\mathcal{P}_{2}(\brn)\times \brd\times[0,T]\rightarrow \brn;\\
	&\text{Diffusion:}\quad\qquad \sigma=(\sigma^1,\dots,\sigma^N): \brn\times\mathcal{P}_{2}(\brn)\times \brd\times [0,T]\rightarrow \br^{n\times n};\\
	&\text{Running cost:}\quad\  g:\brn\times\mathcal{P}_{2}(\brn)\times \brd\times [0,T]\rightarrow \br;\\
	&\text{Terminal cost:}\quad\  g_T:\brn\times\mathcal{P}_{2}(\brn)\rightarrow \br.
\end{align*}
For initial data $(t,m)\in [0,T]\times\pr_2(\br^n)$, and $X\in \hr_m$, which is independent of $\mathcal{W}_t$, we consider the following mean field type control problem:
\begin{problem}\label{intr_4}
	$\inf_{v_{X_\cdot}\in L^2_{\mathcal{W}_{tX}}(t,T;\ur_m)} \ J_{Xmt}(v_{X_\cdot})$, 
	\begin{equation*}
		\begin{aligned}
			&\text{where}\quad J_{Xmt}(v_{X_\cdot}):=\int_{t}^{T}\e\left[\int_{\brn}g(X^v_{X_x}(s),X_{X}^v(s)\otimes m),v_{X_x}(s),s)dm(x)\right]ds\\
			&\qquad\qquad\qquad\qquad\qquad+\e\left[\int_{\brn}g_T(X^v_{X_x}(T),X_{X}^v(T)\otimes m)dm(x)\right],\\
			&\text{and}\quad X_{X_x}^v(s):=X_x+\int_t^s f(X_{X_x}^v(r),X_{X}^v(r)\otimes m,v_{X_x}(r),r)dr\\
			&\ \ \qquad\qquad\qquad\qquad +\sum_{j=1}^{N} \int_t^s \sigma^j(X_{X_x}^v(r),X_{X}^v(r)\otimes m,v_{X_x}(r),r)dw_j(r).
		\end{aligned}
	\end{equation*}
\end{problem}
The controlled state process $X_\cdot^v(\cdot)$ belongs to $L^2_{\mathcal{W}_{tX}}(t,T;\hr_m)$ under Assumptions (B)'s on coefficients $f$ and $\sigma$ to be stated in Section~\ref{sec:MFC}. For $s\in[t,T]$, $X^v(s)\otimes m$ is the probability law of $X^v_\eta(s)$ when $\eta$ is equipped with the probability $m\otimes\mathbb{P}$. We define the Lagrangian $L:\brn\times\pr_2(\brn)\times\brd\times[0,T]\times\brn\times\br^{n\times n}\to\br$ as
\begin{align}\label{H'}
	L(x,m,v,s;p,q):=p^* f(x,v,m,s)+\sum_{j=1}^n (q^{j})^* \sigma^j(x,v,m,s)+g(x,v,m,s),
\end{align}
while the minimizer $\hat{v}(x,m,s;p,q)$ (under suitable assumptions on the coefficient functions to be stated in Section~\ref{sec:MFC}) satisfies:
\begin{align}\label{hv'}
	D_v L (x,m,\hat{v}(x,m,s;p,q),s;p,q)=0,
\end{align}
and the Hamiltonian is simply:
\begin{equation}\label{bh'}
	\begin{split}
		&H (x,m,s;p,q):=L(x,m,\hat{v}(x,m,s;p,q),s;p,q).
	\end{split}
\end{equation}

\subsection{Lifting method and maximum principle}\label{subsec:MP}

In our previous work \cite{AB8}, we use a lifting method to embed the mean field type control problem into a Hilbert space. For a fixed $m\in\pr_2(\brn)$, we define maps $F,A^j:\hr_m\times \ur_m\times[0,T]\to \hr_m$, functionals $G:\hr_m\times \ur_m\times[0,T]\to \br$ and $G_T:\hr_m\to \br$ as, for $(X,V,s)\in \hr_m\times \ur_m\times[0,T]$, $x\in\brn$ and $1\le j\le n$,
\begin{equation}\label{GAF}
	\begin{split}
		&F(X,V,s)\Big|_x:=f(X_x,X\otimes m,V_x,s);\quad A^j(X,V,s)\Big|_x:=\sigma^j(X_x,X\otimes m,V_x,s);\\
		&G(X,V,s):=\e\left[\int_\brn g(X_x,X\otimes m,V_x,s)dm(x)\right];\\
		&G_T(X):=\e\left[\int_\brn g_T(X_x,X\otimes m)dm(x)\right].
	\end{split}
\end{equation}
This formulation inspires us to study the following control problems on Hilbert spaces $\hr_m$ and $\ur_m$: for any initial time $t\in[0,T]$, and $X\in\hr_m$ which is independent of $\mathcal{W}_t$:
\begin{problem}\label{intr_5}
	$\inf_{V\in L^2_{\mathcal{W}_{tX}}(t,T;\ur_m)} \ J_{Xmt}(V)$,
	\begin{equation*}
		\begin{aligned}
			&\text{where} \quad J_{Xmt}(V):=\int_t^T G\left(X^{V}(s),V(s),s\right) ds+G_T\left(X^{V}(T)\right),\\
			&\text{and}\quad X^{V}(s):=X+\int_t^s F\left(X^{V}(r),V(r),r\right)dr+\sum_{j=1}^{n} A^j\left(X^{V}(r),V(r),r\right)dw_j(r).
		\end{aligned}
	\end{equation*}
\end{problem}
We define $\lr:\hr_m\times \ur_m\times[0,T]\times \hr_m\times \hr_m^n\to\br$ and $\hv:\hr_m \times[0,T]\times\hr_m\times\hr_m^n\to \ur_m$ as
\begin{align}
	&\text{\small (Lift Lagrangian)}\quad \lr(X,V,s;P,Q):=\e\left[\int_\brn {L}(X_x,X\otimes m,V_x,s;P_x,Q_x)dm(x)\right],\label{H''}\\
	&\text{\small (Lift control)}\quad \left.\hv(X,s;P,Q)\right|_x:=\hat{v}(X_x,X\otimes m,s;P_x,Q_x),\quad x\in\brn.\label{v}
\end{align}
We know from \eqref{H''} and \eqref{v} that 
\begin{align}
	&\lr(X,V,s;P,Q):=\left(P,F(X,V,s)\right)_{\hr_m}+\sum_{j=1}^n \left(Q^j,A^j(X,V,s)\right)_{\hr_m} +G(X,V,s),\label{H}\\
	&\text{\small (Lift optimality condition)}\quad D_V \lr (X,s,\hv(X,s;P,Q),s;P,Q)=0.\label{hv}
\end{align}
From Theorems 1 and 2 in \cite{AB8}, we know that Problem~\ref{intr_5} has a unique optimal control under Assumptions (A)'s to be stated in Section~\ref{sec:FB_Hilbert}, and the optimal control is
\begin{align}\label{max_Hilbert}
	\hv(Y_{Xt}(s),s;P_{Xt}(s),Q_{Xt}(s)),\quad \text{a.e.}\ s\in[t,T],
\end{align}
where $Y_{Xt}$ is the corresponding controlled state process and $(P_{Xt},Q_{Xt})$ is the corresponding adjoint process defined by
\begin{equation*}\label{p}
	\begin{aligned}
		P_{Xt}(s)=&D_X G_T\left(Y_{Xt}(T)\right)+\int_s^T D_X\lr \left(Y_{Xt}(r),\hat{V}_{Xt}(s),r;P_{Xt}(r),Q_{Xt}(r)\right)dr\\
		&-\sum_{j=1}^n\int_s^T P_{Xt}^j(r)dw_j(r),\quad s\in[t,T].
	\end{aligned}
\end{equation*}
In the rest of this article, without much ambiguity, we simply suppress the arguments of $\hv$ of \eqref{max_Hilbert}, and denote it as $\hat{V}_{Xt}(s)$. In Theorems 3 and 4 of \cite{AB8}, we apply our results for control problems in Hilbert spaces back to the mean field type control problem, by then we deduce that Problem~\ref{intr_4} has a unique optimal control such that, $\text{a.e.}\ s\in[t,T],\  \text{a.s.}\ dm(x), \ \text{a.s.}\ d\mathbb{P}(\omega)$,
\begin{equation}\label{max_MFC}
	\begin{split}
		\hat{v}_{X_xmt}(s):=\hat{v}(Y_{X_xmt}(s),Y_{Xmt}(s)\otimes m,s;P_{X_xmt}(s),Q_{X_xmt}(s)),
	\end{split}
\end{equation}
where $Y_{X_\cdot t}(\cdot)$ is the corresponding controlled state process and $(P_{X_\cdot t}(\cdot),Q_{X_\cdot t}(\cdot))$ is the corresponding adjoint process defined as
\begin{equation*}\label{p'}
	\begin{split}
		&P_{X_xmt}(s)=D_xg_T(Y_{X_xmt}(T),Y_{Xmt}(T)\otimes m)\\
		&\quad\qquad\qquad +\bar{\e}\left[\int_\brn D_\xi\frac{d g_T}{d\nu}\left(\bar{Y}_{\bx_ymt}(T),Y_{Xmt}(T)\otimes m\right)(Y_{X_xmt}(T)) dm(y)\right]\\
		&\quad\qquad\qquad +\int_s^T \bigg[D_x {L}(Y_{X_xmt}(r),Y_{Xmt}(r)\otimes m,\hat{v}_{X_xmt}(r),r;P_{X_xmt}(r),Q_{X_xmt}(r))\\
		&\quad\qquad\qquad\qquad\qquad +\bar{\e}\bigg(\int_\brn D_\xi \frac{d L}{d\nu}\Big(\bar{Y}_{\bx_ymt}(r),Y_{Xmt}(r)\otimes m,\bar{\hat{v}}_{\bx_ymt}(r),r;\\
		&\quad\qquad\qquad\qquad\qquad\qquad\qquad\qquad\qquad \bar{P}_{\bx_ymt}(r),\bar{Q}_{\bx_ymt}(r)\Big)\ (X_{X_xmt}(r)) dm(y)\bigg)\bigg]dr\\
		&\quad\qquad\qquad -\sum_{j=1}^n \int_s^T Q^j_{X_xmt}(r)dw_j(r),\quad s\in[t,T].
	\end{split}
\end{equation*}
Here, the process $\left(\bar{Y}_{\bx mt}(s),\bar{\hat{v}}_{\bx mt}(s),\bar{P}_{\bx mt}(s),\bar{Q}_{\bx mt}(s)\right)$ is an independent copy of the process $(Y_{Xmt}(s),\hat{v}_{Xmt}(s),P_{Xmt}(s),Q_{Xmt}(s))$; from now on, for any random variable $\xi$, we use $\bar{\xi}$ to stand for its independent copy, and use $\bar{\e}[\bar{\xi}]$ for its expectation. 


\subsection{FBSDEs and value function}

From the maximum principle \eqref{max_Hilbert}, we know that Problem~\ref{intr_5} is associated with the following system of FBSDEs: for $s\in[t,T]$,
\begin{equation}\label{FB:1}
	\left\{
	\begin{aligned}
		&Y_{Xt}(s)=X+\int_t^s D_P\bh (Y_{Xt}(r),r;P_{Xt}(r),Q_{Xt}(r))dr\\
		&\qquad\qquad +\sum_{j=1}^n \int_t^s  D_{Q^j}\bh(Y_{Xt}(r),r;P_{Xt}(r),Q_{Xt}(r))dw_j(r),\\
		&P_{Xt}(s)=D_X G_T(Y_{Xt}(T))+\int_s^T D_X \bh(Y_{Xt}(r),r;P_{Xt}(r),Q_{Xt}(r))dr\\
		&\qquad\qquad -\sum_{j=1}^n\int_s^T Q_{Xt}^j(r)dw_j(r),
	\end{aligned}
	\right.
\end{equation}
where $\bh:\hr_m\times[0,T]\times \hr_m\times\hr_m^n\to \br$ is defined as
\begin{align}\label{bh''}
	\bh(X,s;P,Q):=\e\left[\int_\brn H(X_x,X\otimes m,s;P_x,Q_x)dm(x)\right],
\end{align}
and satisfies $\bh(X,s;P,Q)=\lr(X,\hv(X,s;P,Q),s;P,Q)$. We shall study the G\^ateaux differentiability of FBSDEs \eqref{FB:1} with respect to initial $(X,t)$. We then get back to Problem~\ref{intr_4}, and study the well-posedness of the value function $\mathcal{V}:\pr_2(\brn)\times[0,T]\to\br$ defined as
\begin{align}\label{def:V}
	\mathcal{V}(X\otimes m,t):=\inf_{V\in L_{\mathcal{W}_{tX}}^2(t,T;\ur_m)} J_{Xmt}(V).
\end{align}
In particular, we shall discuss its G\^ateaux differentiability and linear functional differentiability in $m$ via the linear functional derivative of the following FBSDEs: for $s\in[t,T]$,
\begin{equation}\label{FB:4}
	\left\{
	\begin{aligned}
		&Y_{X_xmt}(s)=X_x+\int_t^s D_p H(Y_{X_xmt}(r),Y_{Xmt}(r)\otimes m,r;P_{X_xmt}(r),Q_{X_xmt}(r))dr\\
		&\quad\qquad\qquad +\sum_{j=1}^n\int_t^s  D_{q^j}H(Y_{X_xmt}(r),Y_{Xmt}(r)\otimes m,r;P_{X_xmt}(r),Q_{X_xmt}(r))dw_j(r),\\
		&P_{X_xmt}(s)=D_xg_T(Y_{X_xmt}(T),Y_{Xmt}(T)\otimes m)\\
		&\quad\qquad\qquad +\bar{\e}\left[\int_\brn D_\xi\frac{d g_T}{d\nu}\left(\bar{Y}_{\bx_ymt}(T),Y_{Xmt}(T)\otimes m\right)(Y_{X_xmt}(T)) dm(y)\right]\\
		&\quad\qquad\qquad +\int_s^T\bigg[D_x H(Y_{X_xmt}(r),Y_{Xmt}(r)\otimes m,r;P_{X_xmt}(r),Q_{X_xmt}(r))\\
		&\qquad\qquad\qquad\qquad+\bar{\e}\bigg(\int_\brn D_\xi \frac{d H}{d\nu}\Big(\bar{Y}_{\bx_ymt}(r),Y_{Xmt}(r)\otimes m,r;\bar{P}_{\bx_ymt}(r),\bar{Q}_{\bx_ymt}(r)\Big)\\
		&\qquad\qquad\qquad\qquad\qquad\qquad\qquad\qquad (Y_{X_xmt}(r)) dm(y)\bigg)\bigg]dr,\\
		&\quad\qquad\qquad -\sum_{j=1}^n \int_s^T Q^j_{X_xmt}(r)dw_j(r),
	\end{aligned}
	\right.
\end{equation}
where $\left(\bar{Y}_{\bx mt}(s),\bar{P}_{\bx mt}(s),\bar{Q}_{\bx mt}(s)\right)$ is an independent copy of $(Y_{Xmt}(s),P_{Xmt}(s),Q_{Xmt}(s))$. With these results, we shall show that $\mathcal{V}$ is the unique classical solution of the Bellman equation \eqref{Bellman'}; while the linear functional derivative of $\V$ is the unique classical solution of the mean field type master equation \eqref{master}.

\section{Well-posedness of FBSDEs \eqref{FB:1} and Regularity of the Solution}\label{sec:FB_Hilbert}

In this section, we study the regularity of the FBSDEs \eqref{FB:1} defined on Hilbert space with respect to initial data $(t,X)\in[0,T]\times \hr_m$. We state our assumptions for the maps $F,A_j,G$ and $G_T$ on Hilbert spaces. For the sake of convenience, we use the same constant $L>0$ for all conditions below. \\
\textbf{(A1)} The maps $F$ and $A^j$ (for $1\le j\le n$) are linear in $X$ and $V$. That is,
\begin{align*}
	&F(X,V,s)=F_0(s)+\mathcal{F}_1(s)X+\mathcal{F}_2(s)V,\\
	&A^j(X,V,s)=A^j_0(s)+\mathcal{A}^j_1(s)X+\mathcal{A}^j_2(s)V,
\end{align*}
where $\mathcal{F}_1(s)$ and $\mathcal{A}^j_1(s)$ are linear maps on $\hr_m$, $\mathcal{F}_2(s)$ and $\mathcal{A}^j_2(s)$ are linear maps on $\ur_m$; $\|F_0(s)\|_{\hr_m}$, $\|\mathcal{F}_1(s)\|_{L(\hr_m;\hr_m)}$,  $\|\mathcal{F}_2(s)\|_{L(\ur_m;\hr_m)}$,  $\|A^j_0(s)\|_{\hr_m}$, $\|\mathcal{A}^j_1(s)\|_{L(\hr_m;\hr_m)}$ and $\|\mathcal{A}^j_2(s)\|_{L(\ur_m;\hr_m)}$ are all bounded above by $L$.\\
\textbf{(A2)} The running cost functional $G$ satisfies
\begin{align*}
	|G(X,V,s)|\le L(1+\|X\|^2_{\hr_m}+\|V\|^2_{\ur_m}),\quad (X,V,s)\in \hr_m\times \ur_m\times[0,T].
\end{align*}
For any $s\in [0,T]$, the functional $\hr_m\times \ur_m\ni(X,V)\mapsto G(X,V,s)\in\br$ is differentiable, with the derivatives satisfying for $s\in[0,T]$, $X,X',\tx\in\hr_m$ and $V,V',\tx\in\ur_m$,
\begin{align*}
	&\|D_X G (X,V,s)\|_{\hr_m}+\|D_V G (X,V,s)\|_{\ur_m} \le L\left(1+\|X\|_{\hr_m}+\|V\|_{\ur_m}\right),\\
	&\left\|D_X G \left(X',V',s\right)-D_X G (X,V,s)\right\|_{\hr_m}+\left\|D_V G \left(X',V',s\right)-D_V G (X,V,s)\right\|_{\ur_m} \\
	&\ \ \quad\qquad\qquad\qquad\qquad\qquad\qquad\qquad\qquad \le L\left(\left\|X'-X\right\|_{\hr_m}+\left\|V'-V\right\|_{\ur_m}\right),
\end{align*}
and the maps $(X,V)\mapsto D_X G (X,V,s), D_V G (X,V,s)$ are G\^ateaux differentiable, with the continuous derivatives satisfying
\begin{align*}
	\left\|D^2_X G \left(X,V,s\right)\left(\tx\right)\right\|_{\hr_m}+\left\|D_XD_V G (X,V,s)\left(\tx\right)\right\|_{\ur_m}\le &\  L\left\|\tx\right\|,\\
	\left\|D^2_V G (X,V,s)\left(\tv\right)\right\|_{\ur_m}+\left\|D_VD_X G (X,V,s)\left(\tv\right)\right\|_{\hr_m}\le &\  L\left\|\tv\right\|.
\end{align*}
The terminal cost functional $G_T$ satisfies $|G_T(X)|\le L(1+\|X\|^2_{\hr_m}),\ X\in \hr_m$. The functional $\hr_m\ni X\mapsto G_T(X)\in\br$ is differentiable, with the derivative satisfying for $X,X',\tx\in\hr_m$, 
\begin{align*}
	\|D_X G_T (X)\|_{\hr_m} \le L\left(1+\|X\|_{\hr_m}\right),\ \left\|D_X G \left(X'\right)-D_X G (X)\right\|_{\hr_m}\le L \left\|X'-X\right\|_{\hr_m},
\end{align*}
and the map $X\mapsto D_X G_T (X)$ is G\^ateaux differentiable with a continuous derivative satisfying 
\begin{align*}
	\left\|D^2_X G (X,V,s)\left(\tx\right)\right\|_{\hr_m}\le L\left\|\tx\right\|.
\end{align*}
\textbf{(A3)} (i) There exists $\lambda>0$ such that, for any $X\in \hr_m$, $V,V'\in \ur_m$ and $s\in[0,T]$, 
\begin{align*}
	G(X,V',s)-G(X,V,s) \geq  \left(D_V G(X,V,s),V'-V\right)_{\ur_m}+\lambda\left\|V'-V\right\|^2_{\ur_m}.
\end{align*}		
(ii) There exists $\lambda> 0$ such that, for any $X,X'\in \hr_m$, $V,V'\in \ur_m$ and $s\in[0,T]$, 
\begin{align*}
	&G(X',V',s)-G(X,V,s) \geq  \left(D_X G(X,V,s),X'-X\right)_{\hr_m}+\left(D_V G(X,V,s),V'-V\right)_{\ur_m}+\lambda\left\|V'-V\right\|^2_{\ur_m},\\
	&\left(D_X G_T(X')-D_X G_T(X),X'-X\right)_{\hr_m}\geq 0.
\end{align*}
We refer to Remark 3.1 of \cite{AB5} and Section 5.3 of \cite{AB9} and Graber-M\'{e}sz\'{a}ros \cite{GM} and Remark 1 of our previous work \cite{AB8} for further discussion on the relation between the aforementioned convexity assumption (A3) and the displacement monotonicity condition and Lasry-Lions monotonicity condition for mean field type control problem and mean field game.

\subsection{Solvability, $L^2$-boundedness and continuity for FBSDEs \eqref{FB:1}}

We begin by giving the solvability of general FBSDEs defiend on Hilber spaces under the following monotonicity condition \eqref{monotonicity}, which will be used for not only the solvability of FBSDEs \eqref{FB:1}, but also that of the Jacobian flow of FBSDEs \eqref{FB:1} in Subsection~\ref{subsec:Derivatives}. For initial $(t,X)\in[0,T]\times \hr_m$, we consider the following FBSDEs: for $s\in[t,T]$,
\begin{equation}\label{FB:11}
	\left\{
	\begin{aligned}
		&Y_{Xt}(s)=X+\int_t^s \mathbb{F} (Y_{Xt}(r),r;P_{Xt}(r),Q_{Xt}(r))dr\\
		&\qquad\qquad +\sum_{j=1}^n\int_t^s \mathbb{A}^j(Y_{Xt}(r),r;P_{Xt}(r),Q_{Xt}(r))dw_j(r),\\
		&P_{Xt}(s)=\mathbb{G}_T(Y_{Xt}(T))-\int_s^T \mathbb{G}(Y_{Xt}(r),r;P_{Xt}(r),Q_{Xt}(r))dr-\sum_{j=1}^n \int_s^T Q_{Xt}^j(r)dw_j(r),
	\end{aligned}
	\right.
\end{equation}
where $\mathbb{F},\mathbb{A}^j,\mathbb{G}:\hr_m\times[0,T]\times\hr_m\times\hr_m^n\to \hr_m$ and $\mathbb{G}_T:\hr_m\to \hr_m$. We study the solvability of FBSDEs \eqref{FB:11} under the following monotonicity condition: suppose that there exists a map $\beta:\hr_m\times[0,T]\times\hr_m\times\hr_m^n\ni(X,s;P,Q)\mapsto\beta(X,s;P,Q)\in\ur_m$ and constants $\Lambda>0$, $\alpha\geq 0$, such that for any $X,X',P,P',Q^j,Q'^{j}\in\hr_m$,
\begin{equation}\label{monotonicity}
	\begin{split}
		&\left(\mathbb{G}(X',s;P',Q')-\mathbb{G}(X,s;P,Q),X'-X\right)_{\hr_m}\\
		&+\left(\mathbb{F}(X',s;P',Q')-\mathbb{F}(X,s;P,Q),P'-P\right)_{\hr_m}\\
		& +\sum_{j=1}^n \left(\mathbb{A}^j (X',s;P',Q')-\mathbb{A}^j (X,s;P,Q),Q'^{j}-Q^j\right)_{\hr_m}\\
		\le & -\Lambda \left\| \beta(X',s;P',Q')-\beta(X,s;P,Q)\right\|^2_{\ur_m}\\
		&+\alpha\left(\|X'-X\|^2_{\hr_m}+\|P'-P\|^2_{\hr_m}+\sum_{j=1}^n\|Q'^{j}-Q^j\|^2_{\hr_m}\right).
	\end{split}
\end{equation}
Condition \eqref{monotonicity} is satisfied by the coefficients of FBSDEs \eqref{FB:1} under Assumptions (A1)-(A3); also explained further in Theorem~\ref{thm2}. We also need the following Lipschitz-continuity: there exists a constant $\mathbb{L}>0$, such that,
\small
\begin{align}
	&\left\|\mathbb{F}(X',s;P',Q')-\mathbb{F}(X,s;P,Q)\right\|_{\hr_m} +\sum_{j=1}^n\left\|\mathbb{A}^j(X',s;P',Q')-\mathbb{A}^j(X,s;P,Q)\right\|_{\hr_m} \notag\\
	\le & \ \mathbb{L} \Big(\left\|X'-X\right\|_{\hr_m}+\left\| \beta(X',s;P',Q')-\beta(X,s;P,Q)\right\|_{\ur_m}\Big);\label{bh_2}\\
	&\left\|\mathbb{G}(X',s;P',Q')-\mathbb{G}(X,s;P,Q)\right\|_{\hr_m}+\left\|\mathbb{G}_T(X')-\mathbb{G}_T(X)\right\|_{\hr_m} \notag\\
	\le & \ \mathbb{L} \Big(\left\|X'-X\right\|_{\hr_m}+\left\|P'-P\right\|_{\hr_m}+\sum_{j=1}^n\left\|Q'^{j}-Q^j\right\|_{\hr_m} +\left\| \beta(X',s;P',Q')-\beta(X,s;P,Q)\right\|_{\ur_m}\Big).\label{bh_1}
\end{align}
\normalsize
We then have the following solvability of FBSDEs \eqref{FB:11}, whose proof is given in Appendix~\ref{pf:thm:1}.

\begin{lemma}\label{thm:1}
	Under the monotonicity condition \eqref{monotonicity} and the  Lipschitz-continuity conditions \eqref{bh_2}-\eqref{bh_1}, there exists a constant $c(\alpha,\mathbb{L},T)$ depending only on $(\alpha,\mathbb{L},T)$, such that when $\Lambda\geq c(\alpha,\mathbb{L},T)$, there is a unique adapted solution $(Y_{Xt},P_{Xt},Q_{Xt})$ of the FBSDEs \eqref{FB:11}. And for $X,X'\in\hr_m$, both are independent of $\mathcal{W}_t$, we have 
	\begin{align}
		&\sup_{t\le s\le T}\|Y_{Xt}(s)\|_{\hr_m}+\sup_{t\le s\le T}\|P_{Xt}(s)\|_{\hr_m}+\sum_{j=1}^n\|Q^j_{Xt}\|_{L^2(t,T;\hr_m)} \label{thm1_6}\\
		\le& \ C(\alpha,\mathbb{L},T,\Lambda)(1+\|X\|_{\hr_m}); \notag\\
		&\sup_{t\le s\le T}\left\|Y_{X't}(s)-Y_{Xt}(s)\right\|_{\hr_m}+\sup_{t\le s\le T}\left\|P_{X't}(s)-P_{Xt}(s)\right\|_{\hr_m}+\sum_{j=1}^n \left\|Q^j_{X't}-Q^j_{Xt}\right\|_{L^2(t,T;\hr_m)}\notag\\
		\le& \ C(\alpha,\mathbb{L},T,\Lambda)\|X'- X\|_{\hr_m}, \label{thm1_1}
	\end{align}
	where $C(\alpha,\mathbb{L},T,\Lambda)$ is a constant depending only on $(\alpha,\mathbb{L},T,\Lambda)$. Furthermore, if the parameter (in the monotonicity condition \eqref{monotonicity}) $\alpha=0$, and $\mathbb{G}_T$ satisfies the following monotonicity condition
	\begin{align}\label{thm1_7}
		&\left(\mathbb{G}_T(X')-\mathbb{G}_T(X),X'-X\right)_{\hr_m}\geq 0,\quad X,X'\in\hr_m,
	\end{align} 
	then, for any $\Lambda> 0$, FBSDEs \eqref{FB:11} has a unique adapted solution satisfying \eqref{thm1_6} and \eqref{thm1_1}.
\end{lemma}
From Assumptions (A1)-(A3), the map $\hv$ in \eqref{hv} is well-defined and satisfies
\begin{equation}\label{hv2}
	\begin{split}
		&\mathcal{F}^*_2(s)P+\sum_{j=1}^n {\mathcal{A}^j}^*_2(s) Q^j+D_V G\left(X,\hv(X,s;P,Q),s\right)=0.
	\end{split}
\end{equation}
Moreover, the regularity for $\hv$, due to Assumptions (A1)-(A3), is (also see \cite{AB6,AB5,AB9}):
\begin{equation}\label{hv1}
	\begin{split}
		&\left\|\hv(0,s;0,0)\right\|_{\ur_m}\le \frac{L}{2\lambda},\\
		&\left\|\hv(X',s;P',Q')-\hv(X,s;P,Q)\right\|_{\ur_m}\\
		\le& \frac{C(L)}{\lambda} \left(\left\|X'-X\right\|_{\hr_m}+\left\|P'-P\right\|_{\hr_m}+\sum_{j=1}^n\left\|Q'^{j}-Q^j\right\|_{\hr_m}\right).
	\end{split}	
\end{equation}
We then give the solvability of FBSDEs \eqref{FB:1}, and also the $L^2$-boundedness and continuity of its solution in the initial data $X\in\hr_m$. Its proof is given in Appendix~\ref{pf:thm2}.

\begin{theorem}\label{thm2}
	Under Assumptions (A1), (A2) and (A3)-(i), there exists a constant $c(L,T)$ depending only on $(L,T)$, such that when $\lambda\geq c(L,T)$, there is a unique adapted solution $(Y_{Xt},P_{Xt},Q_{Xt})$ of FBSDEs \eqref{FB:1}. For $X,X'\in\hr_m$, both are independent of $\mathcal{W}_t$, we then have
	\begin{align}
		&\sup_{t\le s\le T}\|Y_{Xt}(s)\|_{\hr_m}+\sup_{t\le s\le T}\|P_{Xt}(s)\|_{\hr_m}+\sum_{j=1}^n\|Q^j_{Xt}\|_{L^2(t,T;\hr_m)} \notag\\
		\le\ & C(L,T,\lambda)(1+\|X\|_{\hr_m}), \label{thm2_1}\\
		&\sup_{t\le s\le T}\left\|Y_{X't}(s)-Y_{Xt}(s)\right\|_{\hr_m}+\sup_{t\le s\le T}\left\|P_{X't}(s)-P_{Xt}(s)\right\|_{\hr_m}+\sum_{j=1}^n \|Q^j_{X't}-Q^j_{Xt}\|_{L^2(t,T;\hr_m)} \notag\\
		\le\ & C(L,T,\lambda)\left\|X'- X\right\|_{\hr_m}.\label{thm2_2}
	\end{align}	
	Furthermore, if Assumption (A3)-(ii) is satisfied, the FBSDEs \eqref{FB:1} has a unique solution satisfying \eqref{thm2_1} and \eqref{thm2_2} for any $\lambda> 0$.
\end{theorem}

\subsection{Jacobian flow of FBSDEs \eqref{FB:1}}\label{subsec:Derivatives}

We next consider the differentiability of $(Y_{Xt}(s),P_{Xt}(s),Q_{Xt}(s))$ in the initial data $X\in\hr_m$. In view of Assumptions (A2) and (A3)-(i) and the Young's inequality with a weight of $\frac{1}{2\sqrt{\lambda}}$, for any $X,Y\in\hr_m$, $V,U\in\ur_m$ and $s\in[0,T]$,
\begin{equation}\label{mono_11}
	\begin{split}
		&\left(D_X^2 G(X,V,s)(Y),Y\right)_{\hr_m}+2\left(D_VD_X G(X,V,s)(U),Y\right)_{\hr_m}+\left(D_V^2 F(X,V,s)(U),U\right)_{\ur_m}\\
		\geq\ & \lambda \|U\|^2_{\ur_m} -C(L)\left(1+\frac{1}{\lambda}\right)\|Y\|^2_{\hr_m},
	\end{split}
\end{equation}
Furthermore, if Assumption (A3)-(ii) is satisfied, then,
\begin{equation}\label{mono_12}
	\begin{split}
		&\left(D_X^2 G(X,V,s)(Y),Y\right)_{\hr_m}+2\left(D_VD_X G(X,V,s)(U),Y\right)_{\hr_m}+\left(D_V^2 F(X,V,s)(U),U\right)_{\ur_m}\geq  2\lambda \|U\|^2_{\ur_m};\\
		&\left(D_X^2 G_T(X)(Y),Y\right)_{\hr_m} \geq 0.
	\end{split}
\end{equation}
From Assumption (A2) and Condition \eqref{hv2}, we know that the map $\hv$ (defined in \eqref{hv}) is continuously differentiable in $(X,P,Q)\in\hr_m\times\hr_m\times\hr_m^n$ and satisfies, for any $\tx,\tp,\tq^j\in\hr_m$,
\begin{equation}\label{optimal_condition}
	\begin{split}
		&D_XD_V G\left(X,\hv(X,s;P,Q),s\right)\left(\tx\right)\\
		&\quad+D^2_V G\left(X,\hv(X,s;P,Q),s\right)\left(D_X\hv(X,s;P,Q)\left(\tx\right)\right)=0;\\
		&\mathcal{F}^*_2(s)\tp+D^2_V G\left(X,\hv(X,s;P,Q),s\right)\left(D_P\hv(X,s;P,Q)\left(\tp\right)\right)=0;\\
		&{\mathcal{A}_2^j}^*(s)\tq^j+D^2_V G\left(X,\hv(X,s;P,Q),s\right)\left(D_{Q^j}\hv(X,s;P,Q)\left(\tq^j\right)\right)=0.
	\end{split}
\end{equation}
For the initial $(t,X)$ and $\tx\in\hr_m$, which is independent of $\mathcal{W}_t$, consider the following FBSDEs: for $s\in[t,T]$,
\begin{equation}\label{FB:6}
	\left\{
	\begin{aligned}
		&\mathcal{D}_{\tx} Y_{Xt}(s)=\tx+\int_t^s \left[\mathcal{F}_1(r)\mathcal{D}_{\tx} Y_{X t}(r)+\mathcal{F}_2(r)\mathcal{D}_{\tx} \hat{V}_{Xt}(r)\right]dr\\
		&\qquad\qquad\qquad +\sum_{j=1}^n \int_t^s\left[\mathcal{A}^j_1(r)\mathcal{D}_{\tx} Y_{X t}(r)+\mathcal{A}^j_2(r)\mathcal{D}_{\tx} \hat{V}_{X t}(r)
		\right]dw_j(r),\\
		&\mathcal{D}_{\tx}P_{Xt}(s)=D^2_X G_T(Y_{Xt}(T))\left(\mathcal{D}_{\tx}Y_{Xt}(T)\right)\\
		&\qquad\qquad\qquad +\int_s^T \bigg[\mathcal{F}^*_1(r)\mathcal{D}_{\tx}P_{Xt}(r)+\sum_{j=1}^n {\mathcal{A}^j_1}^*(r)\mathcal{D}_{\tx}Q^j_{Xt}(r)\\
		&\quad\qquad\qquad\qquad\qquad +D_X^2 G\left(Y_{Xt}(r),\hat{V}_{Xt}(r),r\right)\left(\mathcal{D}_{\tx}Y_{Xt}(r)\right)\\
		&\quad\qquad\qquad\qquad\qquad +D_V D_X G\left(Y_{Xt}(r),\hat{V}_{Xt}(r),r\right)\left(\mathcal{D}_{\tx}\hat{V}_{Xt}(r)\right)\bigg]dr\\
		&\qquad\qquad\qquad -\sum_{j=1}^n \int_s^T \mathcal{D}_{\tx}Q_{Xt}^j(r)dw_j(r),
	\end{aligned}
	\right.
\end{equation}
where $\hat{V}_{Xt}(s):=\hv(Y_{Xt}(s),s;P_{Xt}(s),Q_{Xt}(s))$ and
\begin{align*}
	\mathcal{D}_{\tx}\hat{V}_{Xt}(s):=&D_X\hv(Y_{Xt}(s),s;P_{Xt}(s),Q_{Xt}(s))\left(\mathcal{D}_{\tx}Y_{Xt}(s)\right)\\
	& +D_P\hv(Y_{Xt}(s),s;P_{Xt}(s),Q_{Xt}(s))\left(\mathcal{D}_{\tx}P_{Xt}(s)\right)\\
	& +\sum_{j=1}^nD_{Q^j}\hv(Y_{Xt}(s),s;P_{Xt}(s),Q_{Xt}(s))\left(\mathcal{D}_{\tx}Q^j_{Xt}(s)\right).
\end{align*}
We have the following solvability of FBSDEs \eqref{FB:6} and the boundedness of the solution, whose proof is given in Appendix~\ref{pf:thm4}.

\begin{lemma}\label{thm4}
	Under Assumptions (A1), (A2) and (A3)-(i), there exists a constant $c(L,T)$ depending only on $(L,T)$, such that when $\lambda\geq c(L,T)$, there is a unique solution $\left(\mathcal{D}_{\tx}Y_{Xt},\mathcal{D}_{\tx}P_{Xt},\mathcal{D}_{\tx}Q_{Xt}\right)$ of the FBSDEs \eqref{FB:6}, such that
	\begin{equation}\label{thm4_0}
		\begin{split}
			&\sup_{t\le s\le T}\left\|\mathcal{D}_{\tx}Y_{Xt}(s)\right\|_{\hr_m}+\sup_{t\le s\le T}\left\|\mathcal{D}_{\tx}P_{Xt}(s)\right\|_{\hr_m}+\sum_{j=1}^n\left\|\mathcal{D}_{\tx}Q^j_{Xt}\right\|_{L^2(t,T;\hr_m)}\\
			\le\ & C(L,T,\lambda)\left\|\tx\right\|_{\hr_m}.
		\end{split}
	\end{equation}   
	Furthermore, if Assumption (A3)-(ii) is satisfied, then, for any $\lambda> 0$, FBSDEs \eqref{FB:6} has a unique solution satisfying \eqref{thm4_0}.
\end{lemma}
We then have the following G\^ateaux differentiability of $(Y_{Xt}(s),P_{Xt}(s),Q_{Xt}(s))$ in the initial data $X$, which is proven in Appendix~\ref{pf:thm5}; we also refer to \cite{AB5} for further discussion on their very nature as Fr\'echet derivatives.

\begin{theorem}\label{thm8}
	Under Assumptions (A1), (A2) and (A3)-(i), there exists a constant $c(L,T)$ depending only on $(L,T)$, such that when $\lambda\geq c(L,T)$, for any $\tx\in\hr_m$, which is independent of $\mathcal{W}_t$, we have
	\begin{equation}\label{thm5_0}
		\begin{split}
			\lim_{\epsilon\to0}\Bigg[ & \sup_{t\le s\le T}\left\|\frac{ Y_{X^\epsilon t}(s)-Y_{Xt}(s)}{\epsilon}-\mathcal{D}_{\tx}Y_{Xt}(s)\right\|^2_{\hr_m}\\
			&+ \sup_{t\le s\le T}\left\|\frac{ P_{X^\epsilon t}(s)-P_{Xt}(s)}{\epsilon}-\mathcal{D}_{\tx}P_{Xt}(s)\right\|^2_{\hr_m}\\
			&+\sum_{j=1}^n \int_t^T \left\|\frac{ Q^j_{X^\epsilon t}(s)-Q^j_{Xt}(s)}{\epsilon}-\mathcal{D}_{\tx}Q^j_{Xt}(s)\right\|^2_{\hr_m} ds\Bigg]=0,
		\end{split}
	\end{equation}
	where $X^\epsilon=X+\epsilon\tx$ for $\epsilon\in(0,1)$. That is, $\left(\mathcal{D}_{\tx}Y_{Xt}(s),\mathcal{D}_{\tx}P_{Xt}(s),\mathcal{D}_{\tx}Q_{Xt}(s)\right)$ defined in \eqref{FB:6} are the G\^ateaux derivatives $\left(D_{\tx}Y_{Xt}(s),D_{\tx}P_{Xt}(s),D_{\tx}Q_{Xt}(s)\right)$ of $(Y_{Xt}(s),P_{Xt}(s),Q_{Xt}(s))$ in $X$ along the direction $\tx$. Moreover, the G\^ateaux derivatives satisfy the boundedness property \eqref{thm4_0}, and are linear in $\tx$ and continuous in $X$. Furthermore, if Assumption (A3)-(ii) is satisfied, for any $\lambda> 0$, the same assertion holds.  
\end{theorem}

\section{Regularity of the Value Function in $(t,X)$}\label{sec:value_Hilbert}

In view of the lifting method \eqref{GAF} and the FBSDEs \eqref{FB:1}, for $(t,m)\in [0,T]\times\pr_2(\br^n)$ and $X\in \hr_m$, which is independent of $\mathcal{W}_t$, we have
\begin{align}\label{v_1}
	\V(X\otimes m,t)=\int_t^T F\left(Y_{Xt}(s),\hv_{Xt}(s),s\right)ds+F_T(Y_{Xt}(T)),
\end{align}
where $\hv_{Xt}(s)=\hv(Y_{Xt}(s),s;P_{Xt}(s),Q_{Xt}(s)),\  s\in(t,T)$. We then establish in the following the G\^ateaux differentiability of $\V$ in $X$, and its proof  is given in Appendix~\ref{pf:lem:8}.

\begin{lemma}\label{lem:8}
	Under Assumptions (A1), (A2) and (A3)-(i), there exists a constant $c(L,T)$ depending only on $(L,T)$, such that when $\lambda\geq c(L,T)$, the value function $\V$ is twice G\^ateaux differentiable in $X$ with the derivatives
	\begin{align}
		D_X \V(X\otimes m,t)=&\ P_{Xt}(t),\label{lem8_1}\\
		D^2_X \V(X\otimes m,t)\left(\tx\right)=&\ {D}_{\tx}P_{Xt}(t),\label{lem8_1.1}
	\end{align}
	and they satisfy the growth conditions:
	\begin{align}
		|\V(X\otimes m,t)|\le&\  C(L,T,\lambda)\left(1+\|X\|^2_{\hr_m}\right),\label{lem8_3}\\
		\|D_X \V(X\otimes m,t)\|_{\hr_m}\le&\  C(L,T,\lambda)(1+\|X\|_{\hr_m}),\label{lem8_4}\\
		\left\|D^2_X \V(X\otimes m,t)\left(\tx\right)\right\|_{\hr_m}\le&\  C(L,T,\lambda)\left\|\tx\right\|_{\hr_m},\label{lem8_4'}
	\end{align}
	and we also have the continuity condition
	\begin{equation}\label{lem8_5}
		\begin{split}
			&\left\|D_X \V(X'\otimes m,t)-D_X \V(X\otimes m,t)\right\|_{\hr_m} \le C(L,T,\lambda)\left\|X'-X\right\|_{\hr_m},
		\end{split}
	\end{equation}
	and $D_X^2 \V$ is linear in $\tx$ and continuous in $X$. Furthermore, if Assumption (A3)-(ii) is satisfied, then, for any $\lambda> 0$, all the above assertions remain valid.
\end{lemma}
To give the regularity of $\V$ in time $t$, we need the following additional assumption on coefficients $\left(F,A^j,G\right)$ in $t$.\\
\textbf{(A4)} The coefficients $F$ and $G$ are continuous in $t\in[0,T]$, and $\mathcal{A}^j_2=0$ for $1\le j\le n$. 

\begin{remark}
	In general, even for classical control problems without mean-field term, the Bellman equation barely has a classical solution when the diffusion depends on the control. We refer to Pham-Wei \cite{PH} and Tang-Zhang \cite{TS} for viscosity solution of Bellman equations when the diffusion depends on the control in different settings.
\end{remark}
Next we give the regularity of $\V$ in $t$ in the following, and its proof  is given in Appendix~\ref{pf:lem9}.

\begin{lemma}\label{lem9}
	Under Assumptions (A1), (A2) (A3)-(i) and (A4), there exists a constant $c(L,T)$ depending only on $(L,T)$, such that when $\lambda\geq c(L,T)$, for $0\le t\le t'\le T$, 
	\begin{equation}\label{lem9_1}
		\begin{split}
			&\sup_{t\le s\le T}\left\|Y_{Xt'}(s)-Y_{Xt}(s)\right\|_{\hr_m}+\sup_{t\le s\le T}\left\|P_{Xt'}(s)-P_{Xt}(s)\right\|_{\hr_m} \\
			&+\sum_{j=1}^n \left\|Q^j_{Xt'}-Q^j_{Xt}\right\|_{L^2(t,T;\hr_m)}\le C(L,T,\lambda)(1+\|X\|_{\hr_m})|t'-t|^{\frac{1}{2}},
		\end{split}
	\end{equation}
	and $\left(\mathcal{D}_{\tx}Y_{Xt}(s),\mathcal{D}_{\tx}P_{Xt}(s),\mathcal{D}_{\tx}Q_{Xt}(s)\right)$ is continuous in $t$. The derivative $D_X\V(X\otimes m,t)$ satisfies
	\begin{equation}\label{lem9_2}
		\left\|D_X\V\left(X\otimes m,t'\right)-D_X\V(X\otimes m,t)\right\|_{\hr_m} \le C(L,T,\lambda)\left(1+\|X\|_{\hr_m} \right)\left|t'-t\right|^{\frac{1}{2}},
	\end{equation}	
	and the derivative $D^2_X\V(X\otimes m,t)$ is continuous in $t$. The value function $\V$ is $C^1$ in $t$ with the derivative
	\begin{equation}\label{lem9_3}
		\begin{split}
			\frac{\dd\V}{\dd t}(X\otimes m,t)=&-G(X,\hv_{Xt}(t),t)-\left(D_X\V(X\otimes m,t),\ G\left(X, \hv_{Xt}(t),t\right)\right)_{\hr_m}\\
			&-\frac{1}{2}\left(D_X^2\V(X\otimes m,t)\left(\sum_{j=1}^n A^j(X,t)\mathcal{N}_t^j\right),\   \sum_{j=1}^n A^j(X,t)\mathcal{N}_t^j\right)_{\hr_m},
		\end{split}
	\end{equation}
	where $\mathcal{N}_t^j$'s, for $j=1,2,\dots,n$, are $n$ independent standard normal random variables which are also independent of $X$ and the underlying Wiener Process. The derivative $\frac{\dd\V}{\dd t}$ also satisfies the growth condition
	\begin{align}\label{lem9_4}
		&\left|\frac{\dd\V}{\dd t}(X\otimes m,t)\right|\le C(L,T,\lambda)\left(1+\|X\|^2_{\hr_m}\right).
	\end{align}    
	Furthermore, if Assumption (A3)-(ii) is satisfied, then, for any $\lambda> 0$, the above assertion remains valid.
\end{lemma}

Now we introduce the following Bellman equation: for $(t,m)\in[0,T]\times\pr_2(\brn)$, $X\in\hr_m$, which is also independent of  $\mathcal{W}_t$,
\begin{equation}\label{Bellman}
	\left\{
	\begin{aligned}
		&\frac{\dd \V}{\dd t}(X\otimes m,t)+\frac{1}{2}\left(D_X^2\V(X\otimes m,t)\left(\sum_{j=1}^n A^j(X,t)\mathcal{N}_t^j\right),\   \sum_{j=1}^n A^j(X,t)\mathcal{N}_t^j\right)_{\hr_m}\\
		&\qquad\qquad\qquad +\left(D_X\V(X\otimes m,t),\ F\left(X, \hv(X,t;D_X \V(X\otimes m,t)),t\right)\right)_{\hr_m}\\
		&\qquad\qquad\qquad +G\left(X,\hv(X,t;D_X \V(X\otimes m,t)),t\right)=0,\\
		&\V(X\otimes m,T)=G_T(X),
	\end{aligned}
	\right.
\end{equation}
where these $\mathcal{N}_t^j$'s are specified in Lemma~\ref{lem9}. We therefore have the solvability of the Bellman equation \eqref{Bellman}, and its proof is given in Appendix~\ref{pf:thm:9}.

\begin{theorem}\label{thm:9}
	Under Assumptions (A1), (A2), (A3)-(i) and (A4), there exists a constant $c(L,T)$ depending only on $(L,T)$, such that when $\lambda\geq c(L,T)$, the value function $\V$ defined in \eqref{def:V} is the unique solution of the Bellman equation \eqref{Bellman} possessing the standard regularity of  \eqref{lem8_3}-\eqref{lem8_5}, \eqref{lem9_2} and \eqref{lem9_4}. Furthermore, if Assumption (A3)-(ii) is satisfied, then, for any $\lambda> 0$, the same assertion remains valid.
\end{theorem}

As a consequence of Lemmas~\ref{lem:8} and \ref{lem9}, we have the following sensitivity relations (as called in Bonnet-Frankowska \cite{FRA1}) between the maximum principle and the G\^ateaux derivatives of the value function $\V$.

\begin{corollary}\label{cor1}
	Under assumptions in Lemma~\ref{lem9}, for any $s\in[t,T]$, the following sensitivity relations hold:
	\begin{align}
		D_X\V(Y_{Xt}(s)\otimes m,s)=&P_{Xt}(s),\label{cor1_1}\\
		D^2_X\V(Y_{Xt}(s)\otimes m,s)\left(\sum_{j=1}^nA^j(Y_{Xt}(s),s)w_j(s)\right)=&\sum_{j=1}^nQ^j_{Xt}(s)w_j(s).\label{cor1_2}
	\end{align}
\end{corollary}

\section{Regularity of Solution of FBSDEs \eqref{FB:4}}\label{sec:MFC}

We now go back to the Problem~\ref{intr_4}. From the respective definitions for $G,A,F,\lr,\hv$ and $\bh$ in \eqref{GAF}, \eqref{H''}, \eqref{v} and \eqref{bh''}, we know that FBSDEs \eqref{FB:4} concide with FBSDEs \eqref{FB:1}. In this section, we apply our results in Section~\ref{sec:FB_Hilbert} to study the regularity of FBSDEs \eqref{FB:4}. Our assumptions for coefficients $f,\sigma$ and $g$ are as follows. For the sake of convenience, we use the same constant $l>0$ for all the conditions below. \\
\textbf{(B1)} The functions $f$ and $\sigma^j$ (for $1\le j\le n$) are linear. That is,
\begin{align*}
	&(i)\ f(x,m,v,s)=f_0(s)+f_1(s)x+f_2(s)\int_\brn ydm(y)+f_3(s)v,\\
	&(ii)\ \sigma^j(x,m,v,s)=\sigma^j_0(s)+\sigma^j_1(s)x+\sigma^j_2(s)\int_\brn ydm(y)+\sigma^j_3(s)v,
\end{align*}
with $f_0(s),\sigma^j_0(s)\in \brn$, $f_1(s),f_2(s),\sigma^j_1(s),\sigma^j_2(s)\in\br^{n\times n}$ and $f_3(s),\sigma^j_3(s)\in\br^{n\times d}$, and all their norms are bounded above by $l$.\\
\textbf{(B2)} The running cost function $g$ and the terminal cost function $g_T$ satisfy: for $(x,m,v,s)\in \brn\times \pr_2(\brn)\times \brd\times[0,T]$,
\begin{align*}
	|g(x,m,v,s)|\le& l\left(1+|x|^2+W^2_2(m,\delta_0)+|v|^2\right),\\  
	|g_T(x,m)|\le& l\left(1+|x|^2+W^2_2(m,\delta_0)\right).
\end{align*}
The following derivatives exist, and they are continuous in all of their own arguments, and their norms are bounded above by $l$:
\begin{align*}
	&D_x^2 g, \ D_vD_x g, \ D_v^2 g, \ D_x D_\xi \frac{dg}{d\nu}, \ D_v D_\xi \frac{dg}{d\nu},\ D_\xi^2 \frac{dg}{d\nu}, \  D_{\xi'}D_\xi\frac{d^2g}{d\nu^2},\\
	&D_x^2 g_T,\ D_x D_\xi \frac{dg_T}{d\nu},\ D_\xi^2 \frac{dg_T}{d\nu}, \ D_{\xi'}D_\xi\frac{d^2g_T}{d\nu^2}.
\end{align*}
\textbf{(B3)} (i) There exists $\lambda> 0$ such that for any $(x,m,v,v')\in\brn\times\pr_2(\brn)\times\brd\times\brd$ and $s\in[0,T]$,
\begin{align*}
	&g(x,m,v',s)-g(x,m,v,s)\geq \left(D_v g (x,m,v,s)\right)^* (v'-v)+\lambda |v'-v|^2.
\end{align*}
(ii) There exists $\lambda>0 $ such that, for any $x,x',\xi,\xi'\in\brn$, $v,v'\in\brd$ and $s\in[0,T]$,
\begin{align*}
	&(\alpha1)\ g\left(x',m',v',s\right)-g(x,m,v,s)\geq \left(D_x g (x,m,v,s)\right)^* \left(x'-x\right)\\
	&\qquad\qquad\qquad\qquad\qquad\qquad\qquad\qquad +\int_\brn\frac{d g}{d\nu}(x,m,v,s)(\xi)d\left(m'-m\right)(\xi)\\
	&\qquad\qquad\qquad\qquad\qquad\qquad\qquad\qquad +\left(D_v g (x,m,v,s)\right)^* \left(v'-v\right)+\lambda \left|v'-v\right|^2;\\
	&(\alpha2)\ \frac{dg}{d\nu}(x,m,v,s)\left(\xi'\right)-\frac{dg}{d\nu}(x,m,v,s)(\xi)\geq \left(D_\xi \frac{dg}{d\nu}(x,m,v,s)(\xi)\right)^* \left(\xi'-\xi\right);\\
	&(\beta1)\ g_T\left(x',m'\right)-g_T(x,m)\geq \left(D_x g_T (x,m)\right)^* \left(x'-x\right)+\int_\brn\frac{d g_T}{d\nu}(x,m)(\xi)d\left(m'-m\right)(\xi);\\
	&(\beta2)\ \frac{dg_T}{d\nu}(x,m)\left(\xi'\right)-\frac{dg_T}{d\nu}(x,m)(\xi)\geq \left(D_\xi \frac{dg_T}{d\nu}(x,m)(\xi)\right)^*  \left(\xi'-\xi\right).
\end{align*}
We can see that, under Assumptions (B1) and (B3)-(i), the maps $F$ and $A$ defined in \eqref{GAF} satisfy (A1), and $G$ and $G_T$ defined in \eqref{GAF} satisfy (A3)-(i), with a constant $L=C(l)$, where $C(l)$ is a fixed polynomial of $l$. Furthermore, Assumption (B3)-(ii) yields (A3)-(ii). We refer to Section 6 in \cite{AB8} for more details. Moreover, for $s\in[0,T]$, $X,\tx\in \hr_m$, $V,\tv\in \ur_m$, we have for $x\in\brn$,
\begin{align*}
	&(\gamma1)\ \left.D_V^2 G(X,V,s)\left(\tv\right)\right|_x=\left(D^2_v g(X_x,X\otimes m,V_x,s)\right)^*\tv_x,\\
	&(\gamma2)\ \left.D_VD_X G(X,V,s)\left(\tv\right)\right|_x=\left( D_vD_x g(X_x,X\otimes m,V_x,s)\right)^*\tv_x\\
	&\qquad\qquad\qquad\qquad\qquad\qquad\qquad +\bar{\e}\left[\int_\brn \left( D_v D_\xi \frac{d g}{d\nu}\left(\bar{X}_y,X\otimes m,\bar{V}_y,s\right)(X_x)\right)^* \bar{\tv}_y dm(y)\right],\\
	&(\gamma3)\ \left.D_XD_V G(X,V,s)\left(\tx\right)\right|_x=\left(D_v D_x g(X_x,X\otimes m,V_x,s)\right)^*\tx_x\\
	&\qquad\qquad\qquad\qquad\qquad\qquad\qquad +\bar{\e}\left[\int_\brn \left(D_v D_\xi\frac{d g}{d\nu}(X_x,X\otimes m,V_x,s)\left(\bar{X}_y\right)\right)^*\bar{\tx}_y dm(y)\right],\\
	&(\gamma4)\ \left.D_X^2 G(X,V,s)\left(\tx\right)\right|_x\\
	&\quad =\left(D^2_x g(X_x,X\otimes m,V_x,s)\right)^*\tx_x+\bar{\e}\left[\int_\brn \left(D_\xi^2 \frac{d g}{d\nu}\left(\bar{X}_y,X\otimes m,\bar{V}_y,s\right)(X_x)\right)^* dm(y)\right]\tx_x\\
	&\qquad +\bar{\e}\left[\int_\brn \left(D_x D_\xi\frac{d g}{d\nu}(X_x,X\otimes m,V_x,s)\left(\bar{X}_y\right)\right)^*\bar{\tx}_y dm(y)\right]\\
	&\qquad+\bar{\e}\left[\int_\brn \left( D_x D_\xi \frac{d g}{d\nu}\left(\bar{X}_y,X\otimes m,\bar{V}_y,s\right)(X_x)\right)^* \bar{\tx}_y dm(y)\right]\\
	&\qquad+\bar{\e}\left[\int_\brn \check{\e}\left(\int_\brn \left( D_{\xi'}D_\xi \frac{d^2 g}{d\nu^2}\left(\bar{X}_y,X\otimes m,\bar{V}_y,s\right)\left(X_x,\cx_z\right)\right)^* \check{\tx}_z dm(z)\right) dm(y)\right],
\end{align*}
where $\left(\bar{X},\bar{\tx}\right)$ and $\left(\cx,\check{\tx}\right)$ are two independent copies of $\left(X,\tx\right)$, and $\left(\bar{V},\bar{\tv}\right)$ is an independent copy of $\left(V,\tv\right)$. Under Assumption (B2), all these items $(\gamma1)$-$(\gamma4)$ with $G$ and $G_T$ defined in \eqref{GAF} can be easily shown to satisfy (A2).

From Assumption (B2) and the definitions of the maps $\hat{v}$ in \eqref{bh'} and $\hat{V}$ in \eqref{v}, since $\hat{v}$ is continuously differentiable in $(x,m,p,q)$, then $\hv$ can be shown to be G\^ateaux differentiable in $(X,P,Q)$ (see \cite{AB6,AB8,AB5,AB9}). Moreover, for $s\in[0,T]$, $X,P,\tx,\tp\in\hr_m$ and $Q,\tq\in\hr_m^n$,
\begin{align*}
	&(\delta1)\ \left.D_X \hv(X,s;P,Q)\left(\tx\right)\right|_x=\left(D_x \hat{v}(X_x,X\otimes m,s;P_x,Q_x)\right)^*\tx_x\\
	&\qquad\qquad\qquad\qquad\qquad\qquad\qquad +\bar{\e}\left[\int_\brn \left(D_\xi\frac{d \hat{v}}{d\nu}(X_x,X\otimes m,s;P_x,Q_x)(\bar{X}_y)\right)^* \bar{\tx}_y dm(y)\right],\\
	&(\delta2)\ \left.D_P \hv(X,s;P,Q)\left(\tp\right)\right|_x=\left(D_p \hat{v}(X_x,X\otimes m,s;P_x,Q_x)\right)^*\tp_x,\\
	&(\delta3)\ \left.D_{Q^j} \hv(X,s;P,Q)\left(\tq^j\right)\right|_x=\left(D_{q^j} \hat{v}(X_x,X\otimes m,s;P_x,Q_x)\right)^*\tq^j_x.
\end{align*} 
In view of the above calculation and FBSEDs \eqref{FB:6}, the Jacobian flow
\begin{align*}
	\left(D_{\tx_x}Y_{X_x mt}(s),D_{\tx_x}P_{X_x mt}(s),D_{\tx_x}Q_{X_x mt}(s)\right)
\end{align*}
of $(Y_{X_x mt}(s),P_{X_x mt}(s),Q_{X_x mt}(s))$ in $X$ along the direction $\tx\in\hr_m$ can be characterized as the solution of the following FBSDEs: for $(s,x)\in[t,T]\times\brn$,
\small
\begin{align}
	&D_{\tx_x}Y_{X_xmt}(s) \notag\\
	=&\tx_x+\int_t^s \left[\mathcal{F}_1(r)D_{\tx_x}Y_{X_xmt}(r)+\mathcal{F}_2(r)D_{\tx_x}\hvv_{X_xmt}(r)\right]dr \notag\\
	& +\sum_{j=1}^n \int_t^s\left[\mathcal{A}^j_1(r)D_{\tx_x}Y_{X_xmt}(r)+\mathcal{A}^j_2(r)D_{\tx_x}\hvv_{X_xmt}(r)\right]dw_j(r); \notag\\
	&D_{\tx_x}P_{X_xmt}(s)\notag\\
	=&\left(D^2_x g_T(Y_{X_x mt}(T),Y_{Xmt}(T)\otimes m)\right)^*D_{\tx_x}Y_{X_xmt}(T) \notag\\
	& +\bar{\e}\left[\int_\brn \left(D_\xi^2 \frac{d g_T}{d\nu}\left(\bar{Y}_{\bx_y mt}(T),Y_{Xmt}(T)\otimes m \right)(Y_{X_x mt}(T))\right)^* dm(y)\right]D_{\tx_x}Y_{X_x mt}(T) \notag\\
	&+\bar{\e}\left[\int_\brn \left( D_x D_\xi\frac{d g_T}{d\nu}\left(Y_{X_x mt}(T),Y_{Xmt}(T)\otimes m \right)\left(\bar{Y}_{\bx_y mt}(T)\right)\right)^* D_{\bar{\tx}_y}\bar{Y}_{\bx_y mt}(T)  dm(y)\right] \notag\\
	&+\bar{\e}\left[\int_\brn \left(D_x D_\xi \frac{d g_T}{d\nu}\left(\bar{Y}_{\bx_y mt}(T),Y_{Xmt}(T)\otimes m\right)(Y_{X_x mt}(T))\right)^* D_{\bar{\tx}_y}\bar{Y}_{\bx_y mt}(T) dm(y)\right] \notag\\
	&+\bar{\e}\bigg[\int_\brn \check{\e}\bigg(\int_\brn \left(D_{\xi'}D_\xi \frac{d^2 g_T}{d\nu^2}\left(\bar{Y}_{X_y mt}(T),Y_{Xmt}(T)\otimes m \right)\left(Y_{X_x mt}(T),\cy_{\cx_z mt}(T)\right)\right)^* \notag\\
	&\quad\qquad\qquad\qquad\qquad D_{\check{\tx}_z} \cy_{\cx_z mt}(T) dm(z)\bigg) dm(y)\bigg] \notag\\
	&+\int_s^T \Bigg\{\mathcal{F}^*_1(r) D_{\tx_x}P_{X_x mt}(r)+\sum_{j=1}^n {\mathcal{A}^j_1}^*(r) D_{\tx_x}Q^j_{X_x mt}(r) \notag\\
	&\quad\quad\qquad +\left(D^2_x g\left(Y_{X_x mt}(r),Y_{Xmt}(r)\otimes m,\hvv_{X_x mt}(r),r\right)\right)^*D_{\tx_x}Y_{X_x mt}(r) \notag\\
	&\qquad\qquad+\bar{\e}\bigg[\int_\brn \left(D_\xi^2 \frac{d g}{d\nu}(\bar{Y}_{\bx_y mt}(r),Y_{Xmt}(r)\otimes m,\bar{\hvv}_{\bx_y mt}(r),r)(Y_{X_x mt}(r))\right)^* D_{\tx_x}Y_{X_x mt}(r) dm(y)\bigg] \notag\\
	&\qquad\qquad+\bar{\e}\bigg[\int_\brn \left(D_x D_\xi\frac{d g}{d\nu}\left(Y_{X_x mt}(r),Y_{Xmt}(r)\otimes m,\hvv_{X_x mt}(r),r\right)\left(\bar{Y}_{\bx_y mt}(r)\right)\right)^* D_{\bar{\tx}_y} \bar{Y}_{\bx_y mt}(r) dm(y)\bigg] \notag\\
	&\qquad\qquad+\bar{\e}\bigg[\int_\brn \left(D_x D_\xi \frac{d g}{d\nu}\left(\bar{Y}_{\bx_y mt}(r),Y_{Xmt}(r)\otimes m,\bar{\hvv}_{\bx_y mt}(r),r\right)(Y_{X_x mt}(r))\right)^* D_{\bar{\tx}_y}\bar{Y}_{\bx_y mt}(r)  dm(y)\bigg] \notag\\
	&\qquad\qquad+\bar{\e}\bigg[\int_\brn \check{\e}\bigg[\int_\brn \bigg(D_{\xi'}D_\xi \frac{d^2 g}{d\nu^2}\left(\bar{Y}_{\bx_y mt}(r),Y_{Xmt}(r)\otimes m,\bar{\hvv}_{\bx_y mt}(r),r\right) \notag\\
	&\ \ \quad\qquad\qquad\qquad\qquad\qquad \left(Y_{X_x mt}(r),\cy_{\cx_z mt}(r)\right)\bigg)^* D_{\check{\tx}_z} \cy_{\cx_z mt}(r) dm(z)\bigg] dm(y)\bigg] \notag\\
	&\qquad\qquad +\left(D_vD_x g\left(Y_{X_x mt}(r),Y_{Xmt}(r)\otimes m,\hvv_{X_x mt}(r),r\right)\right)^* D_{\tx_x}u_{X_x mt}(r) \notag\\
	&\qquad\qquad +\bar{\e}\bigg[\int_\brn \left(D_v D_\xi \frac{d g}{d\nu}\left(\bar{Y}_{\bx_y mt}(r),Y_{Xmt}(r)\otimes m,\bar{\hvv}_{\bx_y mt}(r),r\right)(Y_{X_x mt}(r))\right)^*  \notag\\
	&\quad\qquad\qquad\qquad\qquad D_{\bar{\tx}_y}\bar{\hvv}_{\bx_y mt}(r) dm(y)\bigg]\Bigg\}dr -\int_s^T \sum_{j=1}^n D_{\tx_x}Q^j_{X_x mt}(r)dw_j(r), \label{FB:7}
\end{align}
\normalsize
where $\hvv_{X_x mt}(s)$ is defined in \eqref{max_MFC}, $D_{\tx_x}\hvv_{X_x mt}(s)$ can be shown to be
\begin{align*}
	&D_{\tx_x}\hvv_{X_x mt}(s)\\
	=&(D_x \hat{v}(Y_{X_x mt}(s),Y_{Xmt}(s)\otimes m,s;P_{X_x mt}(s),Q_{X_x mt}(s)))^*D_{\tx_x}Y_{X_x mt}(s)\\
	&+\bar{\e}\bigg[\int_\brn \bigg(D_\xi\frac{d \hat{v}}{d\nu}(Y_{X_x mt}(s),Y_{Xmt}(s)\otimes m,s;P_{X_x mt}(s),Q_{X_x mt}(s))\left(\bar{Y}_{\bx_y mt}(s)\right)\bigg)^* \\
	&\quad\qquad\qquad D_{\bar{\tx}_y} \bar{Y}_{\bx_y mt}(s) dm(y)\bigg],\\
	&+(D_p \hat{v}(Y_{X_x mt}(s),Y_{Xmt}(s)\otimes m,s;P_{X_x mt}(s),Q_{X_x mt}(s)))^*D_{\tx_x}P_{X_x mt}(s)\\
	& +\sum_{j=1}^n (D_{q^j} \hat{v}(Y_{X_x mt}(s),Y_{Xmt}(s)\otimes m,s;P_{X_x mt}(s),Q_{X_x mt}(s)))^*D_{\tx_x}Q^j_{X_x mt}(s),
\end{align*}
and $\left(\bx_y, \bar{\tx}_y, \bar{Y}_{\bx_y mt}(s), D_{\bar{\tx}_y} \bar{Y}_{\bx_y mt}(s) \right)$ is an independent copy of $\left(X_y, \tx_y, {Y}_{X_y mt}(s), D_{{\tx}_y} {Y}_{X_y mt}(s) \right)$ for all $y\in\brn$, which is viable due to the continuity of each component process; similarly, process $\left(\cx_z,\check{\tx}_z,\cy_{\cx_z mt}(s), D_{\check{\tx}_z} \cy_{\cx_z mt}(s) \right)$ is another independent copy of $\left(X_z, \tx_z, {Y}_{X_z mt}(s), D_{{\tx}_z} {Y}_{X_z mt}(s) \right)$ for all $z\in\brn$, and $\left(\bar{\hvv}_{\bx_y mt}(s), D_{\bar{\tx}_y} \bar{\hvv}_{\bx_y mt}(s) \right)$ is an independent copy of $\left(\hvv_{X_y mt}(s), D_{{\tx}_y} \hvv_{X_y mt}(s) \right)$ for all $y\in\brn$. We also refer to \cite{AB5} for further discussion on their Fr\'echet derivative nature of these G\^ateaux derivatives. As a consequence of Theorems~\ref{thm2} and \ref{thm8}, we have the following solvability and regularity of FBSDEs \eqref{FB:4}.

\begin{theorem}\label{thm3}
	Under Assumptions (B1), (B2) and (B3)-(i), there exists a constant $c(l,T)$ depending only on $(l,T)$, such that when $\lambda\geq c(l,T)$, there is a unique solution $(Y_{X_\cdot mt}(\cdot),P_{X_\cdot mt}(\cdot),Q_{X_\cdot mt}(\cdot))$ of FBSDEs \eqref{FB:4} satisfying: for $X,X'\in\hr_m$ which are both independent of $\mathcal{W}_t$,
	\begin{align}
		&\sup_{t\le s\le T}\|Y_{X_\cdot mt}(s)\|_{\hr_m}+\sup_{t\le s\le T}\|P_{X_\cdot mt}(s)\|_{\hr_m}+\sum_{j=1}^n\left\|Q^j_{X_\cdot mt}(\cdot)\right\|_{L^2(t,T;\hr_m)} \notag\\
		\le \ & C(l,T,\lambda)\left(1+\|X\|_{\hr_m}\right),  \label{thm3_1}\\
		&\sup_{t\le s\le T}\|Y_{X'_\cdot mt}(s)-Y_{X_\cdot mt}(s)\|_{\hr_m}+\sup_{t\le s\le T}\|P_{X'_\cdot mt}(s)-P_{X_\cdot mt}(s)\|_{\hr_m} \notag\\
		&+\sum_{j=1}^n \left\|Q^j_{X'_\cdot mt}(\cdot)-Q^j_{X_\cdot mt}(\cdot)\right\|_{L^2(t,T;\hr_m)}\le C(l,T,\lambda)\left\|X'- X\right\|_{\hr_m}.  \label{thm3_2}
	\end{align}	
	For any $\tx\in\hr_m$, there is a unique solution $\left({D}_{\tx_\cdot}Y_{X_\cdot mt}(\cdot),\mathcal{D}_{\tx_\cdot}P_{X_\cdot mt},\mathcal{D}_{\tx_\cdot}Q_{X_\cdot mt}\right)$ of the FBSDEs \eqref{FB:7} in $L^2_{\mathcal{W}_{tX\tx}}(t,T;\hr_m)\times L^2_{\mathcal{W}_{tX\tx}}(t,T;\hr_m)\times \left(L^2_{\mathcal{W}_{tX\tx}}(t,T;\hr_m)\right)^n$, such that 
	\begin{align*}
		&\sup_{t\le s\le T}\left\|{D}_{\tx_\cdot}Y_{X_\cdot mt}(s)\right\|_{\hr_m}+\sup_{t\le s\le T}\left\|{D}_{\tx_\cdot}P_{X_\cdot mt}(s)\right\|_{\hr_m}\\
		&+\sum_{j=1}^n\left\|{D}_{\tx_\cdot}Q^j_{X_\cdot mt}(\cdot)\right\|_{L^2(t,T;\hr_m)}\le C(l,T,\lambda)\left\|\tx\right\|_{\hr_m},
	\end{align*}   
	Moreover, $D_{\tx_\cdot}Y_{X_\cdot mt}(s),D_{\tx_\cdot}P_{X_\cdot mt}(s)$ and $D_{\tx_\cdot}Q_{X_\cdot mt}(s)$ are linear in $\tx$, and they are continuous in $X$, and are respectively the G\^ateaux derivatives of $Y_{X_\cdot mt}(s),P_{X_\cdot mt}(s)$ and $Q_{X_\cdot mt}(s))$ in $X$ along the direction $\tx$. Furthermore, if Assumption (B3)-(ii) is satisfied, then, for any $\lambda> 0$, the same assertion remains valid.
\end{theorem}

\section{Bellman Equation}\label{sec:Bellman}

In this section, we apply our results in Section~\ref{sec:value_Hilbert} to study the regularity of the value function $\V$ defined in \eqref{def:V} in $(t,m)\in [0,T]\times\pr_2(\br^n)$ and $X\in \hr_m$ (being independent of $\mathcal{W}_t$), and then provide the solvability of the corresponding Bellman equation for mean field type control problem. From Subsection~\ref{subsec:MP}, by flow property,
\begin{align*}
	\mathcal{V}(X\otimes m,t)=&J_{Xmt}\left(\hvv_{X_\cdot mt}(\cdot)\right)\\
	=&\int_t^T \e\left[\int_\brn g\left(Y_{xmt}(s),Y_{\cdot mt}(s)\#(X\otimes m),\hvv_{xmt}(s),s \right)d(X\otimes m)(x)\right] ds\\
	& +\e\left[\int_\brn g_T(Y_{xmt}(T),Y_{\cdot mt}(T)\#(X\otimes m))d(X\otimes m))(x)\right],
\end{align*}
therefore, $\V$ depends only on $X\otimes m$ and $t$ and it is well-defined. We next give the regularity of $\mathcal{V}$ in $m\in\pr_2(\brn)$ and $X\in \hr_m$, and its proof is provided in Appendix~\ref{pf:thm9}.

\begin{theorem}\label{thm9}
	Under Assumptions (B1), (B2) and (B3)-(i), there exists a constant $c(l,T)$ depending only on $(l,T)$, such that when $\lambda\geq c(l,T)$, the value function $\V$ is twice G\^ateaux differentiable in $X$ with the derivatives
	\begin{align}
		\left.D_X \V(X\otimes m,t)\right|_x=P_{X_x mt}(t),\quad \left.D^2_X \V(X\otimes m,t)\left(\tx\right)\right|_x={D}_{\tx_x}P_{X_x mt}(t),\label{thm9_10}
	\end{align}
	and they satisfy the growth conditions	
	\begin{align}
		|\mathcal{V}(X\otimes m,t)|\le& C(l,T,\lambda)\left(1+\|X\|^2_{\hr_m}\right),\label{lem10_3}\\
		\|D_X \V(X\otimes m,t)\|_{\hr_m}\le& C(l,T,\lambda)(1+\|X\|_{\hr_m}),\label{lem10_1}\\
		\left\|D^2_X \V(X\otimes m,t)\left(\tx\right)\right\|_{\hr_m}\le& C(l,T,\lambda)\left\|\tx\right\|_{\hr_m},\label{lem10_2}
	\end{align}
	and the following continuity condition is valid:
	\begin{equation}
		\begin{split}
			&\left\|D_X \V(X'\otimes m,t)-D_X \V(X\otimes m,t)\right\|_{\hr_m} \le C(l,T,\lambda)\left\|X'-X\right\|_{\hr_m},
		\end{split}
	\end{equation}
	and $D_X^2 \V$ is linear in $\tx$ and continuous in $X$. Moreover, $\V$ is twice linearly functional differentiable in $m$, and the linear functional derivative satisfies the relations
	\begin{align}
		\left.D_X\mathcal{V}(X\otimes m,t)\right|_x=&\ D_\xi \frac{d\mathcal{V}}{d\nu}(X\otimes m,t)(X_x),\label{lem10_4}\\
		\left.D_X^2 \mathcal{V}(X\otimes m,t)\left(\tx\right)\right|_x=&\ \left(D_\xi^2 \frac{d \mathcal{V}}{d\nu}(X\otimes m,t)(X_x)\right)^*\tx_x \notag\\
		&\ +\bar{\e}\left[\int_\brn \left(D_{\xi'}D_\xi \frac{d^2 \mathcal{V}}{d\nu^2}(X\otimes m)\left(X_x,\bx_z\right)\right)^*\bar{\tx}_z dm(z)\right], \label{lem10_5}
	\end{align}
	where $\left(\bx_z,\bar{\tx}_z\right)$ is an independent copy of $\left(X_z,\tx_z\right)$ for $z\in\brn$. Furthermore, if Assumption (B3)-(ii) is satisfied, for any $\lambda> 0$, the same assertion remains valid.
\end{theorem}

We make  the following  assumption on coefficients $(f,\sigma,g)$ in $t$:\\
\textbf{(B4)} The coefficients $\left(f,g\right)$ are continuous in $t\in[0,T]$, and $\sigma^j_3=0$ for $1\le j\le n$. \\
We now give the regularity of $\V$ in time $t$ and the solvability of the Bellman equation, and its proof is given in Appendix~\ref{pf:lem10}.

\begin{theorem}\label{lem10}
	Under Assumptions (B1), (B2), (B3)-(i) and (B4), there exists a constant $c(l,T)$ depending only on $(l,T)$, such that when $\lambda\geq c(l,T)$, the G\^ateaux derivative $D_X\V$ satisfies the following continuity condition in $t$,
	\begin{equation}\label{lem10_8}
		\begin{split}
			&\|D_X \mathcal{V}(X\otimes m,t')-D_X \mathcal{V}(X\otimes m,t)\|_{\hr_m} \le C(l,T,\lambda)\left(1+\|X\|_{\hr_m}\right) |t'-t|^{\frac{1}{2}},
		\end{split}
	\end{equation}	
	and $D_X^2 \mathcal{V}$ is continuous in $t$. The value function $\V$ is $C^1$ in $t$ with the time derivative 
	\begin{align}
		&\frac{\dd\mathcal{V}}{\dd t}(X\otimes m,t)=-\e \bigg[\int_\brn H\bigg(X_x,X\otimes m,t; D_\xi \frac{d\mathcal{V}}{d\nu}(X\otimes m,t)(X_x), \notag\\
		&\quad\qquad\qquad\qquad\qquad\qquad\qquad \frac{1}{2}D_\xi^2 \frac{d \mathcal{V}}{d\nu}(X\otimes m,t)(X_x) \sigma(X_x,X\otimes m,t)\bigg) dm(x)\bigg]. \label{lem10_6}
	\end{align}
	satisfying the growth condition
	\begin{align}
		&\left|\frac{\dd\mathcal{V}}{\dd t}(X\otimes m,t)\right|\le C(l,T,\lambda)\left(1+\|X\|^2_{\hr_m}\right), \label{lem10_7}
	\end{align}
	Moreover, $\mathcal{V}$ is the unique solution (possessing the standard regularity of \eqref{lem10_3}-\eqref{lem10_7}) of the following mean field type Bellman equation: for $ (m,t)\in\pr_2(\brn)\times[0,T]$,
	\begin{equation}\label{Bellman'}
		\left\{
		\begin{aligned}
			&\frac{\dd\mathcal{V}}{\dd t}(m,t)+\int_\brn H\Big(x,m, t; D_\xi \frac{d\mathcal{V}}{d\nu}(m,t)(x),\frac{1}{2}D_\xi^2 \frac{d \mathcal{V}}{d\nu}(m,t)(x) \sigma(x,m,t)\Big) dm(x)=0,\\
			&\mathcal{V}(m,T)= \int_\brn g_T(x, m) dm(x).
		\end{aligned}
		\right.
	\end{equation}	
	Furthermore, if Assumption (B3)-(ii) is satisfied, then, for any $\lambda> 0$, the same assertion remains valid.
\end{theorem}

\section{Mean Field Type Master equation}\label{sec:Master}

We here study the master equation for $\mathcal{U}$ defined as, for $(x,m,t)\in\brn\times\pr_2(\brn)\times[0,T]$,
\begin{align}\label{master_2}
	\mathcal{U}(x,m,t):=\frac{d\V}{d\nu}(m,t)(x)
\end{align}
corresponding to Problem~\ref{intr_4}. We first study the higher regularity of the linear functional derivative of the value function $\V$. From \eqref{thm9_10}, \eqref{lem10_4} and \eqref{master_2}, under assumptions in Theorem~\ref{thm9}, for $t\in[0,T]$,
\begin{align}\label{master_1}
	D_x \mathcal{U}(x,m,t)=D_x\frac{d\V}{d\nu}(m,t)(x)=P_{x mt}(t).
\end{align}
Here, $(Y_{xmt},P_{xmt},Q_{xmt})$ is the unique solution of FBSDEs \eqref{FB:4} with the initial data $X_x:=\mathcal{I}_x=x$, namely, for $(s,x)\in[t,T]\times\brn$,
\small
\begin{align}
	Y_{xmt}(s) =\ & x+\int_t^s f\left(Y_{xmt}(r),Y_{\cdot mt}(r)\otimes m,\hvv_{xmt}(r),r\right)dr \notag\\
	&+\sum_{j=1}^n\int_t^s \sigma^j\left(Y_{xmt}(r),Y_{\cdot mt}(r)\otimes m,\hvv_{xmt}(r),r\right)dw_j(r); \notag\\
	P_{xmt}(s) =\ & -\sum_{j=1}^n \int_s^T Q^j_{xmt}(r)dw_j(r)+D_xg_T(Y_{xmt}(T),Y_{\cdot mt}(T)\otimes m) \notag\\
	&+\bar{\e}\left[\int_\brn D_\xi\frac{d g_T}{d\nu}\left(\bar{Y}_{ymt}(T),Y_{\cdot mt}(T)\otimes m\right)(Y_{xmt}(T)) dm(y)\right] \notag\\
	& +\int_s^T\Bigg\{f_1^*(r)P_{xmt}(r)+\sum_{j=1}^n {\sigma_1^j}^*(r)Q^j_{xmt}(r) +D_x g(Y_{xmt}(r),Y_{\cdot mt}(r)\otimes m,\hvv_{xmt}(r),r) \notag\\
	&\quad\qquad+\bar{\e}\bigg[\int_\brn\bigg( f_2^*(r)\bar{P}_{ymt}(r)+\sum_{j=1}^n {\sigma^j_2}^*(r)\bar{Q}^j_{ymt}(r) \notag\\
	&\qquad\qquad\qquad\qquad +D_\xi \frac{d g}{d\nu}\left(\bar{Y}_{ymt}(r),Y_{\cdot mt}(r)\otimes m,\bar{\hvv}_{ymt}(r),r\right) (Y_{xmt}(r))\bigg) dm(y)\bigg]\Bigg\}dr,  \label{FB:5}
\end{align}
\normalsize
where for $(s,x)\in[t,T]\times\brn$,
\begin{equation}\label{master_7}
	\begin{split}
		\hvv_{xmt}(s)=\hat{v}(Y_{xmt}(s),Y_{\cdot mt}(s)\otimes m,s;P_{xmt}(s),Q_{xmt}(s)),
	\end{split}
\end{equation}
and for all $y\in\brn$, the process $\left(\bar{Y}_{y mt}(s),\bar{P}_{y mt}(s),\bar{Q}_{y mt}(s),\bar{\hvv}_{y mt}(s)\right)$ is an independent copy of process $(Y_{ymt}(s),P_{ymt}(s),Q_{ymt}(s),\hvv_{ymt}(s))$. Therefore, from \eqref{master_1} and FBSDEs \eqref{FB:5}, we know that
\begin{align}\label{master_4}
	D_x^2 \mathcal{U}(x,m,t)=D_x^2\frac{d\V}{d\nu}(m,t)(x)=D_xP_{x mt}(t),
\end{align}
where $(D_xY_{xmt},D_xP_{xmt},D_xQ_{xmt})$ is the unique solution of the following FBSDEs: 
for $(s,x)\in[t,T]\times\brn$,
\begin{align}
	D_xY_{xmt}(s) =\ & \mathcal{I}+\int_t^s [f_1(r)D_xY_{xmt}(r)+f_3(r)D_x\hvv_{xmt}(r)]dr \notag\\
	& +\sum_{j=1}^n\int_t^s \left[\sigma^j_1(r)D_xY_{xmt}(r)+\sigma^j_3(r)D_x\hvv_{xmt}(r)\right]dw_j(r); \notag\\
	D_xP_{xmt}(s) =\ & -\sum_{j=1}^n \int_s^T D_xQ^j_{xmt}(r)dw_j(r)+\left(D_x^2g_T(Y_{xmt}(T),Y_{\cdot mt}(T)\otimes m)\right)^* D_x Y_{xmt}(T) \notag\\
	& +\bar{\e}\left[\int_\brn \left(D^2_\xi\frac{d g_T}{d\nu}\left(\bar{Y}_{ymt}(T),Y_{\cdot mt}(T)\otimes m\right)(Y_{xmt}(T))\right)^* dm(y)\right]D_x Y_{xmt}(T) \notag\\
	& +\int_s^T\Bigg\{f_1^*(r)D_xP_{xmt}(r)+\sum_{j=1}^n {\sigma_1^j}^*(r)D_xQ^j_{xmt}(r) \notag\\
	&\quad\qquad+\left(D_x^2 g\left(Y_{xmt}(r),Y_{\cdot mt}(r)\otimes m,\hvv_{xmt}(r),r\right)\right)^* D_x Y_{xmt}(r) \notag\\
	&\quad\qquad+\left(D_vD_x g\left(Y_{xmt}(r),Y_{\cdot mt}(r)\otimes m,\hvv_{xmt}(r),r\right)\right)^* D_xu_{xmt}(r) \notag\\
	&\quad\qquad+\bar{\e}\bigg[\int_\brn \left(D_\xi^2 \frac{d g}{d\nu}\left(\bar{Y}_{ymt}(r),Y_{\cdot mt}(r)\otimes m,\bar{\hvv}_{ymt}(r),r\right)(Y_{xmt}(r))\right)^* \notag\\
	&\qquad\qquad\qquad\qquad  D_xY_{xmt}(r) dm(y)\bigg] \Bigg\}dr,\label{FB:8}
\end{align}
where for $(s,x)\in[t,T]\times\brn$, and
\begin{align*}
	D_{x}\hvv_{x mt}(s):=\ &(D_x \hat{v}(Y_{x mt}(s),Y_{\cdot mt}(s)\otimes m,s;P_{x mt},Q_{x mt}))^*D_{x}Y_{x mt}(s)\\
	&+(D_p \hat{v}(Y_{x mt}(s),Y_{\cdot mt}(s)\otimes m,s;P_{x mt},Q_{x mt}))^*D_{x}P_{x mt}(s)\\
	& +\sum_{j=1}^n (D_{q^j} \hat{v}(Y_{x mt}(s),Y_{\cdot mt}(s)\otimes m,s;P_{x mt},Q_{x mt}))^*D_{x}Q^j_{x mt}(s).
\end{align*}

\begin{lemma}
	Under Assumptions (B1), (B2) and (B3)-(i), there exists a constant $c(l,T)$ depending only on $(l,T)$, such that when $\lambda\geq c(l,T)$, there is a unique $\mathcal{W}_{t\mathcal{I}}$-adapted solution of FBSDEs \eqref{FB:8} such that
	\begin{equation}\label{master_3}
		\begin{split}
			&\sup_{t\le s\le T}\e\left[\int_\brn |D_x Y_{xmt}(s)|^2 dm(x)\right]+\sup_{t\le s\le T}\e\left[\int_\brn |D_x P_{xmt}(s)|^2 dm(x)\right]\\
			&+\sum_{j=1}^n\int_t^T\e\left[\int_\brn |D_x Q^j_{xmt}(s)|^2 dm(x)\right]ds\le C(l,T,\lambda).
		\end{split}
	\end{equation}
	Moreover, $D_{x}Y_{x mt}(s),D_{x}P_{x mt}(s)$ and $D_{x}Q_{x mt}(s)$ are the respective derivatives in $x$ of $Y_{x mt}(s),P_{x mt}(s)$ and $Q_{x mt}(s)$. Furthermore, if Assumption (B3)-(ii) is satisfied, then, for any $\lambda> 0$, the same assertion remains valid.
\end{lemma}
Its proof is similar to that of Theorem~\ref{thm3}, and we omit it here. In order to warrant the linear differentiability of $\V$ with derivatives \eqref{master_1} and \eqref{master_4}, we need the following regularity-enhanced version of Assumption (B2):\\
\textbf{(B2')} The functionals $g$ and $g_T$ satisfy (B2). The following derivatives exist, and they are continuous in all their arguments and are bounded by $l$:
\begin{align*}
	&D_x^3 g, \ D_vD_x^2 g, \ D_v^2D_x g,\ D_v^3g, \ D_x^2 D_\xi \frac{dg}{d\nu},\ D_xD_v D_\xi \frac{dg}{d\nu}, \\
	&D_v^2 D_\xi \frac{dg}{d\nu},\ D_x D_\xi^2 \frac{dg}{d\nu},\ D_v D_\xi^2 \frac{dg}{d\nu},\ D_\xi^3 \frac{dg}{d\nu}, \  D_{\xi'}D_\xi^2\frac{d^2g}{d\nu^2},\\
	&D_x^3 g_T, \ D_x^2 D_\xi \frac{dg_T}{d\nu},\ D_x D_\xi^2 \frac{dg_T}{d\nu},\ D_\xi^3 \frac{dg_T}{d\nu}, \  D_{\xi'}D_\xi^2\frac{d^2g_T}{d\nu^2};
\end{align*}
and for $(x,m,v,s,\xi,\xi')\in\brn\times\pr_2(\brn)\times\brd\times[0,T]\times\brn\times\brn$,
\begin{align*}
	&\left|\left(D_x^2\frac{dg}{d\nu},D_vD_x\frac{dg}{d\nu},D_v^2\frac{dg}{d\nu}\right)(x,m,v,s)(\xi)\right|,\ \left|D_x^2\frac{dg_T}{d\nu}(x,m)(\xi)\right|\le l(1+|\xi|);\\
	&\left|D_\xi^2\frac{d^2g}{d\nu^2}(x,m,v,s)(\xi,\xi')\right|,\ \left|D_\xi^2\frac{d^2g_T}{d\nu^2}(x,m)(\xi,\xi')\right|\le l(1+|\xi'|).
\end{align*}

From Assumption (B2') and Estimate \eqref{master_3}, the linear functional derivatives of the processes $\left(Y_{ymt}(s),P_{ymt}(s),Q_{ymt}(s),\hvv_{ymt}(s)\right)$ and $\left(D_{x}Y_{x mt}(s),D_{x}P_{x mt}(s),D_{x}Q_{x mt}(s)\right)$ in $m\in\pr_2(\brn)$ can be characterized as the solutions of FBSDEs. These two corresponding FBSDEs can be checked out in Appendix~\ref{pf:thm10}.

\begin{theorem}\label{thm:10}
	Under Assumptions (B1), (B2') and (B3)-(i), there exists a constant $c(l,T)$ depending only on $(l,T)$, such that when $\lambda\geq c(l,T)$, there is a unique solution
	\begin{align*}
		\left(\frac{dY_{xmt}}{d\nu}(\xi,s),\frac{dP_{xmt}}{d\nu}(\xi,s),\frac{dQ_{xmt}}{d\nu}(\xi,s)\right)
	\end{align*}	
	of FBSDEs \eqref{FB:9} (the unlifted version of FBSDEs \eqref{FB:7}) and a unique solution
	\begin{align*}
		\left(\frac{dD_xY_{xmt}}{d\nu}(\xi,s),\frac{dD_xP_{xmt}}{d\nu}(\xi,s),\frac{dD_xQ_{xmt}}{d\nu}(\xi,s)\right)
	\end{align*}
	of FBSDEs \eqref{FB:10} (the Jacobian flow of FBSDEs \eqref{FB:9}), such that for $t\in[0,T]$,
	\begin{align*}
		&(i)\ \sup_{t\le s\le T}\e\left[\int_\brn \left|\frac{dY_{xmt}}{d\nu}(\xi,s)\right|^2 dm(x)\right]+\sup_{t\le s\le T}\e\left[\int_\brn \left|\frac{dP_{xmt}}{d\nu}(\xi,s)\right|^2dm(x)\right]\\
		&\qquad+\sum_{j=1}^n\int_t^T\e\left[\int_\brn\left|\frac{dQ^j_{xmt}}{d\nu}(\xi,s)\right|^2 dm(x)\right]ds\le C(l,T,\lambda)\left(1+|\xi|^2\right);\\
		&(ii)\ \sup_{t\le s\le T}\e\left[\int_\brn \left|\frac{dD_xY_{xmt}}{d\nu}(\xi,s)\right|^2 dm(x)\right]+\sup_{t\le s\le T}\e\left[\int_\brn \left|\frac{dD_xP_{xmt}}{d\nu}(\xi,s)\right|^2dm(x)\right]\\
		&\qquad+\sum_{j=1}^n\int_t^T\e\left[\int_\brn\left|\frac{dD_xQ^j_{xmt}}{d\nu}(\xi,s)\right|^2 dm(x)\right]ds\le C(l,T,\lambda)\left(1+|\xi|^2\right).
	\end{align*}      
	Moreover, processes $\frac{dY_{xmt}}{d\nu}(\xi,s)$, $\frac{dP_{xmt}}{d\nu}(\xi,s)$, $\frac{dQ_{xmt}}{d\nu}(\xi,s)$,  $\frac{dD_xY_{xmt}}{d\nu}(\xi,s)$, $\frac{dD_xP_{xmt}}{d\nu}(\xi,s)$ and $\frac{dD_xQ_{xmt}}{d\nu}(\xi,s)$ are the linear functional derivatives of processes $Y_{x mt}(s)$, $P_{x mt}(s)$, $Q_{x mt}(s)$, $D_xY_{x mt}(s)$, $D_xP_{x mt}(s)$ and $D_xQ_{x mt}(s))$, respectively. Furthermore, if Assumption (B3)-(ii) is satisfied, then, for any $\lambda> 0$, the same assertion remains valid.
\end{theorem}

We finally give the well-posedness of the mean field type master equation, and its proof is given in Appendix~\ref{pf:thm11}.

\begin{theorem}\label{thm:11}
	Under Assumptions (B1), (B2'), (B3)-(i) and (B4), there exists a constant $c(l,T)$ depending only on $(l,T)$, such that when $\lambda\geq c(l,T)$, for $(t,m,x,\xi)\in[0,T]\times\pr_2(\brn)\times\brn\times\brn$, the derivatives
	\begin{equation}\label{master_6}
		\begin{split}
			&D_x \mathcal{U}(x,m,t)=P_{x mt}(t),\quad D_x^2 \mathcal{U}(x,m,t)=D_xP_{x mt}(t),\\
			&D_{x}\frac{d \mathcal{U}}{d\nu}(x,m,t)(\xi)=\frac{dP_{x mt}}{d\nu}(\xi,t),\quad D_{x}^2\frac{d \mathcal{U}}{d\nu}(x,m,t)(\xi)=\frac{dD_x P_{x mt}}{d\nu}(\xi,t),
		\end{split}
	\end{equation}
	satisfy
	\begin{equation}\label{thm11_0}
		\begin{split}
			&\int_\brn \left|D_x \mathcal{U}(x,m,t)\right|^2 dm(x)\le C(l,T,\lambda)\left(1+|x|^2\right),\  \int_\brn \left|D_x^2 \mathcal{U}(x,m,t)\right|^2 dm(x)\le C(l,T,\lambda),\\
			&\int_\brn \left|D_{x}\frac{d\mathcal{U}}{d\nu}(x,m,t)(\xi)\right|^2dm(x),\ \int_\brn \left|D_{x}^2\frac{d \mathcal{U}}{d\nu}(x,m,t)(\xi)\right|^2dm(x)\le C(l,T,\lambda)\left(1+|\xi|^2\right).
		\end{split}
	\end{equation}
	Moreover, $\mathcal{U}$ is differentiable in $t$, and is the unique classical solution of the following mean field type master equation of: for $ (t,m,x)\in[0,T]\times\pr_2(\brn)\times\brn$,
	\begin{equation}\label{master}
		\left\{
		\begin{aligned}
			&\frac{\dd \mathcal{U}}{\dd t}(x,m,t)+H\bigg(x,m,t; D_x \mathcal{U}(x,m,t),\frac{1}{2}D_x^2 \mathcal{U}(x,m,t)(t) \sigma(x,m,t)\bigg) \\
			&+\int_\brn\bigg\{ \left(D_p H(\xi,m,t;D_\xi \mathcal{U}(\xi,m,t),\frac{1}{2}D_\xi^2 \mathcal{U}(\xi,m,t) \sigma(\xi,m,t))\right)^*D_{\xi} \frac{d\mathcal{U}}{d\nu}(\xi,m,t)(x) \\
			&\qquad\qquad +\frac{1}{2}\text{Tr}\left[\sigma\sigma^*(\xi,m,t)D_{\xi}^2\frac{d \mathcal{U}}{d\nu}(\xi,m,t)(x)\right]\\
			&\qquad\qquad +\frac{dH}{d\nu}(\xi,m,t;D_\xi \mathcal{U}(\xi,m,t),\frac{1}{2}D_\xi^2 \mathcal{U}(\xi,m,t) \sigma(\xi,m,t))(x)\bigg\}dm(\xi)=0,\\
			&\mathcal{U}(x,m,T)=g_T(x, m)+ \int_\brn \frac{dg_T}{d\nu}(\xi, m)(x) dm(\xi).
		\end{aligned}
		\right.
	\end{equation}	
	Furthermore, if Assumption (B3)-(ii) is satisfied, then, for any $\lambda> 0$, the same assertion remains valid.
\end{theorem}

\section*{Acknowledgments}

Alain Bensoussan is supported by the National Science Foundation under grant NSF-DMS-2204795. This work also constitutes part of Ziyu Huang's Ph.D. dissertation at Fudan University, China. Shanjian Tang is supported by the National Natural Science Foundation of China under grant nos. 11631004 and 12031009. Phillip Yam acknowledges the financial supports from HKGRF-14301321 with the project title ``General Theory for Infinite Dimensional Stochastic Control: Mean Field and Some Classical Problems'', and HKGRF-14300123 with the project title ``Well-posedness of Some Poisson-driven Mean Field Learning Models and their Applications''.

\appendix

\section{Proof of Statements in Section~\ref{sec:FB_Hilbert}}\label{app:01}

\subsection{Proof of Lemma~\ref{thm:1}}\label{pf:thm:1}

The existence and uniqueness of $(Y_{Xt},P_{Xt},Q_{Xt})$ can be proven with the method of continuation in coefficients of Hu-Peng \cite{YH2}, similar to the proof of Peng-Wu \cite[Theorem 2.3]{SP} (for FBSDEs in Euclidean spaces) and  Ahuja-Ren-Yang \cite[Theorem 1]{SA1} for FBSDEs in an Hilbert space, and is thus omitted. Here, we only prove \eqref{thm1_1}. We set $\de{Y}(s):={Y}_{X't}(s)-{Y}_{Xt}(s)$, $\de{P}(s):={P}_{X't}(s)-{P}_{Xt}(s)$, $\de{Q}(s):={Q}_{X't}(s)-{Q}_{Xt}(s)$ and $ \Delta X:=X'-X$. For $s\in[t,T]$, we set $\beta_{Xt}(s):=\beta(Y_{Xt}(s),s;P_{Xt}(s),Q_{Xt}(s))$, $\beta_{X't}(s):=\beta(Y_{X't}(s),s;P_{X't}(s),Q_{X't}(s))$ and denote by $\de \beta(s):=\beta_{X't}(s)-\beta_{Xt}(s)$. Then $(\Delta Y,\Delta P,\Delta Q)$ satisfy the following FBSDEs: for $s\in[t,T]$,
\begin{equation}\label{thm1_4}
	\begin{aligned}
		&\de Y(s)=\de X+\int_t^s \left[\mathbb{F} (Y_{X't}(r),r;P_{X't}(r),Q_{X't}(r))-\mathbb{F} (Y_{Xt}(r),r;P_{Xt}(r),Q_{Xt}(r))\right]dr\\
		&\qquad\quad\ \  +\sum_{j=1}^n\int_t^s  [\mathbb{A}^j(Y_{X't}(r),r;P_{X't}(r),Q_{X't}(r))-\mathbb{A}^j(Y_{Xt}(r),r;P_{Xt}(r),Q_{Xt}(r))]dw_j(r),\\
		&\de P(s)=\mathbb{G}_T(Y_{X't}(T))-\mathbb{G}_T(Y_{Xt}(T))-\sum_{j=1}^n\int_s^T \de Q^j(r)dw_j(r)\\
		&\qquad\qquad+\int_s^T[\mathbb{G}(Y_{X't}(r),r;P_{X't}(r),Q_{X't}(r))-\mathbb{G}(Y_{Xt}(r),r;P_{Xt}(r),Q_{Xt}(r))]dr.
	\end{aligned}
\end{equation}
Using standard arguments of SDEs, from Condition \eqref{bh_2}, we deduce
\begin{equation}\label{thm1_2}
	\sup_{t\le s\le T}\|\de Y(s)\|_{\hr_m}\le C(\mathbb{L},T)\left(\|\de X\|_{\hr_m}+\|\de \beta\|_{L^2(t,T;\ur_m)}\right).
\end{equation}
Meanwhile, by Condition \eqref{bh_1} and Estimate \eqref{thm1_2}, we have	
\begin{equation}\label{thm1_3}
	\begin{split}
		&\sup_{t\le s\le T}\|\de P(s)\|_{\hr_m}+\sum_{j=1}^n\|\de {Q}^j\|_{L^2(t,T;\hr_m)}\le C(\mathbb{L},T)\left(\|\de X\|_{\hr_m}+\|\de \beta\|_{L^2(t,T;\ur_m)}\right).
	\end{split}
\end{equation}
From FBSDEs \eqref{thm1_4} and the monotonicity condition \eqref{monotonicity}, we have for $s\in[t,T]$,
\begin{equation}\label{thm1_8}
	\begin{split}
		&\frac{d}{ds}\left(\de P(s),\de Y(s)\right)_{\hr_m}\\
		=&\left(\de P(s), \mathbb{F} (Y_{X't}(s),s;P_{X't}(s),Q_{X't}(s))-\mathbb{F} (Y_{Xt}(s),s;P_{Xt}(s),Q_{Xt}(s))\right)_{\hr_m}\\
		&+\sum_{j=1}^n \left(\de Q^j(s), \mathbb{A}^j(Y_{X't}(s),s;P_{X't}(s),Q_{X't}(s))-\mathbb{A}^j(Y_{Xt}(s),s;P_{Xt}(s),Q_{Xt}(s))\right)_{\hr_m}\\
		&-\left(\de Y(s), \mathbb{G}(Y_{X't}(s),s;P_{X't}(s),Q_{X't}(s))-\mathbb{G}(Y_{Xt}(s),s;P_{Xt}(s),Q_{Xt}(s))\right)_{\hr_m}\\
		\le& -\Lambda \|\de \beta (s)\|^2_{\ur_m}+\alpha \left(\|\de Y(s)\|^2_{\hr_m}+\|\de P(s)\|^2_{\hr_m}+\sum_{j=1}^n\|\de Q^j(s)\|^2_{\hr_m}\right).
	\end{split}
\end{equation}
Therefore, from Condition \eqref{bh_1} and the Young's inequality with a weight of $\frac{1}{2}$, we have
\begin{equation}\label{thm1_5}
	\begin{split}
		&\Lambda \|\de \beta\|^2_{L^2(t,T;\ur_m)}\\
		\le& \left(\de P(t),\de X\right)_{\hr_m}-\left(D_X \mathbb{G}_T(Y_{X't}(T))-D_X \mathbb{G}_T(Y_{Xt}(T)),\de Y(T)\right)_{\hr_m}\\
		& +\alpha\int_t^T  \left(\|\de Y(s)\|^2_{\hr_m}+\|\de P(s)\|^2_{\hr_m}+\sum_{j=1}^n\|\de Q^j(s)\|^2_{\hr_m} \right)ds\\
		\le& \  C(\alpha,\mathbb{L},T) \left(\sup_{t\le s\le T}\|\de Y(s)\|^2_{\hr_m}+\sup_{t\le s\le T}\|\de P(s)\|^2_{\hr_m}+\sum_{j=1}^n\|\de {Q}^j\|^2_{L^2(t,T;\hr_m)}\right).
	\end{split}
\end{equation}
We then substitute \eqref{thm1_2} and \eqref{thm1_3} into \eqref{thm1_5}, we can deduce 
\begin{align*}
	\Lambda \|\de \beta (s)\|^2_{L^2(t,T;\ur_m)} \le C(\alpha,\mathbb{L},T) \left( \|\de X\|^2_{\hr_m}+\|\de \beta\|^2_{L^2(t,T;\ur_m)}\right).
\end{align*}
So there exists a constant $c(\alpha,\mathbb{L},T)$, such that we have  for $\Lambda \geq c(\alpha,\mathbb{L},T)$, $\|\de \beta\|_{L^2(t,T;\ur_m)} \le C(\alpha,\mathbb{L},T,\Lambda) \|\de X\|_{\hr_m}$, which is substituted back to \eqref{thm1_2} and \eqref{thm1_3}, we obtain \eqref{thm1_1}. We now consider the case when the parameter $\alpha=0$ and $\mathbb{G}_T$ satisfy the monotonicity condition \eqref{thm1_7}. In this case, from \eqref{thm1_8}, we have	$\frac{d}{ds}\left(\de P(s),\de Y(s)\right)_{\hr_m}\le -\Lambda \|\de \beta (s)\|^2_{\ur_m}$. Therefore, from Condition \eqref{thm1_7} and the Young's inequality with a weight of $\frac{1}{2\sqrt{\epsilon}}$, for any $\epsilon>0$,
\begin{equation}\label{thm1_9}
	\begin{split}
		\Lambda  \|\de \beta \|^2_{L^2(t,T;\ur_m)} &\le \left(\de P(t),\de X\right)_{\hr_m}\le \frac{1}{4\epsilon}\|\de X\|^2_{\hr_m}+\epsilon\sup_{t\le s\le T}\|\de P(s)\|^2_{\hr_m}.
	\end{split}
\end{equation}
We then substitute \eqref{thm1_3} into \eqref{thm1_9}, we can deduce
\begin{align*}
	\Lambda \|\de \beta \|^2_{L^2(t,T;\ur_m)} \le \frac{1}{4\epsilon}\|\de X\|^2_{\hr_m}+\epsilon C(\mathbb{L},T)\left(\|\de X\|^2_{\hr_m}+\|\de \beta\|^2_{L^2(t,T;\ur_m)}\right).
\end{align*}
We then obtain \eqref{thm1_1} by choosing $\epsilon$ small.

\subsection{Proof of Theorem~\ref{thm2}}\label{pf:thm2}

In view of Lemma~\ref{thm:1}, we only need to check the monotonicity condition \eqref{monotonicity} and the Lipschitz-continuity conditions \eqref{bh_2}-\eqref{bh_1} under our settings. From Assumption (A1), \eqref{hv} and \eqref{hv2}, for $s\in[t,T]$, for any $X,X',P,P'\in\hr_m$ and $Q,Q'\in\hr_m^n$,
\begin{align}
	&\left(D_X \bh(X',s;P',Q')-D_X \bh(X,s;P,Q),X'-X\right)_{\hr_m} \notag\\
	& -\left(D_P \bh(X',s;P',Q')-D_P \bh(X,s;P,Q),P'-P\right)_{\hr_m} \notag\\
	& -\sum_{j=1}^n \left(D_{Q^j} \bh(X',s;P',Q')-D_{Q^j} \bh(X,s;P,Q),Q'^{j}-Q^j\right)_{\hr_m} \notag\\
	=&\left(D_X G\left(X',\hv(X',s;P',Q'),s\right)-D_X G\left(X,\hv(X,s;P,Q),s\right),\ X'-X\right)_{\hr_m} \notag\\
	&-\left(\mathcal{F}_2(s)\left(\hv(X',s;P',Q')-\hv(X,s;P,Q)\right),P'-P\right)_{\hr_m} \notag\\
	& -\sum_{j=1}^n \left(\mathcal{A}^j_2(s)\left(\hv(X',s;P',Q')-\hv(X,s;P,Q)\right),Q'^{j}-Q^j\right)_{\hr_m} \notag\\
	=&\left(D_X G\left(X',\hv(X',s;P',Q'),s\right)-D_X G\left(X,\hv(X,s;P,Q),s\right),\ X'-X\right)_{\hr_m} \notag\\
	& +\Big(D_V G\left(X',\hv(X',s;P',Q'),s\right)-D_V G\left(X,\hv(X,s;P,Q),s\right),\notag\\
	&\qquad \hv(X',s;P',Q')-\hv(X,s;P,Q) \Big)_{\ur_m}. \label{mono_1}
\end{align}
From Assumption (A3)-(i), we have 
\begin{equation}\label{mono_2}
	\begin{split}
		&\Big(D_V G\left(X,\hv(X',s;P',Q'),s\right)-D_V G\left(X,\hv(X,s;P,Q),s\right),\\
		&\quad \hv(X',s;P',Q')-\hv(X,s;P,Q) \Big)_{\ur_m}\geq  2\lambda \left\|\hv(X',s;P',Q')-\hv(X,s;P,Q)\right\|^2_{\ur_m}. 
	\end{split}
\end{equation}
From Assumption (A2) and the Young's inequality with a weight of $\frac{1}{2\sqrt{\lambda}}$, we have
\begin{equation}\label{mono_3}
	\begin{split}
		&\left(D_X G\left(X',\hv(X',s;P',Q'),s\right)-D_X G\left(X,\hv(X,s;P,Q),s\right),\ X'-X\right)_{\hr_m}\\
		& +\Big(D_V G\left(X',\hv(X',s;P',Q'),s\right)-D_V G\left(X,\hv(X',s;P',Q'),s\right),\\
		&\qquad \hv(X',s;P',Q')-\hv(X,s;P,Q)\Big)_{\ur_m}\\
		\geq& -C(L)\left\|X'-X\right\|_{\hr_m}\left(\left\|X'-X\right\|_{\hr_m}+\left\|\hv(X',s;P',Q')-\hv(X,s;P,Q)\right\|_{\ur_m}\right)\\
		\geq& -\lambda \left\|\hv(X',s;P',Q')-\hv(X,s;P,Q)\right\|^2_{\ur_m} -C(L)\left(1+\frac{1}{\lambda}\right) \left\|X'-X\right\|^2_{\hr_m}.
	\end{split}
\end{equation}
From \eqref{mono_2} and \eqref{mono_3}, we have
\begin{equation}\label{mono_10}
	\begin{split}
		&\left(D_X G\left(X',\hv(X',s;P',Q'),s\right)-D_X G\left(X,\hv(X,s;P,Q),s\right),\ X'-X\right)_{\hr_m} \\
		& +\Big(D_V G\left(X',\hv(X',s;P',Q'),s\right)-D_V G\left(X,\hv(X,s;P,Q),s\right),\\
		&\qquad \hv(X',s;P',Q')-\hv(X,s;P,Q) \Big)_{\ur_m}\\
		\geq&\  \lambda \left\|\hv(X',s;P',Q')-\hv(X,s;P,Q)\right\|^2_{\ur_m} -C(L)\left(1+\frac{1}{\lambda}\right) \left\|X'-X\right\|^2_{\hr_m}.
	\end{split}
\end{equation}	
Substituting \eqref{mono_10} into \eqref{mono_1}, the monotonicity condition \eqref{monotonicity} is satisfied with the map $\beta(\cdot)=\hv(\cdot)$ by setting $\Lambda=\lambda$ and $\alpha=C(L)\left(1+\frac{1}{\lambda}\right)$. From \eqref{hv2} and \eqref{hv1}, we know that the Lipschitz-continuity conditions \eqref{bh_2}-\eqref{bh_1} are also satisfied with $\mathbb{L}=C(L)$. Hence, if Assumption (A3)-(ii) is satisfied, the right hand side of \eqref{mono_1}:
\begin{align*}
	&\left(D_X G\left(X',\hv(X',s;P',Q'),s\right)-D_X G\left(X,\hv(X,s;P,Q),s\right),\ X'-X\right)_{\hr_m} \notag\\
	& +\Big(D_V G\left(X',\hv(X',s;P',Q'),s\right)-D_V G\left(X,\hv(X,s;P,Q),s\right),\notag\\
	&\qquad \hv(X',s;P',Q')-\hv(X,s;P,Q) \Big)_{\ur_m}\geq  \ 2\lambda \left\|\hv(X',s;P',Q')-\hv(X,s;P,Q)\right\|^2_{\ur_m},
\end{align*}
that is, the case is valid for $\alpha=0$. As a consequence of Lemma~\ref{thm:1}, we obtain the solvability for FBSDEs \eqref{FB:1}, we conclude with Estimates \eqref{thm2_1} and \eqref{thm2_2}.

\subsection{Proof of Lemma~\ref{thm4}}\label{pf:thm4}

Note the fact that $\left(Y_{Xt},P_{Xt},Q_{Xt},\hat{V}_{Xt}\right)$ is already known, and \eqref{FB:6} is a system of FBSDEs for $\left(\mathcal{D}_{\tx}P_{Xt},\mathcal{D}_{\tx}Q_{Xt},\mathcal{D}_{\tx}Y_{Xt}\right)$. We first use Lemma~\ref{thm:1} to prove the solvability. We define the map $\hat{Z}:[0,T]\times\hr_m\times\hr_m\times\hr_m^n\to\ur_m$ as
\begin{align*}
	\hat{Z}\left(s;\ty,\tp,\tq\right):=&D_X\hv(Y_{Xt}(s),s;P_{Xt}(s),Q_{Xt}(s))\left(\ty\right)+D_P\hv(Y_{Xt}(s),s;P_{Xt}(s),Q_{Xt}(s))\left(\tp\right)\\
	&+\sum_{j=1}^nD_{Q^j}\hv(Y_{Xt}(s),s;P_{Xt}(s),Q_{Xt}(s))\left(\tq^j\right).
\end{align*}
We now compute the monotonicity condition for FBSDEs \eqref{FB:6}. For $s\in[t,T]$, $\ty,\ty'\tp,\tp'\in\hr_m$ and $\tq,\tq'\in\hr^n_m$, we have
\small
\begin{align}
	\de\mathbb{I}(s) :=&-\bigg(\mathcal{F}^*_1(s)\left(\tp'-\tp\right)+\sum_{j=1}^n {\mathcal{A}^j_1}^*(s)\left(\tq'^{j}-\tq^j\right)+D_X^2 G\left(Y_{Xt}(s),\hat{V}_{Xt}(s),s\right)\left(\ty'-\ty\right) \notag\\
	&\qquad +D_V D_X G\left(Y_{Xt}(s),\hv_{Xt}(s),s\right)\left(\hat{Z}\left(s,\ty',\tp',\tq'\right)-\hat{Z}\left(s,\ty,\tp,\tq\right)\right),\quad \ty'-\ty\bigg)_{\hr_m} \notag\\
	&+\left(\mathcal{F}_1(s)\left(\ty'-\ty\right)+\mathcal{F}_2(s)\left(\hat{Z}\left(s,\ty',\tp',\tq'\right)-\hat{Z}\left(s,\ty,\tp,\tq\right)\right), \quad \tp'-\tp\right)_{\hr_m} \notag\\
	&+\sum_{j=1}^n\left(\mathcal{A}^j_1(s)\left(\ty'-\ty\right)+\mathcal{A}^j_2(s)\left(\hat{Z}\left(s,\ty',\tp',\tq'\right)-\hat{Z}\left(s,\ty,\tp,\tq\right)\right), \quad \tq'^{j}-\tq^j\right)_{\hr_m} \notag\\
	=&-\left(D_X^2 G\left(Y_{Xt}(s),\hv_{Xt}(s),s\right)\left(\ty'-\ty\right),\quad \ty'-\ty\right)_{\hr_m} \notag\\
	& -\left(\hat{Z}\left(s,\ty',\tp',\tq'\right)-\hat{Z}\left(s,\ty,\tp,\tq\right),\quad D_X D_V G\left(Y_{Xt}(s),\hv_{Xt}(s),s\right)\left(\ty'-\ty\right)\right)_{\ur_m} \notag\\
	&+\left(\hat{Z}\left(s,\ty',\tp',\tq'\right)-\hat{Z}\left(s,\ty,\tp,\tq\right), \quad \mathcal{F}^*_2(s)\left(\tp'-\tp\right)\right)_{\ur_m} \notag\\
	&+\sum_{j=1}^n\left(\hat{Z}\left(s,\ty',\tp',\tq'\right)-\hat{Z}\left(s,\ty,\tp,\tq\right), \quad {\mathcal{A}^j_2}^*(s)\left(\tq'^{j}-\tq^j\right)\right)_{\ur_m}. \label{thm4_1}
\end{align}
\normalsize
From Condition \eqref{optimal_condition}, we have	  
\begin{equation}\label{thm4_2}
	\begin{split}
		&D_XD_V G\left(Y_{Xt}(s),\hv_{Xt}(s),s\right)\left(\ty'-\ty\right)\\
		=&-D^2_V G\left(Y_{Xt}(s),\hv_{Xt}(s),s\right)\Big(D_X\hv(Y_{Xt}(s),s;P_{Xt}(s),Q_{Xt}(s))\left(\ty'\right)\\
		&\quad\qquad\qquad\qquad\qquad\qquad\qquad -D_X\hv(Y_{Xt}(s),s;P_{Xt}(s),Q_{Xt}(s))\left(\ty\right)\Big);\\
		&\mathcal{F}^*_2(s)\left(\tp'-\tp\right)\\
		=&-D^2_V G\left(Y_{Xt}(s),\hv_{Xt}(s),s\right)\Big(D_P\hv(Y_{Xt}(s),s;P_{Xt}(s),Q_{Xt}(s))\left(\tp'\right)\\
		&\quad\qquad\qquad\qquad\qquad\qquad\qquad -D_P\hv(Y_{Xt}(s),s;P_{Xt}(s),Q_{Xt}(s))\left(\tp\right)\Big);\\
		&{\mathcal{A}^j_2}^*(s)\left(\tq'^{j}-\tq^j\right)\\
		=&-D^2_V G\left(Y_{Xt}(s),\hv_{Xt}(s),s\right)\Big(D_{Q^j}\hv(Y_{Xt}(s),s;P_{Xt}(s),Q_{Xt}(s))\left(\tq'^{j}\right)\\
		&\quad\qquad\qquad\qquad\qquad\qquad\qquad -D_{Q^j}\hv(Y_{Xt}(s),s;P_{Xt}(s),Q_{Xt}(s))\left(\tq^j\right)\Big).
	\end{split}
\end{equation}
Substituting \eqref{thm4_2} into \eqref{thm4_1}, we have
\small
\begin{align}
	\de\mathbb{I}(s)=&-\left(D_X^2 G\left(Y_{Xt}(s),\hv_{Xt}(s),s\right)\left(\ty'-\ty\right),\  \ty'-\ty\right)_{\hr_m} \notag \\
	& -2\Big(D_X D_V G\left(Y_{Xt}(s),\hv_{Xt}(s),s\right)\left(\ty'-\ty\right),\ \hat{Z}\left(s,\ty',\tp',\tq'\right)-\hat{Z}\left(s,\ty,\tp,\tq\right)\Big)_{\ur_m} \notag\\
	& -\Big(D^2_V G\left(Y_{Xt}(s),\hv_{Xt}(s),s\right)\left(\hat{Z}\left(s,\ty',\tp',\tq'\right)-\hat{Z}\left(s,\ty,\tp,\tq\right) \right), \notag\\
	&\ \quad\qquad\qquad\qquad\qquad\qquad\qquad \hat{Z}\left(s,\ty',\tp',\tq'\right)-\hat{Z}\left(s,\ty,\tp,\tq\right)\Big)_{\ur_m}. \label{thm4_3}
\end{align}
\normalsize
From \eqref{mono_11}, we know that
\begin{align*}
	\de \mathbb{I}(s)\le -\lambda \left\|\hat{Z}\left(s,\ty',\tp',\tq'\right)-\hat{Z}\left(s,\ty,\tp,\tq\right)\right\|^2_{\ur_m}-C(L)\left(1+\frac{1}{\lambda}\right)\left\|\ty'-\ty\right\|^2_{\hr_m}.
\end{align*}
Therefore, FBSDEs \eqref{FB:6} satisfy the monotonicity condition \eqref{monotonicity} with the map $\beta(\cdot)=\hat{Z}(\cdot)$ by setting $\Lambda=\lambda$ and $\alpha=C(L)\left(1+\frac{1}{\lambda}\right)$. From Assumptions (A1) and (A2), we know that 
\begin{align*}
	&\left\|\mathcal{F}_1(s)\left(\ty'-\ty\right)+\mathcal{F}_2(s)\left(\hat{Z}(s,\ty',\tp',\tq')-\hat{Z}(s,\ty,\tp,\tq)\right)\right\|_{\hr_m}\\
	&\ +\sum_{j=1}^n\left\|\mathcal{A}^j_1(s)\left(\ty'-\ty\right)+\mathcal{A}^j_2(s)\left(\hat{Z}\left(s,\ty',\tp',\tq'\right)-\hat{Z}\left(s,\ty,\tp,\tq\right)\right)\right\|_{\hr_m} \\
	\le& \ C(L) \left(\left\|\ty'-\ty\right\|_{\hr_m}+\left\| \hat{Z}\left(s,\ty',\tp',\tq'\right)-\hat{Z}\left(s,\ty,\tp,\tq\right)\right\|_{\ur_m}\right);\\
	&\bigg\|\mathcal{F}^*_1(s)\left(\tp'-\tp\right)+\sum_{j=1}^n {\mathcal{A}^j_1}^*(s)\left(\tq'^{j}-\tq^j\right)+D_X^2 G(Y_{Xt}(s),\hv_{Xt}(s),s)\left(\ty'-\ty\right) \\
	& \ +D_V D_X G\left(Y_{Xt}(s),\hv_{Xt}(s),s\right)\left(\hat{Z}\left(s,\ty',\tp',\tq'\right)-\hat{Z}\left(s,\ty,\tp,\tq\right)\right)\bigg\|_{\hr_m}\\
	\le&\  C(L)\bigg(\left\|\ty'-\ty\right\|_{\hr_m}+\left\|\tp'-\tp\right\|_{\hr_m}+\sum_{j=1}^n \left\|\tq'^{j}-\tq^j\right\|_{\hr_m}\\
	&\ \quad\qquad+\left\| \hat{Z}\left(s,\ty',\tp',\tq'\right)-\hat{Z}\left(s,\ty,\tp,\tq\right)\right\|_{\ur_m}\bigg);\\
	&\left\|D^2_X G_T(Y_{Xt}(T))\left(\ty'-\ty\right)\right\|_{\hr_m}\le L \left\|\ty'-\ty\right\|_{\hr_m}.
\end{align*}
That is, the coefficients of FBSDEs \eqref{FB:6} satisfy the Lipschitz-continuity conditions \eqref{bh_2}-\eqref{bh_1} with $\mathbb{L}=C(L)$. Furthermore, if Assumption (A3)-(ii) is satisfied, then, from \eqref{thm4_3} and Condition \eqref{mono_12}, we have the following monotonicity    
\begin{align*}
	&\de \mathbb{I}(s)\le -2\lambda \left\|\hat{Z}\left(s,\ty',\tp',\tq'\right)-\hat{Z}\left(s,\ty,\tp,\tq\right)\right\|^2_{\ur_m},
\end{align*}
that is, the case is valid for $\alpha=0$. Therefore, from Lemma~\ref{thm:1}, we obtain the solvability of FBSDEs \eqref{FB:6}. Next, we prove \eqref{thm4_0}. We can compute that, for $s\in(t,T)$,
\begin{align}
	&\frac{d}{ds}\left(\mathcal{D}_{\tx}Y_{Xt}(s),\mathcal{D}_{\tx}P_{Xt}(s)\right)_{\hr_m} \notag\\
	=&-\left(D_X^2 G\left(Y_{Xt}(s),\hv_{Xt}(s),s\right)\left(\mathcal{D}_{\tx}Y_{Xt}(s)\right),\  \mathcal{D}_{\tx}Y_{Xt}(s)\right)_{\hr_m} \notag\\
	& -2\Big(D_X D_V G\left(Y_{Xt}(s),\hv_{Xt}(s),s\right)\left(\mathcal{D}_{\tx}Y_{Xt}(s)\right), \ \hat{Z}\left(s,\mathcal{D}_{\tx}Y_{Xt}(s),\mathcal{D}_{\tx}P_{Xt}(s),\mathcal{D}_{\tx}Q_{Xt}(s)\right)\Big)_{\ur_m} \notag\\
	& -\Big(D^2_V G\left(Y_{Xt}(s),\hv_{Xt}(s),s\right)\left(\hat{Z}\left(s,\mathcal{D}_{\tx}Y_{Xt}(s),\mathcal{D}_{\tx}P_{Xt}(s),\mathcal{D}_{\tx}Q_{Xt}(s)\right)\right),\notag\\
	&\qquad\qquad\qquad\qquad\qquad\qquad\qquad \hat{Z}\left(s,\mathcal{D}_{\tx}Y_{Xt}(s),\mathcal{D}_{\tx}P_{Xt}(s),\mathcal{D}_{\tx}Q_{Xt}(s)\right)\Big)_{\ur_m} \notag\\
	\le& -\lambda \left\|\hat{Z}\left(s,\mathcal{D}_{\tx}Y_{Xt}(s),\mathcal{D}_{\tx}P_{Xt}(s),\mathcal{D}_{\tx}Q_{Xt}(s)\right)\right\|^2_{\ur_m}  +C(L)\left(1+\frac{1}{\lambda}\right) \left\|\mathcal{D}_{\tx}Y_{Xt}(s)\right\|^2_{\hr_m}. \label{thm4_10}
\end{align}
Therefore, we have
\begin{equation}\label{thm4_6}
	\begin{split}
		&\lambda \int_t^T \left\|\hat{Z}\left(s,\mathcal{D}_{\tx}Y_{Xt}(s),\mathcal{D}_{\tx}P_{Xt}(s),\mathcal{D}_{\tx}Q_{Xt}(s)\right) \right\|^2_{\ur_m} ds\\
		\le& \ C(L,T)\left(1+\frac{1}{\lambda}\right)\left(\sup_{t\le s\le T}\left\|\mathcal{D}_{\tx}Y_{Xt}(s)\right\|^2_{\hr_m}+\sup_{t\le s\le T}\left\| \mathcal{D}_{\tx}P_{Xt}(s)\right\|^2_{\hr_m}\right).
	\end{split}
\end{equation}
By using the standard arguments of SDEs, we have
\begin{equation}\label{thm4_4}
	\begin{split}
		&\sup_{t\le s\le T}\left\|\mathcal{D}_{\tx}Y_{Xt}(s)\right\|^2_{\hr_m}\\
		\le&\  C(L,T)\left(\left\|\tx\right\|^2_{\hr_m}+\int_t^T\left\|\hat{Z}\left(s,\mathcal{D}_{\tx}Y_{Xt}(s),\mathcal{D}_{\tx}P_{Xt}(s),\mathcal{D}_{\tx}Q_{Xt}(s)\right)\right\|^2_{\ur_m} ds\right).
	\end{split}
\end{equation}
Using the standard arguments of BSDEs and from Estimate \eqref{thm4_4}, we also have	
\begin{equation}\label{thm4_5}
	\begin{split}
		&\sup_{t\le s\le T}\left\| \mathcal{D}_{\tx}P_{Xt}(s)\right\|^2_{\hr_m}+\sum_{j=1}^n\left\| \mathcal{D}_{\tx}Q^j_{Xt}\right\|^2_{L^2(t,T;\hr_m)}\\
		\le&\  C(L,T)\left(\left\|\tx\right\|^2_{\hr_m}+\int_t^T\left\|\hat{Z}\left(s,\mathcal{D}_{\tx}Y_{Xt}(s),\mathcal{D}_{\tx}P_{Xt}(s),\mathcal{D}_{\tx}Q_{Xt}(s)\right)\right\|^2_{\ur_m} ds\right).
	\end{split}
\end{equation}
Substituting \eqref{thm4_4} and \eqref{thm4_5} into \eqref{thm4_6}, there is a constant $c(L,T)$ depending only on $(L,T)$, such that when $\lambda\geq c(L,T)$, we have
\begin{equation}\label{thm4_7}
	\begin{split}
		&\int_t^T \left\|\hat{Z}\left(s,\mathcal{D}_{\tx}Y_{Xt}(s),\mathcal{D}_{\tx}P_{Xt}(s),\mathcal{D}_{\tx}Q_{Xt}(s)\right)\right\|^2_{\ur_m} ds\le C(L,T,\lambda) \left\|\tx\right\|^2_{\hr_m}.
	\end{split}
\end{equation}
Substituting \eqref{thm4_7} back into \eqref{thm4_4} and \eqref{thm4_5}, we obtain \eqref{thm4_0}. Furthermore, if Assumption (A3)-(ii) is satisfied, from the equation of \eqref{thm4_10} and Condition \eqref{mono_12}, we further deduce
\begin{equation}\label{thm4_11}
	\begin{split}
		&2 \lambda\int_t^T \left\|\hat{Z}\left(s,\mathcal{D}_{\tx}Y_{Xt}(s),\mathcal{D}_{\tx}P_{Xt}(s),\mathcal{D}_{\tx}Q_{Xt}(s)\right)\right\|^2_{\ur_m}ds \le \left(\mathcal{D}_{\tx} P_{Xt}(t),\tx\right)_{\hr_m}
	\end{split}
\end{equation}
Substituting \eqref{thm4_5} into \eqref{thm4_11}, we deduce that
\begin{align*}
	&2 \lambda\int_t^T \left\|\hat{Z}\left(s,\mathcal{D}_{\tx}Y_{Xt}(s),\mathcal{D}_{\tx}P_{Xt}(s),\mathcal{D}_{\tx}Q_{Xt}(s)\right)\right\|^2_{\ur_m}ds \\
	\le&\  C(L,T)\left\|\tx\right\|_{\hr_m} \Bigg(\left\|\tx\right\|_{\hr_m}+\sqrt{\int_t^T\left\|\hat{Z}\left(s,\mathcal{D}_{\tx}Y_{Xt}(s),\mathcal{D}_{\tx}P_{Xt}(s),\mathcal{D}_{\tx}Q_{Xt}(s)\right)\right\|^2_{\ur_m} ds}\ \Bigg),
\end{align*}
from which we obtain \eqref{thm4_7}, and so \eqref{thm4_0} follows. 

\subsection{Proof of Theorem~\ref{thm8}}\label{pf:thm5}

We first establish \eqref{thm5_0}. For $\epsilon\in(0,1)$ and $s\in[t,T]$, we denote by
\begin{align*}
	&\de^\epsilon Y(s):=\frac{1}{\epsilon}\left(Y_{X^\epsilon t}(s)-Y_{Xt}(s)\right),\quad \delta^\epsilon Y(s):=\de Y^\epsilon(s)-\mathcal{D}_{\tx}Y_{Xt}(s).
\end{align*}
Similar definitions apply to $\de^\epsilon P(s),\de^\epsilon Q(s)$ and $\delta^\epsilon P(s),\delta^\epsilon Q(s)$. For the sake of convenience, denote by $\theta:=(Y,P,Q)\in\hr_m\times\hr_m\times\hr_m^n$ for $Y,P\in\hr_m$ and $Q\in\hr_m^n$, and use the notations $\theta_{Xt}(s):=(Y_{Xt}(s),P_{Xt}(s),Q_{Xt}(s))$, $\mathcal{D}_{\tx}\theta_{Xt}(s):=\left(\mathcal{D}_{\tx}Y_{Xt}(s),\mathcal{D}_{\tx}P_{Xt}(s),\mathcal{D}_{\tx}Q_{Xt}(s)\right)$, $\de^\epsilon\theta(s):=\left(\de^\epsilon Y(s),\de^\epsilon P(s),\de^\epsilon Q(s)\right)$ and also $\delta^\epsilon\theta(s):=\left(\delta^\epsilon Y(s),\delta^\epsilon P(s),\delta^\epsilon Q(s)\right)$. From \eqref{thm2_2}, we see that 
\begin{equation}\label{thm5_4}
	\begin{split}
		&\sup_{t\le s\le T}\|\de^\epsilon Y(s)\|_{\hr_m}+\sup_{t\le s\le T}\|\de^\epsilon P(s)\|_{\hr_m} +\sum_{j=1}^n \left\|\de^\epsilon Q^j\right\|_{L^2(t,T;\hr_m)}\le  C(L,T,\lambda)\left\|\tx\right\|_{\hr_m}.
	\end{split}
\end{equation}	
Define the map $\hat{Z}:(\hr_m\times\hr_m\times\hr_m^n)\times[0,T]\times(\hr_m\times\hr_m\times\hr_m^n)\to\ur_m$ as
\begin{align*}
	\hat{Z}\left(\theta,s;\tilde{\theta}\right):=D_X\hv(\theta,s)\left(\ty\right)+D_P\hv(\theta,s)\left(\tp\right)+\sum_{j=1}^nD_{Q^j}\hv(\theta,s)\left(\tq^j\right).
\end{align*}
Then, the process $\delta^\epsilon\theta$ satisfies the following FBSDEs: for $s\in[t,T]$,
\begin{align*}
	\delta^\epsilon Y(s)=&\int_t^s \bigg[\mathcal{F}_1(r)\delta^\epsilon Y(r)+\mathcal{F}_2(r)\int_0^1 \Big[\hat{Z}\left(\theta_{Xt}(r)+\gamma\epsilon\de^\epsilon \theta(r),r;\de^\epsilon \theta(r)\right)\\
	&\qquad\qquad\qquad\qquad\qquad\qquad\qquad -\hat{Z}(\theta_{Xt}(r),r;\mathcal{D}_{\tx}{\theta}_{Xt}(r))\Big] d\gamma \bigg]dr\\
	&+\sum_{j=1}^n \int_t^s \bigg[\mathcal{A}^j_1(r)\delta^\epsilon Y(r)+\mathcal{A}^j_2(r)\int_0^1 \Big[\hat{Z}\left(\theta_{Xt}(r)+\gamma\epsilon\de^\epsilon \theta(r),r;\de^\epsilon \theta(r)\right)\\
	&\ \qquad\qquad\qquad\qquad\qquad\qquad\qquad\qquad -\hat{Z}(\theta_{Xt}(r),r;\mathcal{D}_{\tx}{\theta}_{Xt}(r))\Big] d\gamma \bigg]dw_j(r);\\
	\delta^\epsilon P(s)=&\int_0^1 \Big[D^2_X G_T\left(Y_{Xt}(T)+\gamma\epsilon\de^\epsilon Y(T)\right)\left(\de^\epsilon Y(T)\right)-D^2_X G_T(Y_{Xt}(T))\left(\mathcal{D}_{\tx}Y_{Xt}(T)\right)\Big]d\gamma\\
	&+\int_s^T \Bigg[\mathcal{F}^*_1(r)\delta^\epsilon P(r)+\sum_{j=1}^n {\mathcal{A}^j_1}^*(r)\delta^\epsilon Q^j(r)\\
	&\ \quad\qquad +\int_0^1 \Big[D_X^2 G\left(Y_{Xt}(r)+\gamma\epsilon\de^\epsilon Y(r),\hv\left(\theta_{Xt}(r)+\gamma\epsilon\de^\epsilon\theta (r),r\right),r\right)\left(\de^\epsilon Y(r)\right) \\
	&\ \qquad\qquad\qquad -D_X^2 G\left(Y_{Xt}(r),\hv\left(\theta_{Xt}(r),r\right),r\right)\left(\mathcal{D}_{\tx}Y_{Xt}(r)\right)\Big]d\gamma\\
	&\ \quad\qquad +\int_0^1 \Big[D_V D_X G\left(Y_{Xt}(r)+\gamma\epsilon\de^\epsilon Y(r),\hv\left(\theta_{Xt}(r)+\gamma\epsilon\de^\epsilon\theta (r),r\right),r\right)\\
	&\quad\qquad\qquad\qquad\qquad\qquad \left(\hat{Z}\left(\theta_{Xt}(r)+\gamma\epsilon\de^\epsilon \theta(r),r;\de^\epsilon \theta(r)\right)\right) \\
	&\ \qquad\qquad\qquad -D_V D_X G\left(Y_{Xt}(r),\hv(\theta_{Xt}(r),r),r\right)\left(\hat{Z}(\theta_{Xt}(r),r;\tilde{\theta}_{\tx Xt}(r)) \right)\Big]d\gamma \Bigg]dr\\
	&-\sum_{j=1}^n \int_s^T \delta^\epsilon Q^j(r)dw_j(r).
\end{align*}
From the optimal conditions for $\theta_{Xt}(s)$ and $\theta_{X^\epsilon t}(s)$, for $s\in(t,T)$,
\begin{equation}\label{thm5_1}
	\begin{split}
		0=\ &\mathcal{F}^*_2(s)\de^\epsilon P(s)+\sum_{j=1}^n {\mathcal{A}^j_2}^*(s)\de^\epsilon Q^j(s)\\
		& +\int_0^1 D_XD_V G\left(Y_{Xt}(s)+\gamma\epsilon\de^\epsilon Y(s),\hv\left(\theta_{Xt}(s)+\gamma\epsilon\de^\epsilon\theta (s),s\right),s\right)\left(\de^\epsilon Y(s)\right)d\gamma\\
		& +\int_0^1 D_V^2 G\left(Y_{Xt}(s)+\gamma\epsilon\de^\epsilon Y(s),\hv\left(\theta_{Xt}(s)+\gamma\epsilon\de^\epsilon\theta (s),s\right),s\right)\\
		&\qquad\qquad\qquad \left(\hat{Z}(\theta_{Xt}(s)+\gamma\epsilon\de^\epsilon \theta(s),s;\de^\epsilon \theta(s))\right)d\gamma.
	\end{split}
\end{equation}
From Condition \eqref{optimal_condition} for $\mathcal{D}_{\tx}{\theta}_{Xt}(s)$, for $s\in(t,T)$,
\begin{equation}\label{thm5_2}
	\begin{split}
		0=\ &\mathcal{F}^*_2(s)\mathcal{D}_{\tx}P_{Xt}(s)+\sum_{j=1}^n {\mathcal{A}^j_2}^*(s)\mathcal{D}_{\tx}Q_{Xt}^j(s)\\
		&+D_XD_V G\left(Y_{Xt}(s),\hv(\theta_{Xt}(s),s),s\right)\left(\mathcal{D}_{\tx}Y_{Xt}(s)\right)\\
		&+D_V^2 G\left(Y_{Xt}(s),\hv(\theta_{Xt}(s),s),s\right)\left(\hat{Z}\left(\theta_{Xt}(s),s; \mathcal{D}_{\tx}{\theta}_{Xt}(s)\right)\right).
	\end{split}
\end{equation}
From \eqref{thm5_1} and \eqref{thm5_2}, we deduce that, for $s\in(t,T)$,
\begin{equation}\label{thm5_3}
	\begin{split}
		&\mathcal{F}^*_2(s)\delta^\epsilon P(s)+\sum_{j=1}^n {\mathcal{A}^j_2}^*(s)\delta^\epsilon Q^j(s)\\
		=& -\int_0^1 \Big[D_XD_V G\left(Y_{Xt}(s)+\gamma\epsilon\de^\epsilon Y(s),\hv\left(\theta_{Xt}(s)+\gamma\epsilon\de^\epsilon\theta (s),s\right),s\right)\left(\de^\epsilon Y(s)\right)\\
		&\quad\qquad -D_XD_V G\left(Y_{Xt}(s),\hv(\theta_{Xt}(s),s),s\right)\left(\mathcal{D}_{\tx}Y_{Xt}(s)\right)\Big]d\lambda\\
		&  -\int_0^1 \Big[D_V^2 G\left(Y_{Xt}(s)+\gamma\epsilon\de^\epsilon Y(s),\hv\left(\theta_{Xt}(s)+\gamma\epsilon\de^\epsilon\theta (s),s\right),s\right)\\
		&\qquad\qquad\qquad \left(\hat{Z}\left(\theta_{Xt}(s)+\gamma\epsilon\de^\epsilon \theta(s),s;\de^\epsilon \theta(s)\right)\right)\\
		&\quad\qquad -D_V^2 G\left(Y_{Xt}(s),\hv(\theta_{Xt}(s),s),s\right)\left(\hat{Z}\left(\theta_{Xt}(s),s; \mathcal{D}_{\tx}{\theta}_{Xt}(s)\right)\right)\Big]d\gamma.
	\end{split}
\end{equation}
From the FBSDEs for $\delta^\epsilon\theta$ and \eqref{thm5_3}, we can derive that for $s\in[t,T]$,
\begin{align*}
	&\frac{d}{ds}\left(\delta^\epsilon P(s),\delta^\epsilon Y(s) \right)_{\hr_m} \\
	=& -\bigg(\int_0^1 \Big[D_X^2 G\left(Y_{Xt}(s)+\gamma\epsilon\de^\epsilon Y(s),\hv\left(\theta_{Xt}(s)+\gamma\epsilon\de^\epsilon\theta (s),s\right),s\right)\left(\de^\epsilon Y(s)\right) \\
	&\qquad\qquad-D_X^2 G\left(Y_{Xt}(s),\hv\left(\theta_{Xt}(s),s\right),s\right)\left(\mathcal{D}_{\tx}Y_{Xt}(s)\right)\Big] d\gamma,\  \delta^\epsilon Y(s) \bigg)_{\hr_m}\\
	& -\bigg(\int_0^1 \Big[D_V D_X G\left(Y_{Xt}(s)+\gamma\epsilon\de^\epsilon Y(s),\hv\left(\theta_{Xt}(s)+\gamma\epsilon\de^\epsilon\theta (s),s\right),s\right)\\
	&\quad\qquad\qquad\qquad\qquad \left(\hat{Z}(\theta_{Xt}(s)+\gamma\epsilon\de^\epsilon \theta(s),s;\de^\epsilon \theta(s))\right) \\
	&\qquad\qquad-D_V D_X G\left(Y_{Xt}(s),\hv(\theta_{Xt}(s),s),s\right)\left(\hat{Z}\left(\theta_{Xt}(s),s;\mathcal{D}_{\tx}{\theta}_{Xt}(s)\right)\right)\Big]d\gamma,\  \delta^\epsilon Y(s) \bigg)_{\hr_m} \\
	& -\bigg( \int_0^1 \Big[D_XD_V G\left(Y_{Xt}(s)+\gamma\epsilon\de^\epsilon Y(s),\hv\left(\theta_{Xt}(s)+\gamma\epsilon\de^\epsilon\theta (s),s\right),s\right)\left(\de^\epsilon Y(s)\right)\\
	&\qquad\qquad -D_XD_V G\left(Y_{Xt}(s),\hv(\theta_{Xt}(s),s),s\right)\left(\mathcal{D}_{\tx}Y_{Xt}(s)\right)\Big]d\gamma,\\ 
	&\quad\qquad\qquad \int_0^1 \Big[\hat{Z}\left(\theta_{Xt}(s)+\gamma\epsilon\de^\epsilon \theta(s),s;\de^\epsilon \theta(s)\right)-\hat{Z}\left(\theta_{Xt}(s),s;\mathcal{D}_{\tx}{\theta}_{Xt}(s)\right)\Big]d\lambda\bigg)_{\ur_m}\\
	& -\Big(\int_0^1 \Big[D_V^2 G\left(Y_{Xt}(s)+\gamma\epsilon\de^\epsilon Y(s),\hv\left(\theta_{Xt}(s)+\gamma\epsilon\de^\epsilon\theta (s),s\right),s\right)\\
	&\qquad\qquad\qquad \left(\hat{Z}\left(\theta_{Xt}(s)+\gamma\epsilon\de^\epsilon \theta(s),s;\de^\epsilon \theta(s)\right)\right)\\
	&\qquad\qquad -D_V^2 G\left(Y_{Xt}(s),\hv(\theta_{Xt}(s),s),s\right)\left(\hat{Z}(\theta_{Xt}(s),s; \mathcal{D}_{\tx}{\theta}_{Xt}(s))\right)\Big]d\gamma,\\
	&\quad\qquad \int_0^1\left[\hat{Z}\left(\theta_{Xt}(s)+\gamma\epsilon\de^\epsilon \theta(s),s;\de^\epsilon \theta(s)\right)-\hat{Z}\left(\theta_{Xt}(s),s;\mathcal{D}_{\tx}{\theta}_{Xt}(s)\right)\right]d\gamma\Big)_{\ur_m}.
\end{align*}
Therefore, we have
\small
\begin{align}
	&\left( D^2_X G_T(Y_{Xt}(T))\left(\delta^\epsilon Y(T)\right),\  \delta^\epsilon Y(T)\right)_{\hr_m} \notag \\
	&+\int_t^T \bigg[\left(D_X^2 G\left(Y_{Xt}(s),\hv(\theta_{Xt}(s),s),s\right)(\delta^\epsilon Y(s)),\  \delta^\epsilon Y(s) \right)_{\hr_m}  \notag\\
	&\quad\qquad +2\bigg(D_X D_V G\left(Y_{Xt}(s),\hv(\theta_{Xt}(s),s),s\right)(\delta^\epsilon Y(s)),  \notag\\
	&\quad\qquad\qquad \int_0^1 \left[\hat{Z}\left(\theta_{Xt}(s)+\gamma\epsilon\de^\epsilon \theta(s),s;\de^\epsilon \theta(s)\right)-\hat{Z}\left(\theta_{Xt}(s),s;\mathcal{D}_{\tx}{\theta}_{Xt}(s)\right)\right]d\gamma \bigg)_{\ur_m}  \notag\\
	&\quad\qquad +\bigg(D_V^2 G\left(Y_{Xt}(s),\hv(\theta_{Xt}(s),s),s\right)  \notag\\
	&\quad\qquad\qquad\qquad \left(\int_0^1\left[\hat{Z}(\theta_{Xt}(s)+\gamma\epsilon\de^\epsilon \theta(s),s;\de^\epsilon \theta(s))-\hat{Z}\left(\theta_{Xt}(s),s; \mathcal{D}_{\tx}{\theta}_{Xt}(s)\right)\right]d\gamma \right), \notag\\
	&\quad\qquad\qquad \int_0^1\left[\hat{Z}(\theta_{Xt}(s)+\gamma\epsilon\de^\epsilon \theta(s),s;\de^\epsilon \theta(s))-\hat{Z}\left(\theta_{Xt}(s),s;\mathcal{D}_{\tx}{\theta}_{Xt}(s)\right)\right]d\lambda \bigg)_{\ur_m} \bigg] ds \notag\\
	=&-\bigg(\int_0^1 \Big[D^2_X G_T\left(Y_{Xt}(T)+\gamma\epsilon\de^\epsilon Y(T)\right)(\de^\epsilon Y(T)) \notag\\
	&\qquad\qquad -D^2_X G_T(Y_{Xt}(T))(\de^\epsilon Y(T))\Big]d\gamma,\quad \delta^\epsilon Y(T)\bigg)_{\hr_m} \notag\\
	& -\int_t^T \Bigg[\bigg(\int_0^1 \Big[D_X^2 G\left(Y_{Xt}(s)+\gamma\epsilon\de^\epsilon Y(s),\hv(\theta_{Xt}(s)+\gamma\epsilon\de^\epsilon\theta (s),s),s\right)(\de^\epsilon Y(s))  \notag\\
	&\qquad\qquad\qquad -D_X^2 G\left(Y_{Xt}(s),\hv(\theta_{Xt}(s),s),s\right)(\de^\epsilon Y(s))\Big]d\gamma,\quad \delta^\epsilon Y(s) \bigg)_{\hr_m} \notag\\
	&\quad\qquad +\bigg(\int_0^1 \Big[D_V D_X G\left(Y_{Xt}(s)+\gamma\epsilon\de^\epsilon Y(s),\hv(\theta_{Xt}(s)+\gamma\epsilon\de^\epsilon\theta (s),s),s\right) \notag\\
	&\quad\qquad\qquad\qquad\qquad\qquad \left(\hat{Z}(\theta_{Xt}(s)+\gamma\epsilon\de^\epsilon \theta(s),s;\de^\epsilon \theta(s))\right)  \notag\\
	&\quad\qquad\qquad\qquad -D_V D_X G\left(Y_{Xt}(s),\hv(\theta_{Xt}(s),s),s\right) \notag\\
	&\quad\qquad\qquad\qquad\qquad\qquad \left(\hat{Z}(\theta_{Xt}(s)+\gamma\epsilon\de^\epsilon \theta(s),s;\de^\epsilon \theta(s))\right)\Big]d\gamma,\quad \delta^\epsilon Y(s) \bigg)_{\hr_m}  \notag\\
	&\quad\qquad +\bigg( \int_0^1 \Big[D_XD_V G\left(Y_{Xt}(s)+\gamma\epsilon\de^\epsilon Y(s),\hv(\theta_{Xt}(s)+\gamma\epsilon\de^\epsilon\theta (s),s),s\right)(\de^\epsilon Y(s)) \notag\\
	&\quad\qquad\qquad\qquad -D_XD_V G\left(Y_{Xt}(s),\hv(\theta_{Xt}(s),s),s\right)(\de^\epsilon Y(s))\Big]d\gamma, \notag\\ 
	&\quad\qquad\qquad \int_0^1\left[\hat{Z}(\theta_{Xt}(s)+\gamma\epsilon\de^\epsilon \theta(s),s;\de^\epsilon \theta(s))-\hat{Z}(\theta_{Xt}(s),s;\mathcal{D}_{\tx}{\theta}_{Xt}(s))\right] d\gamma\bigg)_{\ur_m} \notag\\
	&\quad\qquad +\bigg(\int_0^1 \Big[D_V^2 G\left(Y_{Xt}(s)+\gamma\epsilon\de^\epsilon Y(s),\hv(\theta_{Xt}(s)+\gamma\epsilon\de^\epsilon\theta (s),s),s\right) \notag\\
	&\qquad\qquad\qquad\qquad\qquad \left(\hat{Z}(\theta_{Xt}(s)+\gamma\epsilon\de^\epsilon \theta(s),s;\de^\epsilon \theta(s))\right) \notag\\
	&\quad\qquad\qquad\qquad -D_V^2 G\left(Y_{Xt}(s),\hv(\theta_{Xt}(s),s),s\right)\left(\hat{Z}(\theta_{Xt}(s)+\gamma\epsilon\de^\epsilon \theta(s),s;\de^\epsilon \theta(s))\right)\Big]d\gamma, \notag\\
	&\quad\qquad\qquad \int_0^1\left[\hat{Z}(\theta_{Xt}(s)+\gamma\epsilon\de^\epsilon \theta(s),s;\de^\epsilon \theta(s))-\hat{Z}(\theta_{Xt}(s),s;\mathcal{D}_{\tx}{\theta}_{Xt}(s))\right] d\gamma\Big)_{\ur_m} \Bigg] ds \notag\\
	&=: \de I(\epsilon).  \label{thm5_12}
\end{align}
\normalsize
From Condition \eqref{mono_11}, we deduce that
\begin{align*}
	&\left( D^2_X G_T(Y_{Xt}(T))(\delta^\epsilon Y(T)),\  \delta^\epsilon Y(T)\right)_{\hr_m} \notag\\
	&+\int_t^T \Bigg[\left(D_X^2 G\left(Y_{Xt}(s),\hv(\theta_{Xt}(s),s),s\right)(\delta^\epsilon Y(s)),\  \delta^\epsilon Y(s) \right)_{\hr_m} \notag\\
	&\quad\qquad +2\bigg(D_X D_V G\left(Y_{Xt}(s),\hv(\theta_{Xt}(s),s),s\right)(\delta^\epsilon Y(s)), \notag\\
	&\quad\qquad\qquad \int_0^1 \left[\hat{Z}\left(\theta_{Xt}(s)+\gamma\epsilon\de^\epsilon \theta(s),s;\de^\epsilon \theta(s)\right)-\hat{Z}\left(\theta_{Xt}(s),s;\mathcal{D}_{\tx}{\theta}_{Xt}(s)\right)\right]d\gamma \bigg)_{\ur_m}  \notag\\
	&\quad\qquad +\bigg(D_V^2 G\left(Y_{Xt}(s),\hv(\theta_{Xt}(s),s),s\right) \notag\\
	&\ \quad\qquad\qquad\qquad \left(\int_0^1\left[\hat{Z}(\theta_{Xt}(s)+\gamma\epsilon\de^\epsilon \theta(s),s;\de^\epsilon \theta(s))-\hat{Z}\left(\theta_{Xt}(s),s; \mathcal{D}_{\tx}{\theta}_{Xt}(s)\right)\right]d\gamma \right), \notag\\
	&\quad\qquad\qquad \int_0^1\left[\hat{Z}(\theta_{Xt}(s)+\gamma\epsilon\de^\epsilon \theta(s),s;\de^\epsilon \theta(s))-\hat{Z}\left(\theta_{Xt}(s),s;\mathcal{D}_{\tx}{\theta}_{Xt}(s)\right)\right]d\gamma \bigg)_{\ur_m} \Bigg] ds  \notag\\
	\geq&\  \lambda \int_t^T \left\|\int_0^1[ \hat{Z}(\theta_{Xt}(s)+\gamma\epsilon\de^\epsilon \theta(s),s;\de^\epsilon \theta(s))-\hat{Z}\left(\theta_{Xt}(s),s;\mathcal{D}_{\tx}{\theta}_{Xt}(s)\right)] d\gamma \right\|_{\ur_m}^2ds  \notag\\
	&\ -C(L,T)\left(1+\frac{1}{\lambda}\right)\sup_{t\le s\le T}\| \delta^\epsilon Y(s)\|_{\hr_m}^2.   
\end{align*}
Therefore, from \eqref{thm5_12} and the last inequality, we obtain
\begin{align}
	&\lambda \int_t^T \left\|\int_0^1[ \hat{Z}(\theta_{Xt}(s)+\gamma\epsilon\de^\epsilon \theta(s),s;\de^\epsilon \theta(s))-\hat{Z}\left(\theta_{Xt}(s),s;\mathcal{D}_{\tx}{\theta}_{Xt}(s)\right)] d\gamma \right\|_{\ur_m}^2ds \notag\\
	\le & \ C(L,T)\left(1+\frac{1}{\lambda}\right)\sup_{t\le s\le T}\| \delta^\epsilon Y(s)\|_{\hr_m}^2+|\de I(\epsilon)|. \label{thm5_5}
\end{align}
From the SDE for $\delta^\epsilon Y(\cdot)$ and Assumption (A1), we have
\begin{equation}\label{thm5_6}
	\begin{split}
		\sup_{t\le s\le T}\| \delta^\epsilon Y(s)\|_{\hr_m}^2\le & C(L,T) \int_t^T \bigg\|\int_0^1\Big[ \hat{Z}(\theta_{Xt}(s)+\gamma\epsilon\de^\epsilon \theta(s),s;\de^\epsilon \theta(s))\\
		&\qquad\qquad\qquad\qquad -\hat{Z}\left(\theta_{Xt}(s),s;\mathcal{D}_{\tx}{\theta}_{Xt}(s)\right)\Big] d\gamma \bigg\|_{\ur_m}^2ds.
	\end{split}	
\end{equation}
Substituting \eqref{thm5_6} into \eqref{thm5_5}, there exists a constant $c(L,T)$ depending only on $(L,T)$, such that when $\lambda\geq c(L,T)$, we have
\begin{equation*}
	\begin{split}
		& \int_t^T \left\|\int_0^1[ \hat{Z}(\theta_{Xt}(s)+\gamma\epsilon\de^\epsilon \theta(s),s;\de^\epsilon \theta(s))-\hat{Z}\left(\theta_{Xt}(s),s;\mathcal{D}_{\tx}{\theta}_{Xt}(s)\right)] d\gamma \right\|_{\ur_m}^2ds \le  C(L,T,\lambda) \cdot |\de I(\epsilon)|.
	\end{split}
\end{equation*}
From Estimates \eqref{thm4_0} and \eqref{thm5_4}, Assumption (A2), and the dominated convergence theorem, for instance, we have
\begin{align*}
	\lim_{\epsilon\to0}&\bigg(\int_0^1 \Big[D^2_X G_T\left(Y_{Xt}(T)+\gamma\epsilon\de^\epsilon Y(T)\right)(\de^\epsilon Y(T)) \notag\\
	&\qquad\qquad -D^2_X G_T(Y_{Xt}(T))(\de^\epsilon Y(T))\Big]d\gamma,\quad \delta^\epsilon Y(T)\bigg)_{\hr_m} =0.
\end{align*}
Using similar approach for all other terms, we can deduce that $\lim_{\epsilon\to0}\de I(\epsilon)=0$, which further implies that
\begin{equation*}
	\begin{split}
		\lim_{\epsilon\to0} \int_t^T \bigg\|\int_0^1\Big[& \hat{Z}(\theta_{Xt}(s)+\gamma\epsilon\de^\epsilon \theta(s),s;\de^\epsilon \theta(s))-\hat{Z}\left(\theta_{Xt}(s),s;\mathcal{D}_{\tx}{\theta}_{Xt}(s)\right)\Big] d\gamma \bigg\|_{\ur_m}^2ds=0,
	\end{split}
\end{equation*}
and then, from \eqref{thm5_6}, $\lim_{\epsilon\to0} \sup_{t\le s\le T}\| \delta^\epsilon Y(s)\|_{\hr_m}=0$. Then, from the BSDE for $(\delta^\epsilon P,\delta^\epsilon Q)$, we obtain \eqref{thm5_0}. For the case when Assumption (A3)-(ii) is satisfied, we can also deduce that
\begin{equation*}
	\begin{split}
		&\lambda  \int_t^T \left\|\int_0^1[ \hat{Z}(\theta_{Xt}(s)+\gamma\epsilon\de^\epsilon \theta(s),s;\de^\epsilon \theta(s))-\hat{Z}\left(\theta_{Xt}(s),s;\mathcal{D}_{\tx}{\theta}_{Xt}(s)\right)] d\gamma \right\|_{\ur_m}^2ds\le  |\de I(\epsilon)|,
	\end{split}
\end{equation*}
from which we conclude the claim on the convergence to zero. The differentiability is then a consequence of \eqref{thm5_0} and Lemma~\ref{thm4}. The proof of continuity of $({D}_{\tx}Y_{Xt}(s),{D}_{\tx}P_{Xt}(s),{D}_{\tx}Q_{Xt}(s))$ in $X$ is similar to that of \eqref{thm5_0}, and is omitted here.

\section{Proof of Statements in Section~\ref{sec:value_Hilbert}}\label{app:02}

\subsection{Proof of Lemma~\ref{lem:8}}\label{pf:lem:8}

In view of Section~\ref{sec:FB_Hilbert}, we only need to give the proof for the case when (A3)-(i) is satisfied. We first prove \eqref{lem8_3} as follows. Based on \eqref{v_1}, \eqref{hv1} and \eqref{thm2_1}, there exists a constant $c(L,T)$ depending only on $(L,T)$, such that when $\lambda\geq c(L,T)$, we have
\begin{align*}
	|\V(X\otimes m,t)|&\le C(L,T)\left(1+\sup_{t\le s\le T}\|Y_{Xt}(s)\|^2_{\hr_m}+\left\|\hv_{Xt}\right\|^2_{L^2(t,T;\ur_m)}\right)\\
	&\le C(L,T,\lambda)\left(1+\|X\|^2_{\hr_m}\right).
\end{align*}
We now prove \eqref{lem8_1}. From the definition of $\V$ in \eqref{def:V}, for any $X,X'\in\hr_m$, both are independent of $\mathcal{W}_t$, we have
\begin{equation}\label{lem8_9}
	J_{X'mt}\left(\hv_{X't}\right)-J_{Xmt}\left(\hv_{X't}\right)\le \V\left(X',t\right)-\V(X,t)\le J_{X'mt}\left(\hv_{Xt}\right)-J_{Xmt}\left(\hv_{Xt}\right).
\end{equation}	
For the upper bound, from Assumption (A2), we have
\begin{equation}\label{lem8_10}
	\begin{split}
		&J_{X'mt}\left(\hv_{Xt}\right)-J_{Xmt}\left(\hv_{Xt}\right)\\
		=&\int_t^T\left[G\left(X_{X't}^{\hv_{Xt}}(s),\hv_{Xt}(s),s\right)-G\left(Y_{Xt}(s),\hv_{Xt}(s),s\right)\right] ds\\
		&\qquad +\left[G_T\left(X_{X't}^{\hv_{Xt}}(T)\right)-G_T\left(Y_{Xt}(T)\right)\right] \\
		\le& \int_t^T \left(D_X G\left(Y_{Xt}(s),\hv_{Xt}(s),s\right),X_{X't}^{\hv_{Xt}}(s)-Y_{Xt}(s)\right)_{\hr_m} ds\\
		& +\left(D_X G_T\left(Y_{Xt}(T)\right),X_{X't}^{\hv_{Xt}}(T)-Y_{Xt}(T)\right)_{\hr_m} \\
		&+C(L,T)\sup_{t\le s\le T}\left\|X_{X't}^{\hv_{Xt}}(s)-Y_{Xt}(s)\right\|^2_{\hr_m}.
	\end{split}
\end{equation}
From Assumption (A1) and the Gr\"onwall's inequality, we can obtain the following estimate
\begin{equation}\label{lem8_11}
	\begin{split}
		\sup_{t\le s\le T}\left\|X_{X't}^{\hv_{Xt}}(s)-Y_{Xt}(s)\right\|_{\hr_m}\le C(L,T)\left\|X'-X\right\|_{\hr_m}.
	\end{split}
\end{equation}
From the BSDE for $(P_{Xt},Q_{Xt})$, for $s\in(t,T)$,
\begin{align*}
	&\frac{d}{ds}\left(P_{Xt}(s),\ X_{X't}^{\hv_{Xt}}(s)-Y_{Xt}(s)\right)_{\hr_m}=-\left(D_X G\left(Y_{Xt}(s),\hv_{Xt}(s),s\right),\ X_{X't}^{\hv_{Xt}}(s)-Y_{Xt}(s)\right)_{\hr_m},
\end{align*}
and then,
\begin{equation}\label{lem8_12}
	\begin{split}
		\left(P_{Xt}(t),\ X'-X\right)_{\hr_m}=&\int_t^T \left(D_X G\left(Y_{Xt}(s),\hv_{Xt}(s),s\right),\ X_{X't}^{\hv_{Xt}}(s)-Y_{Xt}(s)\right)_{\hr_m} ds\\
		& +\left(D_X G_T (Y_{Xt}(T)),X_{X't}^{\hv_{Xt}}(T)-Y_{Xt}(T)\right)_{\hr_m}.
	\end{split}
\end{equation}
Substituting \eqref{lem8_11} and \eqref{lem8_12} back to \eqref{lem8_10}, we have
\begin{equation}\label{lem8_13}
	\begin{split}
		&J_{X'mt}\left(\hv_{Xt}\right)-J_{Xmt}\left(\hv_{Xt}\right)\le \left(P_{Xt}(t),X'-X\right)_{\hr_m} +C(L,T)\left\|X'-X\right\|_{\hr_m}^2.
	\end{split}
\end{equation}
Similarly, we can also obtain the lower bound
\begin{equation}\label{lem8_14}
	\begin{split}
		&J_{X'mt}\left(\hv_{X't}\right)-J_{Xmt}\left(\hv_{X't}\right)\geq \left(P_{X't}(t),X'-X\right)_{\hr_m} -C(L,T)\left\|X'-X\right\|_{\hr_m}^2.
	\end{split}
\end{equation}
From \eqref{thm2_2}, we have
\begin{align*}
	&\left\|P_{X't}(t)-P_{Xt}(t)\right\|_{\hr_m}\le C(L,T,\lambda)\left\|X'-X\right\|_{\hr_m}.
\end{align*}
Substituting the last inequality back to \eqref{lem8_14}, we have
\begin{equation}\label{lem8_15}
	\begin{split}
		&J_{X'mt}\left(\hv_{X't}\right)-J_{Xmt}\left(\hv_{X't}\right)\geq \left(P_{Xt}(t),X'-X\right)_{\hr_m} -C(L,T,\lambda)\left\|X'-X\right\|_{\hr_m}^2.
	\end{split}
\end{equation}
Combining \eqref{lem8_9}, \eqref{lem8_13} and \eqref{lem8_15}, we deduce:
\begin{align}\label{lem8_16}
	\left|\V\left(X'\otimes m,t\right)-\V(X\otimes m,t)-\left(P_{Xt}(t),X'-X\right)_{\hr_m}\right|\le C(L,T,\lambda) \left\|X'-X\right\|_{\hr_m}^2,
\end{align}
from which we obtain \eqref{lem8_1}. Then, \eqref{lem8_1.1} is a consequence of Theorem~\ref{thm8}. Estimates \eqref{lem8_4}-\eqref{lem8_5} are consequences of Theorems~\ref{thm2} and \ref{thm8}.	

\subsection{Proof of Lemma~\ref{lem9}}\label{pf:lem9}

As in the proof of Lemma~\ref{lem:8}, we only prove the case when (A3)-(i) is satisfied. We first prove \eqref{lem9_1}. From the flow property of the FBSDEs \eqref{FB:1}, we have for $0\le t\le t'\le s\le T$, 
\begin{align*}
	Y_{Xt}(s)=Y_{Y_{Xt}(t'),t'}(s),\quad P_{Xt}(s)=P_{Y_{Xt}(t'),t'}(s),\quad Q_{Xt}(s)=Q_{Y_{Xt}(t'),t'}(s).
\end{align*}
Therefore, from \eqref{thm2_2}, we have
\begin{equation}\label{lem9_11'}
	\begin{split}
		&\sup_{t\le s\le T}\|Y_{Xt'}(s)-Y_{Xt}(s)\|_{\hr_m}+\sup_{t\le s\le T}\|P_{Xt'}(s)-P_{Xt}(s)\|_{\hr_m} \\
		&+\sum_{j=1}^n \left\|Q^j_{Xt'}-Q^j_{Xt}\right\|_{L^2(t,T;\hr_m)}\\
		\le& \ C(L,T,\lambda)\left\|Y_{Xt}(t')-X\right\|_{\hr_m}.
	\end{split}
\end{equation}	
Note that when $\mathcal{A}^j_2=0$ (for $1\le j\le n$), the map $\hv$ does not depend on $Q$. By Cauchy's inequality and \eqref{hv1}, we have
\begin{equation}\label{lem9_11}
	\begin{split}
		\left\|Y_{Xt}(t')-X\right\|_{\hr_m}^2&=\left\| \int_t^{t'}F\left(Y_{Xt}(s),\hv_{Xt}(s),s\right) ds+\sum_{j=1}^n \int_t^{t'}A^j(Y_{Xt}(s),s) dw_j(s) \right\|^2_{\hr_m}\\
		&\le C(L,T) \int_t^{t'}\left[1+\|Y_{Xt}(s)\|_{\hr_m}^2+\left\|\hv(Y_{Xt}(s),s;P_{Xt}(s))\right\|_{\ur_m}^2\right]ds\\
		&\le C(L,T,\lambda) \int_t^{t'}\left(1+\|Y_{Xt}(s)\|_{\hr_m}^2+\|P_{Xt}(s)\|_{\hr_m}^2\right)ds\\
		&\le C(L,T,\lambda)\left(1+\|X\|_{\hr_m}^2\right)\left|t'-t\right|.
	\end{split}		
\end{equation}
From \eqref{lem9_11'} and \eqref{lem9_11}, we obtain \eqref{lem9_1}. From the flow property, we also have for $0\le t\le t'\le s\le T$, (i) ${D}_{\tx}Y_{Xt}(s)={D}_{{D}_{\tx}Y_{X t}(t')}Y_{Y_{Xt}(t'),t'}(s)$; (ii) ${D}_{\tx}P_{Xt}(s)={D}_{{D}_{\tx}Y_{X t}(t')}P_{Y_{Xt}(t'),t'}(s)$; (iii) ${D}_{\tx}Q_{Xt}(s)={D}_{{D}_{\tx}Y_{X t}(t')}Q_{Y_{Xt}(t'),t'}(s)$, and thus
\begin{align*}
	&\left\|{D}_{\tx}Y_{Xt'}(s)-{D}_{\tx}Y_{Xt}(s)\right\|_{\hr_m}=\left\| D_{\tx}Y_{Xt'}(s)- D_{D_{\tx}Y_{X t}(t')}Y_{Y_{Xt}(t'),t'}(s) \right\|_{\hr_m}, \\
	&\left\|{D}_{\tx}P_{Xt'}(s)-{D}_{\tx}P_{Xt}(s)\right\|_{\hr_m}=\left\| D_{\tx}P_{Xt'}(s)- D_{D_{\tx}Y_{X t}(t')}P_{Y_{Xt}(t'),t'}(s) \right\|_{\hr_m}, \\
	&\left\|{D}_{\tx}Q^j_{Xt'}(s)-{D}_{\tx}Q^j_{Xt}(s)\right\|_{\hr_m}=\left\| D_{\tx}Q^j_{Xt'}(s)- D_{D_{\tx}Y_{X t}(t')}Q^j_{Y_{Xt}(t'),t'}(s) \right\|_{\hr_m}.
\end{align*}	
In a similar way as \eqref{lem9_11}, we can deduce that 
\begin{equation}\label{lem9_12}
	\begin{split}
		\left\|D_{\tx}Y_{Xt}(t')-\tx\right\|_{\hr_m}^2\le C(L,T,\lambda)\|X\|_{\hr_m} ^2\left|t'-t\right|.
	\end{split}		
\end{equation}
From Estimates \eqref{lem9_11} and \eqref{lem9_12} and the fact that $\left(D_{\tx}Y_{Xt}(s),D_{\tx}P_{Xt}(s),D_{\tx}Q_{Xt}(s)\right)$ is linear in $\tx$ and continuous in $X$,  then, $\left(D_{\tx}Y_{Xt}(s),D_{\tx}P_{Xt}(s),D_{\tx}Q_{Xt}(s)\right)$ is continuous in $t$. We now prove \eqref{lem9_2}. Note that
\begin{equation}\label{lem9_13}
	\begin{split}
		&P_{Xt'}\left(t'\right)-P_{Xt}\left(t\right)\\
		=&\ P_{Xt'}\left(t'\right)-\e\left[P_{Xt}\left(t'\right)|X\right]-\int_t^{t'}\e\left[D_X\lr\left(Y_{Xt}(s),\hv_{Xt}(s),s;P_{Xt}(s),Q_{Xt}(s)\right)\Big|X\right]ds.
	\end{split}
\end{equation}  
From \eqref{thm2_1} and \eqref{thm2_2}, we deduce that
\begin{equation}\label{lem9_14}
	\begin{split}
		&\left\|\int_t^{t'} D_X\lr\left(Y_{Xt}(s),\hv_{Xt}(s),s;P_{Xt}(s),Q_{Xt}(s)\right)ds\right\|_{\hr_m}^2\\
		&\le C(L)\left|t'-t\right|\int_t^{t'}\bigg[1+\|Y_{Xt}(s)\|_{\hr_m}^2+\left\|\hv_{Xt}(s)\right\|_{\hr_m}^2+\|P_{Xt}(s)\|_{\hr_m}^2+\sum_{j=1}^n\|Q^j_{Xt}(s)\|_{\hr_m}^2\bigg]ds\\
		&\le C(L,T,\lambda)\left(1+\|X\|^2_{\hr_m}\right)\left|t'-t\right|.
	\end{split}	
\end{equation}
Substituting \eqref{lem9_1} and \eqref{lem9_14} into \eqref{lem9_13}, we have
\begin{align}\label{lem9_15}
	\left\|P_{Xt'}\left(t'\right)-P_{Xt}(t)\right\|_{\hr_m} \le C(L,T,\lambda)(1+\|X\|_{\hr_m})\left|t'-t\right|^{\frac{1}{2}},
\end{align}  
from which we obtain \eqref{lem9_2}. Similarly, since $D^2_X\V(X\otimes m,t)$ is continuous in $t$, we can use this to prove \eqref{lem9_3}. By dynamic programming principle, for any $\epsilon\in[0,T-t]$,
\begin{equation*}
	\V(X\otimes m,t)=\int_t^{t+\epsilon}G\left(Y_{Xt}(s),\hv_{Xt}(s),s\right)ds+\V(Y_{tX}(t+\epsilon)\otimes m,t+\epsilon),
\end{equation*}
so we have
\begin{equation}\label{lem8_17}
	\begin{split}
		&\frac{1}{\epsilon}\left[\V(X\otimes m,t+\epsilon)-\V(X\otimes m,t)\right]\\
		=&\frac{1}{\epsilon}[\V(X\otimes m,t+\epsilon)-\V(Y_{Xt}(t+\epsilon)\otimes m,t+\epsilon)]-\frac{1}{\epsilon}\int_t^{t+\epsilon}G\left(Y_{Xt}(s),\hv_{Xt}(s),s\right)ds.
	\end{split}	
\end{equation}
From Lemma~\ref{lem:8}, Estimate \eqref{lem9_11}, \cite[Theorem 2.2]{AB5}, and the continuity of $D_X\V$, $D^2_X\V$ and $Y_{Xt}$, we deduce that
\begin{equation}\label{lem9_18}
	\begin{split}
		&\lim_{\epsilon\to0} \frac{1}{\epsilon}[\V(X\otimes m,t+\epsilon)-\V(Y_{Xt}(t+\epsilon)\otimes m,t+\epsilon)]\\
		=&-\left(D_X\V(X\otimes m,t),\ F\left(X, \hv_{Xt}(t),t\right)\right)_{\hr_m}\\
		&-\frac{1}{2}\left(D_X^2\V(X\otimes m,t)\left(\sum_{j=1}^n A^j(X,t)\mathcal{N}_t^j\right),\   \sum_{j=1}^n A^j(X,t)\mathcal{N}_t^j\right)_{\hr_m},
	\end{split}
\end{equation} 
where $\mathcal{N}_t^j$'s are $n$ independent standard normal random variables which are independent of $X$ and the Wiener Process. From Assumption (A4) and the continuity of $\left(Y_{Xt},\hv_{Xt}\right)$, we have
\begin{align}\label{lem9_19}
	&\lim_{\epsilon\to0}\frac{1}{\epsilon}\int_t^{t+\epsilon}G\left(Y_{Xt}(s),\hv_{Xt}(s),s\right)ds=G\left(X,\hv_{Xt}(t),t\right).
\end{align}
Substituting \eqref{lem9_18} and \eqref{lem9_19} back to \eqref{lem8_17}, we deduce that 
\begin{align*}
	&\lim_{\epsilon\to0} \frac{1}{\epsilon}[\V(X\otimes m,t+\epsilon)-\V(X\otimes m,t)]\\
	=&-G\left(X,\hv_{Xt}(t),t\right)-\left(D_X\V(X\otimes m,t),\ G\left(X, \hv_{Xt}(t),t\right)\right)_{\hr_m}\\
	&-\frac{1}{2}\left(D_X^2\V(X\otimes m,t)\left(\sum_{j=1}^n A^j(X,t)\mathcal{N}_t^j\right),\   \sum_{j=1}^n A^j(X,t)\mathcal{N}_t^j\right)_{\hr_m},
\end{align*}
from which we obtain \eqref{lem9_3}. Then, \eqref{lem9_4} is a consequence of \eqref{lem9_3} and Estimates \eqref{lem8_4} and \eqref{lem8_4'}.

\subsection{Proof of Theorem~\ref{thm:9}}\label{pf:thm:9}

From \eqref{lem9_3}, the value function $\V$ is seen as a solution of Equation \eqref{Bellman}. We now show the uniqueness. Let $\V'$ be another solution of Equation \eqref{Bellman} possessing the properties \eqref{lem8_3}-\eqref{lem8_5}, \eqref{lem9_2} and \eqref{lem9_4}. For an arbitrary admissible control $V$ for Problem \ref{intr_5}, denote by $X_{Xt}^{V}$ the corresponding controlled state process. Then, from \cite[Theorem 2.2]{AB5}, we have
\begin{equation}\label{thm9_1}
	\begin{split}
		&\frac{d}{ds}\V'\left(X_{Xt}^{V}(s)\otimes m,s\right)\\
		=&\frac{\dd\V'}{\dd s}\left(X_{Xt}^{V}(s)\otimes m,s\right)+\left(D_X \V'\left(X_{Xt}^{V}(s),s\right),F\left(X_{Xt}^{V}(s),V(s),s\right)\right)_{\hr_m}\\
		& +\frac{1}{2}\Bigg(D_X^2\V'\left(X_{Xt}^{V}(s)\otimes m,s\right)\Bigg(\sum_{j=1}^n A^j\left(X_{Xt}^{V}(s),s\right)\mathcal{N}^j_s\Bigg),\sum_{j=1}^n A^j\left(X_{Xt}^{V}(s),s\right)\mathcal{N}^j_s \Bigg)_{\hr_m}\\
		=&\Big(D_X \V'\left(X_{Xt}^{V}(s)\otimes m,s\right),\ \mathcal{F}_2(s)\left(V(s)-\hv\left(X_{Xt}^{V}(s),s;D_X \V'\left(X_{Xt}^{V}(s)\otimes m,s\right)\right)\right)\Big)_{\hr_m}\\
		&-G\left(X_{Xt}^{V}(s),\hv\left(X_{Xt}^{V}(s),s;D_X \V'\left(X_{Xt}^{V}(s)\otimes m,s\right)\right),s\right).
	\end{split}
\end{equation}
From the first order condition for the optimal control $\hv$, we deduce that for $s\in(t,T)$,
\begin{align*}
	&\Big(D_X \V'\left(X_{Xt}^{V}(s)\otimes m,s\right),\  \mathcal{F}_2(s)\left(V(s)-\hv\left(X_{Xt}^{V}(s),s;D_X \V'\left(X_{Xt}^{V}(s)\otimes m,s\right)\right)\right)\Big)_{\hr_m}\\
	=& -\Big(D_V G\left(X_{Xt}^{V}(s),\hv\left(X_{Xt}^{V}(s),D_X \V'\left(X_{Xt}^{V}(s)\otimes m,s\right),s\right),s\right),\\
	&\qquad V(s)-\hv\left(X_{Xt}^{V}(s),s;D_X \V'(X_{Xt}^{V}(s),s)\right)\Big)_{\ur_m},
\end{align*}
which can be substituted back to \eqref{thm9_1} and integrating $s$ over $[t,T]$, from Assumption (A3)-(i), we deduce that
\begin{align*}
	&J_{Xmt}(V)-\V'(X\otimes m,t)\\
	=&\left[\V'\left(X_{Xt}^{V}(T)\otimes m,T\right)+\int_t^T G\left(X_{Xt}^{V}(s),V(s),s\right) ds\right]-\V'(X\otimes m,t) \\
	=&\int_t^T\Big[ G\left(X_{Xt}^{V}(s),V(s),s\right)-G\left(X_{Xt}^{V}(s),\hv\left(X_{Xt}^{V}(s),s;D_X \V'\left(X_{Xt}^{V}(s)\otimes m,s\right)\right),s\right)\\
	&\qquad -\Big(D_V G\left(X_{Xt}^{V}(s),\hv\left(X_{Xt}^{V}(s),s;D_X \V'\left(X_{Xt}^{V}(s)\otimes m,s\right)\right),s\right),\\
	&\qquad\qquad V(s)-\hv\left(X_{Xt}^{V}(s),s;D_X \V'\left(X_{Xt}^{V}(s)\otimes m,s\right)\right)\Big)_{\ur_m}\Big]ds\\
	\geq& \ \lambda\int_t^T\left\| V(s)-\hv\left(X_{Xt}^{V}(s),s;D_X \V'\left(X_{Xt}^{V}(s)\otimes m,s\right)\right) \right\|_{\ur_m}^2 ds.
\end{align*}
Therefore, for any admissible control $V$, we have $\V'(X\otimes m,t)\le J_{Xmt}(V)$. If we set $V=\hv_{Xt}$, we see that $\V'(X\otimes m,t)= J_{Xmt}\left(\hv_{Xt}\right)$. From \eqref{lem8_9}, $\V'$ coincides with the value function $\V$.

\section{Proof of Statements in Section~\ref{sec:Bellman}}\label{app:1}

\subsection{Proof of Theorem~\ref{thm9}}\label{pf:thm9}	

The twice G\^ateaux differentiability of $\mathcal{V}$ in $X$ and Estimates \eqref{lem10_3}-\eqref{lem10_2} are direct consequences of Section~\ref{sec:MFC} and Lemma~\ref{lem:8}. Since $D_X\mathcal{V}(X\otimes m,t)$ exists and satisfies \eqref{lem10_1}, from \cite[Proposition 4.1]{AB} we know that the linear functional derivative of $\mathcal{V}$ exists. From \eqref{lem01_1} and Lemma~\ref{lem:8}, we can obtain \eqref{lem10_4}; 
similarly, from \eqref{lem10_2}, we know that the second order linear functional derivative of $\mathcal{V}$ also exists, and satisfies \eqref{lem10_5}. 

\subsection{Proof of Theorem~\ref{lem10}}\label{pf:lem10}	

The differentiability of $\mathcal{V}$ in $t$ and Estimates \eqref{lem10_8} and \eqref{lem10_7} are direct consequences of Section~\ref{sec:MFC} and Lemma~\ref{lem9}. We now prove \eqref{lem10_6}. From \eqref{lem9_3}, by It\^o's lemma, we have
\begin{equation}\label{lem10_11}
	\begin{split}
		\frac{\dd\mathcal{V}}{\dd t}(X\otimes m,t)=&-G\left(X,\hv(X,t;D_X\mathcal{V}(X\otimes m,t)),t\right)\\
		&-\left(D_X\mathcal{V}(X\otimes m,t),\ G\left(X,\hv(X,t;D_X\mathcal{V}(X\otimes m),t\right)\right)_{\hr_m}\\
		&-\frac{1}{2}\left(D_X^2\V(X\otimes m,t)\left(\sum_{j=1}^n A^j(X,t)\mathcal{N}_t^j\right),\   \sum_{j=1}^n A^j(X,t)\mathcal{N}_t^j\right)_{\hr_m},
	\end{split}		
\end{equation}
where $\mathcal{N}_t^j$'s are $n$ independent standard Gaussian random variables, which are independent of $X$ and the Wiener Process. From \eqref{lem10_5}, we have
\begin{align*}
	&\left(D_X^2\mathcal{V}(X\otimes m,t)\left(\sum_{j=1}^n A^j(X,t)\mathcal{N}_t^j\right),\   \sum_{j=1}^n A^j(X,t)\mathcal{N}_t^j\right)_{\hr_m}\\
	=&\ \e\Bigg[\int_\brn \Bigg(\sum_{j=1}^n \left.A^j(X,t)\right|_x\mathcal{N}_t^j\Bigg)^* D_\xi^2 \frac{d \mathcal{V}}{d\nu}(X\otimes m,t)(X_x)\Bigg(\sum_{j=1}^n (A^j(X,t))_x\mathcal{N}_t^j\Bigg) dm(x)\Bigg]\\
	& +\e\bar{\e}\Bigg[\int_\brn\int_\brn \Bigg(\sum_{j=1}^n \left.A^j\left(\bx,t\right)\right|_z\bar{\mathcal{N}}_t^j\Bigg)^* D_{\xi'}D_\xi \frac{d^2 \mathcal{V}}{d\nu^2}(X\otimes m)\left(X_x,\bx_z\right)  \\
	&\quad\qquad\qquad\qquad \Bigg(\sum_{j=1}^n \left.A^j(X,t)\right|_x\mathcal{N}_t^j\Bigg) dm(z) dm(x) \Bigg]\\
	=&\ \sum_{j=1}^n\e\left[\int_\brn \left(\left.A^j(X,t)\right|_x\right)^* D_\xi^2 \frac{d \mathcal{V}}{d\nu}(X\otimes m,t)(X_x)  \left(\left.A^j(X,t)\right|_x\right) dm(x)\right]\\
	=&\ \sum_{j=1}^n \left(D_\xi^2 \frac{d \mathcal{V}}{d\nu}(X\otimes m,t)(X) A^j(X,t), \  A^j(X,t)\right)_{\hr_m},
\end{align*}	
Substituting the last equation back into \eqref{lem10_11}, from \eqref{lem10_4} and the definition of the lifted Hamiltonian $\bh$ in \eqref{bh''}, we have
\begin{align*}
	&\frac{\dd\mathcal{V}}{\dd t}(X\otimes m,t)\\
	=&-G\left(X,\hv(X,t;D_X\mathcal{V}(X\otimes m,t)),t\right)\\
	&-\left(D_X\mathcal{V}(X\otimes m,t),\ G\left(X,\hv(X,t;D_X\mathcal{V}(X\otimes m)),t\right)\right)_{\hr_m}\\
	& -\sum_{j=1}^n \Big(\frac{1}{2}D_\xi^2 \frac{d \mathcal{V}}{d\nu}(X\otimes m,t)(X) A^j(X,t), \  A^j(X,t)\Big)_{\hr_m}\\
	=&-\bh \left(X,t;D_X\mathcal{V}(X\otimes m,t),\frac{1}{2}D_\xi^2 \frac{d \mathcal{V}}{d\nu}(X\otimes m,t)(X) A(X,t)\right)\\
	=&-\e \bigg[\int_\brn H\bigg(X_x,X\otimes m, t;D_\xi \frac{d\mathcal{V}}{d\nu}(X\otimes m,t)(X_x),\\
	&\qquad\qquad\qquad \frac{1}{2}D_\xi^2 \frac{d \mathcal{V}}{d\nu}(X\otimes m,t)(X_x) \sigma(X_x,X\otimes m,t)\bigg) dm(x)\bigg],
\end{align*}
from which we deduce \eqref{lem10_6}. Then, the solvability of Bellman equation \eqref{Bellman'} is a direct consequence of \eqref{lem10_6} and Theorem~\ref{thm:9}.

\section{Proof of Statements in Section~\ref{sec:Master}}\label{app:2}

\subsection{Proof of Theorem~\ref{thm:10}}\label{pf:thm10}	

The proof of the well-posedness of FBSDEs \eqref{FB:9} and \eqref{FB:10}, and the corresponding convergence results are similar to those of Theorem~\ref{thm3}, and we just omit. We only give the FBSDEs for $\frac{dY_{xmt}}{d\nu}(\xi,s)$, $\frac{dP_{xmt}}{d\nu}(\xi,s)$, $\frac{dQ_{xmt}}{d\nu}(\xi,s)$ and $\frac{dD_xY_{xmt}}{d\nu}(\xi,s)$, $\frac{dD_xP_{xmt}}{d\nu}(\xi,s)$, $\frac{dD_xQ_{xmt}}{d\nu}(\xi,s)$ here. From FBSDEs \eqref{FB:5}, the linear functional derivatives of $\left(Y_{ymt}(s),P_{ymt}(s),Q_{ymt}(s),\hvv_{ymt}(s)\right)$ in $m\in\pr_2(\brn)$ can be characterized as the solution of the following FBSDEs: for $(s,x,\xi)\in[t,T]\times\brn\times\brn$,
\small
\begin{align}
	&\frac{dY_{xmt}}{d\nu}(\xi,s) \notag\\
	=\ &\int_t^s \bigg[f_1(r)\frac{dY_{xmt}}{d\nu}(\xi,r)+f_2(r)\e\left(Y_{\xi mt}(r)+\int_\brn \frac{dY_{ymt}}{d\nu}(\xi,r)dm(y)\right)  +f_3(r)\frac{d\hvv_{xmt}}{d\nu}(\xi,r)\bigg]dr\notag\\
	&+\sum_{j=1}^n\int_t^s \bigg[\sigma^j_1(r)\frac{dY_{xmt}}{d\nu}(\xi,r)+\sigma^j_2(r)\e\left(Y_{\xi mt}(r)+\int_\brn \frac{dY_{ymt}}{d\nu}(\xi,r)dm(y)\right)+\sigma^j_3(r)\frac{d\hvv_{xmt}}{d\nu}(\xi,r)\bigg]dw_j(r);\notag\\
	&\frac{dP_{xmt}}{d\nu}(\xi,s) \notag\\
	=\ &-\sum_{j=1}^n \int_s^T \frac{dQ^j_{xmt}}{d\nu}(\xi,r)dw_j(r)+\left(D_x^2g_T(Y_{xmt}(T),Y_{\cdot mt}(T)\otimes m)\right)^*\frac{dY_{xmt}}{d\nu}(\xi,T)\notag\\
	& +\check{\e}\bigg[\frac{dD_xg_T}{d\nu}(Y_{xmt}(T),Y_{\cdot mt}(T)\otimes m)\left(\cy_{\xi mt}(T)\right) \notag\\
	&\qquad +\int_\brn \left(D_\xi \frac{dD_xg_T}{d\nu}(Y_{xmt}(T),Y_{\cdot mt}(T)\otimes m)\left(\cy_{ymt}(T)\right)\right)^*\frac{d\cy_{ymt}}{d\nu}(\xi,T)dm(y)\bigg] \notag\\
	& +\bar{\e}\left[\int_\brn \left(D_xD_\xi\frac{d g_T}{d\nu}\left(\bar{Y}_{ymt}(T),Y_{\cdot mt}(T)\otimes m\right)(Y_{xmt}(T))\right)^*\frac{d\bar{Y}_{ymt}}{d\nu}(\xi,T) dm(y)\right]  \notag\\
	& +\bar{\e}\bigg\{\int_\brn \check{\e}\bigg[D_\xi\frac{d^2 g_T}{d\nu^2}\left(\bar{Y}_{ymt}(T),Y_{\cdot mt}(T)\otimes m\right)\left(Y_{xmt}(T),\cy_{\xi mt}(T)\right) \notag\\
	&\quad\qquad\qquad +\int_\brn \left( D_{\xi'}D_\xi\frac{d^2 g_T}{d\nu^2}\left(\bar{Y}_{ymt}(T),Y_{\cdot mt}(T)\otimes m\right)\left(Y_{xmt}(T),\cy_{zmt}(T)\right)\right)^* \frac{d\cy_{zmt}}{d\nu}(\xi,T)dm(z)\bigg] dm(y)\bigg\} \notag\\
	& +\bar{\e}\left[\int_\brn \left(D_\xi^2\frac{d g_T}{d\nu}\left(\bar{Y}_{ymt}(T),Y_{\cdot mt}(T)\otimes m\right)(Y_{xmt}(T))\right)^* \frac{dY_{xmt}}{d\nu}(\xi,T) dm(y)\right] \notag\\
	& +\int_s^T\Bigg\{f_1^*(r)\frac{dP_{xmt}}{d\nu}(\xi,r)+\sum_{j=1}^n {\sigma_1^j}^*(r)\frac{dQ^j_{xmt}}{d\nu}(\xi,r) \notag\\
	&\quad\qquad+\left(D_x^2g\left(Y_{xmt}(r),Y_{\cdot mt}(r)\otimes m,\hvv_{xmt}(r),r\right)\right)^*\frac{dY_{xmt}}{d\nu}(\xi,r) \notag\\
	&\quad\qquad +\check{\e}\bigg[\frac{dD_xg}{d\nu}\left(Y_{xmt}(r),Y_{\cdot mt}(r)\otimes m,\hvv_{xmt}(r),r\right)\left(\cy_{\xi mt}(r)\right) \notag\\
	&\quad\qquad\qquad +\int_\brn \left(D_\xi \frac{dD_xg}{d\nu}\left(Y_{xmt}(r),Y_{\cdot mt}(r)\otimes m,\hvv_{xmt}(r),r\right)\left(\cy_{ymt}(r)\right)\right)^*   \frac{d\cy_{ymt}}{d\nu}(\xi,r)dm(y)\bigg] \notag\\
	&\quad\qquad+\left(D_vD_x g\left(Y_{xmt}(r),Y_{\cdot mt}(r)\otimes m,\hvv_{xmt}(r),r\right)\right)^* \frac{d\hvv_{xmt}}{d\nu}(\xi,r) \notag\\
	&\quad\qquad+\bar{\e}\bigg[\int_\brn\bigg( f_2^*(r)\frac{d\bar{P}_{ymt}}{d\nu}(\xi,r)+\sum_{j=1}^n {\sigma^j_2}^*(r)\frac{d\bar{Q}^j_{ymt}}{d\nu}(\xi,r) \bigg) dm(y)\bigg] \notag\\
	&\quad\qquad +\bar{\e}\bigg[\int_\brn \left(D_xD_\xi\frac{d g}{d\nu}\left(\bar{Y}_{ymt}(r),Y_{\cdot mt}(r)\otimes m,\bar{\hvv}_{ymt}(r),r\right)(Y_{xmt}(r))\right)^*  \frac{d\bar{Y}_{ymt}}{d\nu}(\xi,r) dm(y)\bigg] \notag\\
	&\quad\qquad +\bar{\e}\bigg\{\int_\brn \check{\e}\bigg[D_\xi\frac{d^2 g}{d\nu^2}\left(\bar{Y}_{ymt}(r),Y_{\cdot mt}(r)\otimes m,\bar{\hvv}_{ymt}(r),r\right)\left(Y_{xmt}(r),\cy_{\xi mt}(r)\right) \notag\\
	&\qquad\qquad\qquad\qquad +\int_\brn \bigg( D_{\xi'}D_\xi\frac{d^2 g}{d\nu^2}\left(\bar{Y}_{ymt}(r),Y_{\cdot mt}(r)\otimes m,\bar{\hvv}_{ymt}(r),r\right) \notag\\
	&\qquad\qquad\qquad\qquad\qquad\qquad \left(Y_{xmt}(r),\cy_{zmt}(r)\right)\bigg)^* \frac{d\cy_{zmt}}{d\nu}(\xi,r)dm(z)\bigg] dm(y)\bigg\} \notag\\
	&\quad\qquad +\bar{\e}\bigg[\int_\brn \left(D_vD_\xi\frac{d g}{d\nu}\left(\bar{Y}_{ymt}(r),Y_{\cdot mt}(r)\otimes m,\bar{\hvv}_{ymt}(r),r\right)(Y_{xmt}(r))\right)^*  \frac{d\bar{\hvv}_{ymt}}{d\nu}(\xi,r) dm(y)\bigg] \notag\\
	&\quad\qquad +\bar{\e}\bigg[\int_\brn \left(D_\xi^2\frac{d g}{d\nu}\left(\bar{Y}_{ymt}(r),Y_{\cdot mt}(r)\otimes m,\bar{\hvv}_{ymt}(r),r\right)(Y_{xmt}(r))\right)^*  \frac{dY_{xmt}}{d\nu}(\xi,r) dm(y)\bigg]\Bigg\}dr,  \label{FB:9}
\end{align}
\normalsize
where for $(s,x,\xi)\in[t,T]\times\brn\times\brn$,
\begin{align}
	&\frac{d\hvv_{xmt}}{d\nu}(\xi,s) \notag\\
	=\ &\left(D_x\hat{v}(Y_{xmt}(s),Y_{\cdot mt}(s)\otimes m,s;P_{xmt}(s),Q_{xmt}(s))\right)^*\frac{dY_{xmt}}{d\nu}(\xi,s) \notag\\
	&+\check{\e}\bigg[\frac{d\hat{v}}{d\nu}(Y_{xmt}(s),Y_{\cdot mt}(s)\otimes m,s;P_{xmt}(s),Q_{xmt}(s))\left(\cy_{\xi mt}(s)\right) \notag\\
	&\qquad +\int_\brn \left(D_\xi\frac{d\hat{v}}{d\nu}(Y_{xmt}(s),Y_{\cdot mt}(s)\otimes m,s;P_{xmt}(s),Q_{xmt}(s))\left(\cy_{ymt}(s)\right)\right)^* \frac{d\cy_{ymt}}{d\nu}(\xi,s) dm(y)\bigg] \notag\\
	&+\left(D_p\hat{v}(Y_{xmt}(s),Y_{\cdot mt}(s)\otimes m,s;P_{xmt}(s),Q_{xmt}(s))\right)^*\frac{dP_{xmt}}{d\nu}(\xi,s) \notag\\
	&+\sum_{j=1}^n \left(D_{q^j}\hat{v}(Y_{xmt}(s),Y_{\cdot mt}(s)\otimes m,s;P_{xmt}(s),Q_{xmt}(s))\right)^*\frac{dQ^j_{xmt}}{d\nu}(\xi,s), \label{master_5}
\end{align}
and the process $\left(\frac{d\bar{Y}_{ymt}}{d\nu}(\xi,s),\frac{d\bar{P}_{ymt}}{d\nu}(\xi,s),\frac{d\bar{Q}_{ymt}}{d\nu}(\xi,s),\frac{d\bar{\hvv}_{ymt}}{d\nu}(\xi,s) \right)$ is an independent copy of process $\left(\frac{d{Y}_{ymt}}{d\nu}(\xi,s),\frac{d{P}_{ymt}}{d\nu}(\xi,s),\frac{d{Q}_{ymt}}{d\nu}(\xi,s),\frac{d{\hvv}_{ymt}}{d\nu}(\xi,s) \right)$, and $\left(\cy_{y mt}(s), \frac{d\cy_{ymt}}{d\nu}(\xi,s)  \right)$ is another independent copy of $\left({Y}_{y mt}(s), \frac{d{Y}_{ymt}}{d\nu}(\xi,s)  \right)$. From FBSDEs \eqref{FB:8}, the linear functional derivatives of $\left(D_{x}Y_{x mt}(s),D_{x}P_{x mt}(s),D_{x}Q_{x mt}(s)\right)$ in $m\in\pr_2(\brn)$ can be characterized as the solution of the following FBSDEs: for $(s,x,\xi)\in[t,T]\times\brn\times\brn$,
\small
\begin{align}
	&\frac{dD_xY_{xmt}}{d\nu}(\xi,s) \notag\\
	=\ &\int_t^s \left[f_1(r)\frac{dD_xY_{xmt}}{d\nu}(\xi,r)+f_3(r)\frac{dD_x \hvv_{xmt}}{d\nu}(\xi,r)\right]dr \notag\\
	&+\sum_{j=1}^n\int_t^s \left[\sigma^j_1(r)\frac{dD_xY_{xmt}}{d\nu}(\xi,r)+\sigma^j_3(r)\frac{dD_x \hvv_{xmt}}{d\nu}(\xi,r)\right]dw_j(r); \notag\\
	&\frac{dD_xP_{xmt}}{d\nu}(\xi,s) \notag\\
	=\ &-\sum_{j=1}^n \int_s^T \frac{dD_xQ^j_{xmt}}{d\nu}(\xi,r)dw_j(r)+\left(D_x^2g_T(Y_{xmt}(T),Y_{\cdot mt}(T)\otimes m)\right)^* \frac{dD_x Y_{xmt}}{d\nu}(\xi,T) \notag\\
	&+\left(\frac{dY_{xmt}}{d\nu}(\xi,T)\right)^*D_x^3 g_T(Y_{x mt}(T),Y_{\cdot mt}(T)\otimes m)D_{x}Y_{x mt}(T) \notag\\
	&+\check{\e}\bigg[\left(\frac{dD_x^2g_T}{d\nu}(Y_{xmt}(T),Y_{\cdot mt}(T)\otimes m)\left(\cy_{\xi mt}(T)\right)\right)^* \notag\\
	&\qquad +\int_\brn \left(\frac{d\cy_{ymt}}{d\nu}(\xi,T) \right)^* D_\xi\frac{dD_x^2g_T}{d\nu}(Y_{xmt}(T),Y_{\cdot mt}(T)\otimes m) \left(\cy_{ymt}(T)\right) dm(y)\bigg]D_{x}Y_{x mt}(T) \notag\\
	& +\bar{\e}\left[\int_\brn \left(D^2_\xi\frac{d g_T}{d\nu}\left(\bar{Y}_{ymt}(T),Y_{\cdot mt}(T)\otimes m\right)(Y_{xmt}(T))\right)^* dm(y)\right]\frac{dD_x Y_{xmt}}{d\nu}(\xi,T) \notag\\
	& +\bar{\e}\Bigg\{\int_\brn\bigg\{ \left(\frac{d\bar{Y}_{ymt}}{d\nu}(\xi,T)\right)^*D_x D_\xi^2 \frac{dg_T}{d\nu}(\bar{Y}_{y mt}(T),Y_{\cdot mt}(T)\otimes m)(Y_{xmt}(T)) \notag\\
	&\qquad\qquad +\check{\e}\bigg[\left(D_\xi^2\frac{d^2g_T}{d\nu^2}(\bar{Y}_{ymt}(T),Y_{\cdot mt}(T)\otimes m)\left(Y_{xmt}(T),\cy_{\xi mt}(T)\right)\right)^* \notag\\
	&\qquad\qquad\qquad +\int_\brn \left(\frac{d\cy_{zmt}}{d\nu}(\xi,T) \right)^* D_{\xi'}D_\xi^2\frac{d^2g_T}{d\nu^2}(\bar{Y}_{ymt}(T),Y_{\cdot mt}(T)\otimes m)  \left(Y_{xmt}(T),\cy_{zmt}(T)\right) dm(z)\bigg] \notag\\
	&\qquad\qquad +\left(\frac{dY_{xmt}}{d\nu}(\xi,T)\right)^* D_\xi^3 \frac{dg_T}{d\nu}(\bar{Y}_{y mt}(T),Y_{\cdot mt}(T)\otimes m)(Y_{xmt}(T))\bigg\}dm(y)  D_x Y_{xmt}(T)\Bigg\} \notag\\
	& +\int_s^T\Bigg\{f_1^*(r)\frac{dD_xP_{xmt}}{d\nu}(\xi,r)+\sum_{j=1}^n {\sigma_1^j}^*(r)\frac{dD_xQ^j_{xmt}}{d\nu}(\xi,r) \notag\\
	&\quad\qquad+\left(D_x^2 g\left(Y_{xmt}(r),Y_{\cdot mt}(r)\otimes m,\hvv_{xmt}(r),r\right)\right)^* \frac{dD_x Y_{xmt}}{d\nu}(\xi,r) \notag\\
	&\quad\qquad+\left(\frac{dY_{xmt}}{d\nu}(\xi,r)\right)^*D_x^3 g\left(Y_{x mt}(r),Y_{\cdot mt}(r)\otimes m,\hvv_{x mt}(r),r\right)D_x Y_{xmt}(r) \notag\\
	&\quad\qquad+\check{\e}\bigg[\left(\frac{dD_x^2g}{d\nu}\left(Y_{xmt}(r),Y_{\cdot mt}(r)\otimes m,\hvv_{xmt}(r),r\right)\left(\cy_{\xi mt}(r)\right)\right)^* \notag\\
	&\quad\qquad\qquad +\int_\brn \left(\frac{d\cy_{ymt}}{d\nu}(\xi,r) \right)^* D_\xi\frac{dD_x^2g}{d\nu}\left(Y_{xmt}(r),Y_{\cdot mt}(r)\otimes m,\hvv_{xmt}(r),r\right)  \left(\cy_{ymt}(r)\right)  dm(y)\bigg] D_{x}Y_{x mt}(r) \notag\\
	&\quad\qquad+\left(\frac{d\hvv_{xmt}}{d\nu}(\xi,r)\right)^*D_vD_x^2 g\left(Y_{x mt}(r),Y_{\cdot mt}(r)\otimes m,\hvv_{x mt}(r),r\right)D_x Y_{xmt}(r) \notag\\
	&\quad\qquad+\left(D_vD_x g\left(Y_{xmt}(r),Y_{\cdot mt}(r)\otimes m,\hvv_{xmt}(r),r\right)\right)^* \frac{dD_x \hvv_{xmt}}{d\nu}(\xi,r) \notag\\
	&\quad\qquad+\left(\frac{dY_{xmt}}{d\nu}(\xi,r)\right)^*D_vD_x^2 g\left(Y_{x mt}(r),Y_{\cdot mt}(r)\otimes m,\hvv_{x mt}(r),r\right)D_x \hvv_{xmt}(r) \notag\\
	&\quad\qquad+\check{\e}\bigg[\left(\frac{dD_vD_xg}{d\nu}\left(Y_{xmt}(r),Y_{\cdot mt}(r)\otimes m,\hvv_{xmt}(r),r\right)\left(\cy_{\xi mt}(r)\right)\right)^* \notag\\
	&\quad\qquad\qquad +\int_\brn \left(\frac{d\cy_{ymt}}{d\nu}(\xi,r) \right)^* D_\xi\frac{dD_vD_xg}{d\nu}\left(Y_{xmt}(r),Y_{\cdot mt}(r)\otimes m,\hvv_{xmt}(r),r\right) \left(\cy_{ymt}(r)\right) dm(y)\bigg] D_{x}\hvv_{x mt}(r) \notag\\
	&\quad\qquad+\left(\frac{d\hvv_{xmt}}{d\nu}(\xi,r)\right)^*D_v^2D_x g\left(Y_{x mt}(r),Y_{\cdot mt}(r)\otimes m,\hvv_{x mt}(r),r\right)D_x \hvv_{xmt}(r) \notag\\	
	&\quad\qquad+\bar{\e}\bigg[\int_\brn \left(D_\xi^2 \frac{d g}{d\nu}\left(\bar{Y}_{ymt}(r),Y_{\cdot mt}(r)\otimes m,\bar{\hvv}_{ymt}(r),r\right)(Y_{xmt}(r))\right)^*  \frac{dD_xY_{xmt}}{d\nu}(\xi,r)dm(y)\bigg] \notag\\
	&\quad\qquad +\bar{\e}\Bigg\{\int_\brn\bigg\{ \left(\frac{d\bar{Y}_{ymt}}{d\nu}(\xi,r)\right)^*D_x D_\xi^2 \frac{dg}{d\nu}\left(\bar{Y}_{ymt}(r),Y_{\cdot mt}(r)\otimes m,\bar{\hvv}_{ymt}(r),r\right)(Y_{xmt}(r)) \notag\\
	&\quad\qquad\qquad\qquad +\check{\e}\bigg[\left(D_\xi^2\frac{d^2g}{d\nu^2}\left(\bar{Y}_{ymt}(r),Y_{\cdot mt}(r)\otimes m,\bar{\hvv}_{ymt}(r),r\right)\left(Y_{xmt}(r),\cy_{\xi mt}(r)\right)\right)^* \notag\\
	&\quad\qquad\qquad\qquad\qquad +\int_\brn \left(\frac{d\cy_{zmt}}{d\nu}(\xi,r) \right)^* D_{\xi'}D_\xi^2\frac{d^2g}{d\nu^2}\left(\bar{Y}_{ymt}(r),Y_{\cdot mt}(r)\otimes m,\bar{\hvv}_{ymt}(r),r\right)\notag\\
	&\ \qquad\qquad\qquad\qquad\qquad\qquad \left(Y_{xmt}(r),\cy_{zmt}(r)\right) dm(z)\bigg] \notag\\
	&\quad\qquad\qquad\qquad+\left(\frac{d\bar{\hvv}_{ymt}}{d\nu}(\xi,r)\right)^*D_v D_\xi^2 \frac{dg}{d\nu}\left(\bar{Y}_{ymt}(r),Y_{\cdot mt}(r)\otimes m,\bar{\hvv}_{ymt}(r),r\right)(Y_{xmt}(r)) \notag\\
	&\quad\qquad\qquad\qquad +\left(\frac{dY_{xmt}}{d\nu}(\xi,r)\right)^* D_\xi^3 \frac{dg}{d\nu}\left(\bar{Y}_{ymt}(r),Y_{\cdot mt}(r)\otimes m,\bar{\hvv}_{ymt}(r),r\right)(Y_{xmt}(r))\bigg\} \notag\\
	&\quad\qquad\qquad\qquad\qquad dm(y)\Bigg\}D_x Y_{xmt}(r)\Bigg\}dr,  \label{FB:10}
\end{align}
\normalsize
where for $(s,x,\xi)\in[t,T]\times\brn\times\brn$,
\small
\begin{align*}
	&\frac{dD_{x}\hvv_{x mt}}{d\nu}(\xi,s)\\
	=\ &(D_x \hat{v}(Y_{x mt}(s),Y_{\cdot mt}(s)\otimes m,s;P_{x mt}(s),Q_{x mt}(s)))^*\frac{dD_{x}Y_{x mt}}{d\nu}(\xi,s)\\
	&+\left(\frac{dY_{xmt}}{d\nu}(\xi,s)\right)^*D_x^2 \hat{v}(Y_{x mt}(s),Y_{\cdot mt}(s)\otimes m,s;P_{x mt}(s),Q_{x mt}(s))D_{x}Y_{x mt}(s)\\
	&+\check{\e}\bigg[\left(\frac{dD_x\hat{v}}{d\nu}(Y_{xmt}(s),Y_{\cdot mt}(s)\otimes m,s;P_{xmt}(s),Q_{xmt}(s))\left(\cy_{\xi mt}(s)\right)\right)^*\\
	&\qquad +\int_\brn \left(\frac{d\cy_{ymt}}{d\nu}(\xi,s) \right)^* D_\xi\frac{dD_x\hat{v}}{d\nu}(Y_{xmt}(s),Y_{\cdot mt}(s)\otimes m,s;P_{xmt}(s),Q_{xmt}(s)) \\
	&\ \quad\qquad\qquad \left(\cy_{ymt}(s)\right) dm(y)\bigg]D_{x}Y_{x mt}(s)\\
	&+\left(\frac{dP_{xmt}}{d\nu}(\xi,s)\right)^*D_pD_x \hat{v}(Y_{x mt}(s),Y_{\cdot mt}(s)\otimes m,s;P_{x mt}(s),Q_{x mt}(s))D_{x}Y_{x mt}(s)\\
	&+\sum_{j=1}^n\left(\frac{dQ^j_{xmt}}{d\nu}(\xi,s)\right)^*D_{q^j}D_x \hat{v}(Y_{x mt}(s),Y_{\cdot mt}(s)\otimes m,s;P_{x mt}(s),Q_{x mt}(s))D_{x}Y_{x mt}(s)\\
	&+(D_p \hat{v}(Y_{x mt}(s),Y_{\cdot mt}(s)\otimes m,s;P_{x mt}(s),Q_{x mt}(s)))^*\frac{dD_{x}P_{x mt}}{d\nu}(\xi,s)\\
	&+\left(\frac{dY_{xmt}}{d\nu}(\xi,s)\right)^*D_xD_p \hat{v}(Y_{x mt}(s),Y_{\cdot mt}(s)\otimes m,s;P_{x mt}(s),Q_{x mt}(s))D_{x}P_{x mt}(s)\\
	&+\check{\e}\bigg[\left(\frac{dD_p\hat{v}}{d\nu}(Y_{xmt}(s),Y_{\cdot mt}(s)\otimes m,s;P_{xmt}(s),Q_{xmt}(s))\left(\cy_{\xi mt}(s)\right)\right)^*\\
	&\qquad +\int_\brn \left(\frac{d\cy_{ymt}}{d\nu}(\xi,s) \right)^* D_\xi\frac{dD_p\hat{v}}{d\nu}(Y_{xmt}(s),Y_{\cdot mt}(s)\otimes m,s;P_{xmt}(s),Q_{xmt}(s))\\
	&\ \quad\qquad\qquad \left(\cy_{ymt}(s)\right) dm(y)\bigg]D_{x}P_{x mt}(s)\\
	&+\left(\frac{dP_{xmt}}{d\nu}(\xi,s)\right)^*D_p^2 \hat{v}(Y_{x mt}(s),Y_{\cdot mt}(s)\otimes m,s;P_{x mt}(s),Q_{x mt}(s))D_{x}P_{x mt}(s)\\
	&+\sum_{j=1}^n\left(\frac{dQ^j_{xmt}}{d\nu}(\xi,s)\right)^*D_{q^j}D_p \hat{v}(Y_{x mt}(s),Y_{\cdot mt}(s)\otimes m,s;P_{x mt}(s),Q_{x mt}(s))D_{x}P_{x mt}(s)\\
	&+\sum_{i=1}^n\bigg\{(D_{q^i} \hat{v}(Y_{x mt}(s),Y_{\cdot mt}(s)\otimes m,s;P_{x mt}(s),Q_{x mt}(s)))^*\frac{dD_{x}Q^i_{x mt}}{d\nu}(\xi,s)\\
	&\quad\qquad +\left(\frac{dY_{xmt}}{d\nu}(\xi,s)\right)^*D_xD_{q^i} \hat{v}(Y_{x mt}(s),Y_{\cdot mt}(s)\otimes m,s;P_{x mt}(s),Q_{x mt}(s))D_{x}Q^i_{x mt}(s)\\
	&\quad\qquad +\check{\e}\bigg[\left(\frac{dD_{q^i}\hat{v}}{d\nu}(Y_{xmt}(s),Y_{\cdot mt}(s)\otimes m,s;P_{xmt}(s),Q_{xmt}(s))\left(\cy_{\xi mt}(s)\right)\right)^*\\
	&\quad\qquad\qquad +\int_\brn \left(\frac{d\cy_{ymt}}{d\nu}(\xi,s) \right)^* D_\xi\frac{dD_{q^i}\hat{v}}{d\nu}(Y_{xmt}(s),Y_{\cdot mt}(s)\otimes m,s;P_{xmt}(s),Q_{xmt}(s))\\
	&\ \qquad\qquad\qquad\qquad \left(\cy_{ymt}(s)\right) dm(y)\bigg]D_{x}Q^i_{x mt}(s)\\
	&\quad\qquad +\left(\frac{dP_{xmt}}{d\nu}(\xi,s)\right)^*D_pD_{q^i} \hat{v}(Y_{x mt}(s),Y_{\cdot mt}(s)\otimes m,s;P_{x mt}(s),Q_{x mt}(s))D_{x}Q^i_{x mt}(s)\\
	&\quad\qquad +\sum_{j=1}^n\left(\frac{dQ^j_{xmt}}{d\nu}(\xi,s)\right)^*D_{q^j}D_{q^i} \hat{v}(Y_{x mt}(s),Y_{\cdot mt}(s)\otimes m,s;P_{x mt}(s),Q_{x mt}(s)) D_{x}Q^i_{x mt}(s)\bigg\}.
\end{align*}
\normalsize
From the continuity conditions specified in Assumption (B2'), we can deduce the continuity based on the above expressions.

\subsection{Proof of Theorem~\ref{thm:11}}\label{pf:thm11}	

Similar to the proofs of statements in Section~\ref{sec:value_Hilbert}, we only need to prove the case when (B3)-(i) is satisfied. From \eqref{master_1}, \eqref{master_4} and Theorem~\ref{thm:10}, we have \eqref{master_6}. Then, \eqref{thm11_0} is a direct consequence of \eqref{master_6} and Theorem~\ref{thm:11}. We now aim to show that $\mathcal{U}$ is a classical solution of the master equation \eqref{master}. From Theorem~\ref{lem10}, \eqref{master_1} and \eqref{master_4}, we have for $(t,m)\in[0,T)\times\pr_2(\brn)$,
\begin{equation}\label{thm11_1}
	\begin{split}
		-\frac{\dd\mathcal{V}}{\dd t}(m,t)=\ & \int_\brn H\bigg(\xi,m,t; P_{\xi mt}(t),\frac{1}{2}D_x P_{\xi mt}(t) \sigma(\xi,m,t)\bigg) dm(\xi)\\
		=\ &\int_\brn \bigg[P^*_{\xi mt}(t)f(\xi,m,u_{\xi mt}(t),t)+\frac{1}{2}\text{Tr}\left[\sigma\sigma^*(\xi,m,t)D_{x}P_{\xi mt}(t)\right]\\
		&\quad\qquad +g(\xi,m,u_{\xi mt}(t),t) \bigg] dm(\xi).
	\end{split}
\end{equation}
By differentiating both sides of \eqref{thm11_1} with respect to $m$, we can write
\begin{align}
	&-\frac{d}{d\nu}\frac{\dd\V}{\dd t}(m,t)(x) \notag\\
	=\ &H\bigg(x,m,t; P_{x mt}(t),\frac{1}{2}D_x P_{x mt}(t) \sigma(x,m,t)\bigg) \notag\\
	&+\int_\brn\bigg\{ \left(f\left(\xi,m,\hvv_{\xi mt}(t),t\right)\right)^*\frac{dP_{\xi mt}}{d\nu}(x,t)+x^*f^*_2(t)P_{\xi mt}(t) \notag\\
	&\quad\qquad + P^*_{\xi mt}(t)f_3^*(t)\frac{d\hvv_{\xi mt}}{d\nu}(x,t) \notag\\
	&\quad\qquad +\frac{1}{2}\text{Tr}\left[\sigma\sigma^*(\xi,m,t)\frac{dD_{\xi}P_{\xi mt}}{d\nu}(x,t)\right]+\frac{1}{2}\text{Tr}\left[\sigma_2(t)x\sigma^*(\xi,m,t)D_{\xi}P_{\xi mt}(t)\right] \notag\\
	&\quad\qquad +\frac{dg}{d\nu}\left(\xi,m,\hvv_{\xi mt}(t),t\right)(x)+\left(D_vg\left(\xi,m,\hvv_{\xi mt}(t),t\right)\right)^* \frac{d\hvv_{\xi mt}}{d\nu}(x,t) \bigg\}dm(\xi).  \label{thm11_2}
\end{align}
From Assumption (B4) and \eqref{master_7}, we have
\begin{equation}\label{thm11_4}
	\begin{split}
		&f_3(t)P_{\xi mt}(t)+D_vg\left(\xi,m,\hvv_{\xi mt}(t),t\right)=0,
	\end{split}
\end{equation}
and then,
\begin{equation}\label{thm11_3}
	\begin{split}
		&\frac{dH}{d\nu}\left(\xi,m,t;P_{\xi mt}(t),\frac{1}{2}D_\xi P_{\xi mt}(t) \sigma(\xi,m,t)\right)(x)\\
		=\ &x^*f_2(t)^*P_{\xi mt}(t)+\frac{1}{2}\text{Tr}\left[ \sigma_2(t)x\sigma^*(\xi,m,t)D_\xi P_{\xi mt}(t)\right]+\frac{dg}{d\nu}\left(\xi,m,\hvv_{\xi mt}(t),t\right)(x);\\
		&D_p H(\xi,m,t;P_{\xi mt}(t),\frac{1}{2}D_\xi P_{\xi mt}(t) \sigma(\xi,m,t))=\ f\left(\xi,m,\hvv_{\xi mt}(t),t\right).
	\end{split}
\end{equation}
Substituting \eqref{thm11_4} and \eqref{thm11_3} into \eqref{thm11_2}, we deduce that
\begin{align*}
	&-\frac{d}{d\nu}\frac{\dd\V}{\dd t}(m,t)(x)\\
	=\ &H\bigg(x,m,t; P_{x mt}(t),\frac{1}{2}D_x P_{x mt}(t) \sigma(x,m,t)\bigg) \\
	&+\int_\brn\bigg\{ \frac{1}{2}\text{Tr}\left[\sigma\sigma^*(\xi,m,t)\frac{dD_{\xi}P_{\xi mt}}{d\nu}(x,t)\right]\\
	&\qquad\qquad +\left(D_p H\left(\xi,m,t;P_{\xi mt}(t),\frac{1}{2}D_\xi P_{\xi mt}(t) \sigma(\xi,m,t)\right)\right)^*\frac{dP_{\xi mt}}{d\nu}(x,t)\\
	&\qquad\qquad +\frac{dH}{d\nu}\left(\xi,m,t;P_{\xi mt}(t),\frac{1}{2}D_\xi P_{\xi mt}(t) \sigma(\xi,m,t)\right)(x)\bigg\}dm(\xi)
\end{align*}
Substituting \eqref{master_2} and \eqref{master_6} into the last equation, we obtain the master equation \eqref{master}. The proof of the uniqueness is similar to that of the uniqueness result in Theorems~\ref{thm:9}, \ref{lem10} and Proposition 5.4 in \cite{AB5}, and is thus omitted here.	


\end{document}